\documentclass[reqno]{amsart}
 \usepackage{amsbsy,amssymb,amscd,amsfonts,latexsym,amstext,delarray,
 amsmath,graphicx, color, lineno}
 \usepackage[dvips]{epsfig}
\usepackage{hyperref}
\input xypic

\newcommand*\patchAmsMathEnvironmentForLineno[1]{%
  \expandafter\let\csname old#1\expandafter\endcsname\csname #1\endcsname
  \expandafter\let\csname oldend#1\expandafter\endcsname\csname end#1\endcsname
  \renewenvironment{#1}%
     {\linenomath\csname old#1\endcsname}%
     {\csname oldend#1\endcsname\endlinenomath}}%
\newcommand*\patchBothAmsMathEnvironmentsForLineno[1]{%
  \patchAmsMathEnvironmentForLineno{#1}%
  \patchAmsMathEnvironmentForLineno{#1*}}%
\AtBeginDocument{%
\patchBothAmsMathEnvironmentsForLineno{equation}%
\patchBothAmsMathEnvironmentsForLineno{align}%
\patchBothAmsMathEnvironmentsForLineno{flalign}%
\patchBothAmsMathEnvironmentsForLineno{alignat}%
\patchBothAmsMathEnvironmentsForLineno{gather}%
\patchBothAmsMathEnvironmentsForLineno{multline}%
}

\definecolor{Green}{rgb}{0,1,0}
\definecolor{Blue}{RGB}{0,0,191}
\definecolor{mathmodecolor}{RGB}{0,102,0}
\definecolor{keywordcolor}{RGB}{0,51,151}
\definecolor{sourcebackgroundcolor}{RGB}{255,247,223}
\definecolor{unixagred}{RGB}{255,0,0}
\definecolor{lightgray}{RGB}{191,191,191}
\definecolor{green}{RGB}{1,191,191}

\newtheorem{thm}{Theorem}[section]
\newtheorem{prop}[thm]{Proposition}
\newtheorem{cor}[thm]{Corollary}
\newtheorem{lem}[thm]{Lemma}

\newtheorem{defn}[thm]{Definition}
\newtheorem{rem}[thm]{Remark}

\def\qqq{\,,\quad~\forall}

\def\Aut{{\rm Aut}}

\def\End{{\rm End}}

\def\Gal{{\rm Gal}}

\def\Hom{{\rm Hom}}

\def\Ker{{\rm Ker}}

\def\\adsod{{\rm \adsod}}

\def\Spec{{\rm Spec\,}}
\def\Sp{{\rm Spec}\,}

\def\Tr{{\rm Tr}}
\def\tr{{\rm tr}}

\def\A{{\mathbb A}}
\def\B{{\mathbb B}}
\def\C{{\mathbb C}}
\def\F{{\mathbb F}}
\def\K{{\mathbb K}}

\def\N{{\mathbb N}}

\def\Q{{\mathbb Q}}
\def\R{{\mathbb R}}
\def\Z{{\mathbb Z}}
\def\H{{\mathbb H}}
\def\W{{\mathbb W}}

\def\Tr{{\rm Tr}}
\def\tr{{\rm tr}}

\def\cA{{\mathcal A}}

\def\cE{{\mathcal E}}

\def\cH{{\mathcal H}}

\def\cJ{{\mathcal J}}

\def\c\ads{{\mathcal \ads}}

\def\cM{{\mathcal M}}
\def\cO{{\mathcal O}}
\def\cP{{\mathcal P}}

\def\cW{{\mathcal W}}

\newcommand{\ie}{{\it i.e.\/}\ }
\newcommand{\eg}{{\it e.g.\/}\ }
\newcommand{\cf}{{\it cf.\/}\ }

\def\Hom {{\mbox{Hom}}}

\def\End{{\mbox{End}}}
\def\En{{\mbox{\underline{End}}}}

\def\ffp{\mathfrak{p}}

\def\\adso{\mathfrak{\adso}}
\def\An{\mathfrak{ Ring}}

\def\\adsr{\mathfrak{ \adsR}}

\def\Se{\mathfrak{ Sets}}

\def\ads{\H_\K}
\def\kras{\mathbf{K}}

\def\ad{ad\`ele }
\def\te{Teichm\"uller }

\def\cun{\widehat{\Q_p^{\rm ur}}}
\def\zun{\cO_{\widehat{\Q_p^{\rm ur}}}}
\def\zuni{\cO}

\def\qcy{{\Q^{\rm cyc}}}

\def\fr{{\rm Fr}}
\def\sym{{\rm Symm}}
\def\winf{\W_{p^\infty}}
\def\gh{{\rm gh}}
\def\qfr{{\Q^{\rm    cyc, p}_\fr}}
\def\orb{{\mathbb O}}
\def\fpmp{{\cA}}
\def\qsl{{\Q^{\rm    cycl}_\Delta}}
\def\bcsl{{\Z[\Q/\Z]^\Delta}}
\def\qqsl{{\Q[\Q/\Z]^\Delta}}
\def\trr{\Sigma}
\def\mupin{\Q_p/\Z_p}

\def\emup{\mu^{(p)}}
\def\mup{{(\Q/\Z)^{(p)}}}
\def\qcp{\Q^{\rm    cyc, p}}
\def\tilmup{{\cM{(p)}}}

\newcommand{\nil}[1]{}

\title
{On the arithmetic of the BC-system}

\author[Connes]{Alain Connes}
\author[Consani]{Caterina Consani}
\address{A.~Connes: Coll\`ege de France \\
3, rue d'Ulm \\ Paris, F-75005 France
\\ I.H.E.S. and Vanderbilt
University} \email{alain\@@connes.org}
\address{C.~Consani: Mathematics Department \\ Johns Hopkins
University \\ Baltimore, MD 21218 USA} \email{kc\@@math.jhu.edu}
\thanks{The authors are partially supported by the NSF grant
DMS-FRG-0652164. The first author is grateful to Yves Andr\'e for a helpful suggestion. The second author would like to thank the Institut des Hautes \'Etudes Scientifiques and the Coll\`ege de France for some financial support.}


\begin{document}

\maketitle

\begin{abstract} For each prime $p$ and each embedding $\sigma$ of the multiplicative group  of an algebraic closure of $\F_p$ as complex roots of unity, we construct a $p$-adic indecomposable representation $\pi_\sigma$ of the integral BC-system as additive endomorphisms of the big Witt ring of $\bar\F_p$. The obtained representations are the $p$-adic analogues of the complex, extremal KMS$_\infty$ states of the BC-system. The role of the Riemann zeta function, as partition function of the BC-system over $\C$ is replaced, in the $p$-adic case, by the $p$-adic $L$-functions and the polylogarithms whose values at roots of unity encode the KMS states. We use Iwasawa theory to extend the KMS theory to a covering of the completion $\C_p$ of an algebraic closure of $\Q_p$. We show that our previous work on the hyperring structure of the \ad class space, combines with $p$-adic analysis to refine the space of valuations on the cyclotomic extension of $\Q$ as a noncommutative space intimately related to the integral BC-system and whose arithmetic geometry comes close to fulfill the expectations of the ``arithmetic site". Finally, we explain how the integral BC-system  appears naturally also in de Smit and Lenstra construction of the standard model of $\bar\F_p$ which singles out the subsystem associated to the $\hat\Z$-extension of $\Q$.
\end{abstract}

\tableofcontents

\section{Introduction}
This paper describes several arithmetic properties of the BC-system, showing new and interesting connections with the theory of Witt vectors over the algebraic closure of finite fields and with p-adic analysis.

The BC-system is a system of quantum statistical mechanics defined by a noncommutative Hecke algebra  of double classes in $P^+(\Q)$ with respect to the subgroup $P^+(\Z)$, where $P\subset GL_2$ is the ``$ax+b$" algebraic group (\cf \cite{BC,Co-Takagi}). The complex Hecke algebra ${\mathcal H}_\C$ of the system has a highly non-trivial structure  since  its regular representation, in the Hilbert space of one sided classes, generates a factor of type III$_1$ and a canonical ``time evolution" $\sigma_t\in \Aut({\mathcal H}_\C)$. The study of the  KMS-equilibrium states at different temperatures has revealed the arithmetic nature of this dynamical system in view of the following facts \vspace{.05in}

$\bullet$~The partition function of the system is the Riemann zeta function\vspace{.05in}

$\bullet$~There is a phase transition with spontaneous symmetry breaking at the pole of zeta function\vspace{.05in}

$\bullet$~The zero temperature vacuum states implement the global class field isomorphism for $\Q$.\vspace{.05in}

The study of the BC-system inaugurated the interplay between  number-theory and noncommutative geometry. It is exactly the noncommutativity of the  Hecke algebra of the system which generates its non-trivial dynamics. Moreover,  on the noncommutative space  of ad\`eles classes $\A_\Q/\Q^*$, which is naturally associated to the type II dual of the BC-system, one obtains  the spectral realization of zeros of $L$-functions and the trace formula interpretation of the Riemann-Weil explicit formulas  (\cf \cite{Co-zeta}).

Further study (\cf\cite{ccm}) has shown that the integral Hecke algebra ${\mathcal H}_\Z=\Z[\Q/\Z]\rtimes\N$ supplies an integral model to the BC-system. The endomorphisms $\sigma_n(e(r))=e(nr)$, $n\in\N$ act on the canonical generators $e(r)\in \Z[\Q/\Z]$, for $r\in \Q/\Z$ and have natural linear quasi-inverses
  \begin{equation}\label{newrhos0}
\tilde\rho_n: \Z[\Q/\Z] \to \Z[\Q/\Z]\,, \ \ \tilde\rho_n(e(\gamma))=
\sum_{n\gamma'=\gamma}e(\gamma'),
\end{equation}
which are used in the construction of the crossed product and in the presentation of the algebra.

In this paper we establish, for each prime $p$, a strong relation connecting the integral BC-system and the universal Witt ring $\W_0(\bar\F_p)$ of an algebraic closure of a prime field.
The Witt construction is in fact considered in the following three different forms\vspace{.05in}

$-$~ as a $K$-theory endofunctor $A\mapsto \W_0(A)=K_0(\En_A)/K_0(A)$, in the category of commutative rings (with unit)\vspace{.05in}

$-$~ as the big Witt ring $\W(A)$\vspace{.05in}

$-$~ as the functor $A\mapsto\winf(A)$, for $A$ of characteristic $p$.\vspace{.05in}

In the first two cases, the key structures are provided by the following operators\vspace{.05in}

$-$~the \te multiplicative lift  $\tau:A\to \W_0(A)$\vspace{.05in}

$-$~the Frobenius endomorphisms $F_n: \W_0(A)\to\W_0(A)$, $n\in \N$\vspace{.05in}

$-$~the Verschiebung  additive functorial maps $V_n: \W_0(A)\to \W_0(A)$, $n\in \N$\vspace{.05in}

$-$~the $n$-th ghost components $\gh_n: \W_0(A)\to A$, $n\in \N$.\vspace{.05in}

These basic operators extend from the universal ring $\W_0(A)$ to its completion $\W(A)$ whose elements are expressed by Witt vectors, in terms of which all the algebraic operations can be defined in terms of polynomials with integral coefficients. This integrality property encodes a rich and deep arithmetical information. Moreover, the ring structure restricts to divisor stable subsets of $\N$ yielding,  for the set of powers of a prime $p$,  the functor $\winf$.\vspace{.05in}

 In Proposition \ref{bcrelatenew} and Theorem \ref{bcrelate1} we prove that
 the $p$-primary structure of the integral BC-system    is completely encoded by the universal ring $\W_0(\bar\F_p)$,
with a precise dictionary expressing the key operators $\sigma_n$ and $\tilde\rho_n$ of the  BC-system as respectively the Frobenius $F_n$ and Verschiebung $V_n$ on $\W_0(\bar\F_p)$. The isomorphism connecting these algebraic structures depends upon the choice of a group isomorphism of the multiplicative group of
 $\bar\F_p$ with the group  of {\em complex} roots of unity of order prime to $p$: the ambiguity inherent to this choice is the same as that pertaining to the construction of Brauer lift of characters.

The completion process associated to the inclusion $\W_0(A)\subset \W(A)$ with dense image,
 is then used in Theorem \ref{bcrelate2}  to obtain,  when $A=\bar\F_p$ and  for each injective group homomorphism $\sigma: \bar\F_p^\times \to \C^\times$, a $p$-adic {\em indecomposable} representation $\pi_\sigma$ of the integral BC-system as additive endomorphisms of the big Witt ring $\W(\bar\F_p)$. The construction  uses the identification proven in Theorem \ref{bcrelate1} and the implementation of the Artin--Hasse exponentials. These representations are the $p$-adic analogues of the complex, extremal KMS$_\infty$ states of the BC-system.  In Section \ref{sectL} this analogy  is pursued much further. By implementing the theory of $p$-adic $L$-functions, we construct an analogue, in the $p$-adic case, of the partition function and of the KMS$_\beta$ states.  In particular, we show that the division relations for the $p$-adic polylogarithms at roots of unity correspond to the KMS condition. In \S \ref{sectexten} we prove that the definition of the functionals satisfying such condition extends from the standard ``extended s-disk'' to the natural multiplicative group covering of $\C_p$. These results are the $p$-adic counterparts of the statements proven in \cite{cmqsm} for function fields.  However, we also recognize an important difference with respect to the complex case, namely the presence of an added symmetry  at non-zero temperature, due to the invariance of the states under the natural involution of $\qcy$ which replaces each root of unity by its inverse. This added symmetry is a consequence of the vanishing of the $p$-adic $L$-functions associated to odd Dirichlet characters.

For $p$ a prime number, the set $X_p$  of all injective group
homomorphisms $\sigma:\bar\F_p^\times \to \Q/\Z$ is the parameter space  for the $p$-adic representations of the integral BC-system.
In Section \ref{sectval}, we relate this set with the space ${\rm Val}_p(\qcy)$ of extensions of the $p$-adic valuation to the maximal abelian field extension $\qcy$ of $\Q$. We view  $\qcy$ as an abstract field defined as the quotient of the group ring $\Q[\Q/\Z]$ by the cyclotomic ideal (\cf Definition \ref{acf}).
Let $\mup$ be the subgroup of $\Q/\Z$ of fractions with denominator prime to $p$ and let $\qcp$ be the subfield (\ie the inertia subfield) of $\qcy$ generated by the group $\emup\sim \mup$ of roots of unity of order prime to $p$. We describe canonical isomorphisms of ${\rm Val}_p(\qcy)$ with each of the following spaces\vspace{.05in}

$(1)$~The space of sequences  of irreducible polynomials $P_n(T)\in\F_p[T]$, $n\in \N$,  fulfilling the basic conditions of the Conway polynomials (\cf Theorem \ref{conway0}).\vspace{.05in}

    $(2)$~The space $\Sigma_p$ of bijections of the monoid $\tilmup=\emup\cup \{0\}$ commuting with their conjugates, as in Definition \ref{defnsp} (\cf Proposition \ref{hope}).\vspace{.05in}

    $(3)$~The space $\Hom(\qfr, \Q_p)$ of field homomorphisms,    where $\qfr\subset \qcp$ is the decomposition subfield, \ie the fixed field under the Frobenius automorphism (\cf Proposition \ref{homfield}).\vspace{.05in}

 $(4)$~The quotient of the space $X_p$  by the action of $\Gal(\bar\F_p)$ (\cf Proposition \ref{try}).\vspace{.05in}

 $(5)$~The algebraic spectrum of the quotient algebra $\F_p[\mup]/J_p$, where $J_p$ is the reduction modulo $p$ of the cyclotomic ideal (\cf Definition \ref{acf} and Proposition \ref{cpcp1}).\vspace{.05in}

For a global field $\K$ of positive characteristic (\ie a function field  associated to a projective, non-singular curve  over a finite field $\F_q$) it is a well known fact that the space of valuations of the maximal abelian extension
$\K^{\rm ab}$ of $\K$  has a geometric meaning. In fact, for each finite extension
$E$ of $\bar\F_q\otimes_{\F_q}\K\subset \K^{\rm ab}$, the space
${\rm Val}(E)$ of (discrete) valuations of $E$ is  an algebraic, one-dimensional scheme
  whose non-empty open sets are the complements of the finite subsets $F\subset {\rm Val}(E)$. The structure sheaf  is locally defined by the intersection $\bigcap_F R$ of the valuation rings inside $E$. Then the space ${\rm Val}(\K^{\rm ab})$ is the projective limit of the schemes ${\rm Val}(E)$, $E\subset \K^{\rm ab}$.

For the global field $\K=\Q$ of rational numbers, one can consider its maximal abelian extension  $\qcy$ as an abstract  field  and try to follow a similar idea. In Section \ref{curve}, we show however that the space ${\rm Val}(\qcy)$ provides only a rough analogue,  in characteristic zero, of  ${\rm Val}(\K^{\rm ab})$. Our approach to this problem is guided and motivated by the following three results contained in our previous  work\vspace{.05in}

$(a)$~The adelic interpretation of the loop groupoid $\Pi_1^{\rm ab}(X)'$ of the abelian cover of the algebraic curve $X$ associated to a function field (\cf \cite{wagner} and \S~\ref{subsectprel})\vspace{.05in}

 $(b)$~The determination  of the counting function $N(q)$ (a distribution on $[1,\infty)$) which replaces, for $\K=\Q$, the classical Weil counting function for a function field (\cf  \cite{announc3} and \S~\ref{subsectwhync})\vspace{.05in}

 $(c)$~The interpretation of the counting function  $N(q)$ as an intersection number, using the action of the id\`ele class group on the \ad class space (\cf \cite{japanese}).
 \vspace{.05in}

By applying these results we find that the sought for geometric fiber over a non-archimedean, rational prime $p$ is  the total space of a principal bundle, with base ${\rm Val}_p(\qcy)$ and structure group given by a connected, compact solenoid $S$ whose definition is given in Proposition \ref{solenoid}. Then,  in Proposition \ref{fibration} we derive a natural construction for the fiber as the mapping torus $Y_p$ of the action of the Frobenius on the space $X_p$.  In Section \ref{subsectarch}, we consider  the fiber $Y_\infty$  over the archimedean prime, with the implementation of the theory of multiplicative norms.

The interpretation given in $(c)$ for the counting function as intersection number shows that the fibers $Y_p$ should not be considered in isolation, but as being part of an ambient noncommutative space which is responsible for the transversality factors due to the archimedean contribution to the explicit formulas. This interpretation is explained in details in Section \ref{subsectwhync}.

In Section \ref{subsectendo}, we show that the integral BC-system gives, for each $p$ (including the archimedean prime), a natural embedding of the fiber $Y_p$  into a noncommutative space constructed using the set $\cE(\C_p)$ of the $\C_p$-rational points  of the affine group scheme $\cE$ which defines the abelian part of the system  (\cf\cite{ccm}). Here $\C_p$ denotes the $p$-adic completion of an algebraic closure of $\Q_p$. This result shows that the space
\begin{equation}\label{algadclass21}
 X(\C_p):=  \left( \cE(\C_p)\times (0,\infty)\right)/(\N\times \{\pm 1\})
\end{equation}
 matches,  for any rational prime $p$ including $p=\infty$,  the definition of the \ad class space. In Proposition \ref{groupscheme} we show, using the fact that $\cE$ is a group scheme, that $X(\C_p)$ is a free module of rank one over the  hyperring  $\H_\Q$ of the ad\`ele classes.  The problem of a correct interpretation of the connected factor $(0,\infty)$ in \eqref{algadclass21} remains open.

It is a general principle that in our constructions the noncommutative spaces arise as $X(A)$ for a commutative ring $A$ (\cf\eqref{algadclass21}), while  the classical subspaces of $X(A)$ are defined as the support of the cyclotomic ideal (in the affine scheme $\cE = \Sp(\Z[\Q/\Z])$).

We end the paper by showing in Section \ref{sectlens} the relevance of the  recent work of B. de Smit and H. Lenstra (\cf\cite{smitlenstra}) on the ``standard model" for the algebraic closure of a finite field. When $\K$ is a function field, the intermediate extension $\K\subset L=\bar\F_q\otimes_{\F_q}\K\subset \K^{\rm ab}$ plays an important geometric role, namely the extension of scalars to an algebraically closed field, for the algebraic curve associated to $\K$.
When $\K=\Q$, we show that the intermediate extension $\Q\subset \qsl\subset \qcy$ used by de Smit and Lenstra, comes very close to fulfill the expected properties for a similar intermediate extension $\Q\subset L\subset \qcy$. Their construction provides a conceptual construction of the subfield of $\bar\F_p$ union of all extensions whose degree is prime to $p$.  In the very last part of the paper  we recall one of the first applications provided by E. Witt of his functor, which is a conceptual construction of the missing piece $\bigcup_{n} \F_{p^{p^n}}\subset \bar\F_p$,
using the simple equation $X^p=X+1$ in Witt vectors.

\section{The functor $\W_0$}\label{sectwzero}

In this section we recall the definition and the main properties of the universal ring $\W_0(A)$,  where $A$ is any commutative ring with unit. We refer to \cite{Al} to read more details. The second part of the section describes $\W_0(k)$, for an algebraically closed field $k$.

One lets $\En_A$ (or $\text{End}~\cP(A)$) be the category of endomorphisms of  projective $A$-modules of finite rank. The objects are  pairs $(E,f)$ where $E$ is a finite, projective $A$-module and $f\in \End_A(E)$. The morphisms in this category are required to commute with the endomorphisms $f$.  The following operations of direct sum and tensor product
\begin{equation}\label{operations}
    (E_1,f_1)\oplus (E_2,f_2)=(E_1\oplus E_2,f_1\oplus f_2)\,, \ \ (E_1,f_1)\otimes (E_2,f_2)=(E_1\otimes E_2,f_1\otimes f_2)
\end{equation}
turn the Grothendieck group $K_0(\En_A)$ into  a (commutative) ring. The pairs of the form $(E,f=0)$ generate the ideal $K_0(A)\subset K_0(\En_A)$. We denote the quotient ring by  $\W_0(A)$
\begin{equation}\label{defnwo}
    \W_0(A)=K_0(\En_A)/K_0(A).
\end{equation}

By construction, $\W_0$ is an endofunctor of the category $\An$ of commutative rings with unit. Several key operators and maps act on $\W_0$, the following are the most relevant ones for our applications\vspace{.05in}

$(1)$~ The \te lift $\tau:A\to \W_0(A)$ which is a multiplicative map.\vspace{.05in}

$(2)$~ For $n\in \N$, the Frobenius ring endomorphisms $F_n: \W_0(A)\to \W_0(A)$.\vspace{.05in}

$(3)$~For $n\in \N$, the Verschiebung (shift) additive functorial maps $V_n: \W_0(A)\to \W_0(A)$.\vspace{.05in}

$(4)$~For $n\in \N$, the $n$-th ghost component homomorphisms $\gh_n: \W_0(A)\to A$.\vspace{.05in}

We shortly recall their definitions.\vspace{.05in}

$(1)$~The \te lift $\tau=[\cdot]: A \to \W_0(A)$ is  defined as $f\mapsto \tau(f)=[f]=(A,f)$.\vspace{.05in}

$(2)$~For $n\in\N$, the operations in $\En_A$ of raising an endomorphism $f$ to the $n$-th power induce the Frobenius ring endomorphims in $\W_0(A)$
\begin{equation}\label{frob1}
F_n: \W_0(A)\to \W_0(A),\qquad    F_n(E,f)=(E,f^n).
\end{equation}
$(3)$~For $n\in\N$, the Verschiebung  maps are defined by the following operations on matrices
\begin{equation}\label{vermat}
V_n: \W_0(A)\to \W_0(A), \qquad V_n(E,f) = (E^{\oplus n},\begin{bmatrix}0 & 0&\cdots&\cdots&f\\
1&0&0&\cdots&0\\
0&1&0&\cdots&0\\
\cdots&\cdots&\cdots&\cdots&\cdots\\
0&0&0&1&0
\end{bmatrix}).
\end{equation}
$(4)$~For $n\in\N$,  the ghost components are given by
\begin{equation}\label{ghost1}
\gh_n: \W_0(A)\to A,\quad   \gh_n(E,f)={\rm Trace}(f^n).
\end{equation}\vspace{.05in}

Let $\Lambda(A):=1+tA[[t]]$ be the multiplicative abelian group of formal power series with constant term $1$.
The (inverse of the) characteristic polynomial defines a   homomorphism of abelian groups
\begin{equation}\label{charpol}
 L:\W_0(A)\to \Lambda(A),\quad   L(E,f)=\det(1-tM(f))^{-1}
\end{equation}
where $M(f)= (a_{ij})$ is the matrix associated to $f: E\to E$ (\ie $f~\leftrightarrow~\sum_i x_i^*\otimes x_i$, $x_i^*\in E^*, x_i\in E$, $a_{ij} = \langle x_i^*,x_j\rangle$). By a fundamental result of G. Almkvist (\cite{Al} Theorem 6.4 or \cite{Al1} Main Theorem), one has
 \begin{thm}\label{almthm}
The map $L$ is  injective and its image is the subgroup of $\Lambda(A)$
\begin{equation}\label{rangeL}
{\rm Range}(L)=\{(1+a_1t+\ldots +a_nt^n)/(1+b_1t+\ldots +b_nt^n)\mid a_j, b_j \in A\}.
\end{equation}
\end{thm}
Note in particular that for $E$  a finite, projective $A$-module and $f,g\in \End_A(E)$ one has
\begin{equation}\label{fggf}
    (E,fg)=(E,gf)\in \W_0(A).
\end{equation}
One also has
\begin{equation}\label{rabi9}
    V_{nm}=V_n\circ V_m=V_m\circ V_n\,, \ \ F_{nm}=F_n\circ F_m=F_m\circ F_n.
\end{equation}
The following proposition collects together several standard equations connecting these operators

\begin{prop} \label{rabi10}  Let $A$ be a commutative ring and $x,y\in \W_0(A)$. The following hold\vspace{.05in}

$(1)$~$F_n\circ V_n(x)=nx$.\vspace{.05in}

$(2)$~$V_n(F_n(x)y)=xV_n(y)$.\vspace{.05in}

$(3)$~If $(m,n)=1$, ~~$V_m\circ F_n=F_n\circ V_m$.\vspace{.05in}

$(4)$~For $n\in \N$, ~~ $V_n(x)V_n(y)=nV_n(x y)$.\vspace{.05in}

$(5)$~For $n\in \N$, ~~$F_n(\tau(f))=\tau(f^n)$.\vspace{.05in}

$(6)$~For $n,m\in \N$, ~~$\gh_n(F_m(f))=\gh_{nm}(f)$.\vspace{.05in}

$(7)$~ $\gh_n(V_m(f))=\begin{cases} m \,\gh_{n/m}(f) & \text{if $m\vert n$}\\ 0 &\text{otherwise.}\end{cases}$\vspace{.05in}
\end{prop}
\proof All proofs are straightforward, we just check $(4)$ as an example. For $x\in\End_A(E)$, the action of $X=V_n(E,x)$ on vectors $\xi=(\xi_1,\ldots,\xi_n)\in E^{\oplus n}$ is given by
\begin{equation*}
    (X\xi)_1=x\xi_n\,, \ \ (X\xi)_j=\xi_{j-1}\qqq j,~ 2\leq j\leq n.
\end{equation*}
Similar formulas hold for $Y=V_n(F,y)$, for $y\in\End_A(F)$. By definition, $V_n(x)V_n(y)$ corresponds to $X\otimes Y\in\End_A(E^{\oplus n}\otimes F^{\oplus n})$. This endomorphism decomposes as the direct sum of $n$ endomorphisms of $(E\otimes F)^{\oplus n}$, each of these is of the form
\begin{equation*}
    \begin{bmatrix}0 & 0&\cdots&\cdots&x\otimes 1\\
1&0&0&\cdots&0\\
0&1\otimes y&0&\cdots&0\\
\cdots&\cdots&\cdots&\cdots&\cdots\\
0&0&0&1&0
\end{bmatrix},\quad\text{or}\quad  \begin{bmatrix}0 & 0&\cdots&\cdots&x\otimes y\\
1&0&0&\cdots&0\\
0&1&0&\cdots&0\\
\cdots&\cdots&\cdots&\cdots&\cdots\\
0&0&0&1&0
\end{bmatrix}.
\end{equation*}
By applying \eqref{fggf}, one checks that each of the above endomorphisms is equivalent to $V_n(x\otimes y)$. The equality  $V_n(x)V_n(y)=nV_n(x y)$ follows.\endproof

We shall apply the following proposition to the case $A=k=\bar\F_p$ an algebraic closure of $\F_p$.
\begin{prop}\label{wofield} Let $k$ be an algebraically closed field. Then the map which associates to $(E,f)\in \En_k$ the divisor $\delta(f)$ of non-zero eigenvalues of $f$ (with multiplicity taken into account)  extends to a ring isomorphism
\begin{equation}\label{eigeniso}
 \delta: \W_0(k)\stackrel{\sim}{\to} \Z[k^\times].
\end{equation}
Under the above isomorphism, the Frobenius $F_n$ on $\W_0(k)$ is given on $\Z[k^\times]$ by
the natural linearization of the group endomorphism $k^\times\to k^\times,~g\mapsto g^n$.
\end{prop}

\proof By applying Theorem \ref{almthm}, the characteristic polynomial  extends to a complete invariant on $K_0(\En_k)$ and to an isomorphism of $K_0(\En_k)$ with the ring of quotients of monic polynomials in $k[t]$. Moding out this ring by $K_0(k)$ means that one removes the powers of the variable. Thus the divisor of non-zero eigenvalues of $f$ extends to define a bijection of sets $\W_0(k)\simeq \Z[k^\times]$.

It remains to check that this bijection preserves the ring operations. For addition, the set underlying the divisor $\delta(f_1+f_2)$ is the disjoint union of the two sets of roots of $f_j$  and hence $\delta(f_1+f_2) = \delta(f_1)+\delta(f_2)$. For the product, it is enough and easy to check that the tensor product of two rank one elements $(k,a)\otimes (k,b)$ is given by $(k,ab)$ for non-zero elements of $k$. The statement about $F_n$ is checked in the same way using \eqref{frob1} on elements $(k,a)$.
\endproof
We recall the following formula for $L(f)$ in terms of the divisor $\delta(f)=\sum n(\alpha)[\alpha]\in \Z[k^\times]$
\begin{equation}\label{Ldelta}
    L(f)=\prod (1-\alpha \,t)^{-n(\alpha)}.
\end{equation}

\begin{cor}\label{noncan}
For any given  isomorphism
$
   \sigma: \bar \F_p^\times\stackrel{\sim}{\longrightarrow}\mup
$
of the multiplicative group of the algebraic closure $\bar \F_p$ with the subgroup
 $\mup\subset\Q/\Z$ of fractions with denominator prime to $p$, one derives an isomorphism
\begin{equation}\label{iso1}
    \tilde\sigma: \W_0(\bar \F_p)\stackrel{\sim}{\longrightarrow}\Z[\mup].
\end{equation}
Under the isomorphism $\tilde\sigma$, the Frobenius $F_n$ of $\W_0(\bar \F_p)$ is given on $\Z[\mup]$ by
the natural linearization of the group endomorphism $\mup\to\mup$,  $g\mapsto g^n$ (\ie $\gamma\mapsto n\gamma$ in additive notation).
\end{cor}

\section{The integral BC-system}\label{sectbcsys}

 For each $n\in\N$, one defines group ring endomorphisms
\[
\sigma_n: \Z[\Q/\Z] \to \Z[\Q/\Z],\qquad \sigma_n(e(\gamma))=e(n\gamma)
\]
   and the following additive maps
  \begin{equation}\label{newrhos}
\tilde\rho_n: \Z[\Q/\Z] \to \Z[\Q/\Z], \qquad \tilde\rho_n(e(\gamma))=
\sum_{n\gamma'=\gamma}e(\gamma').
\end{equation}
We recall from \cite{ccm}, Proposition 4.4, the following result
\begin{prop}\label{mapstilderho} The endomorphisms $\sigma_n$ and the maps $\tilde\rho_m$  fulfill the following relations
\begin{equation}\label{rhomult}
 \sigma_{nm}=\sigma_{n}\sigma_{m}\,,\ \
    \tilde\rho_{mn}=\tilde\rho_{m}\tilde\rho_{n}\qqq m,n\in\N
\end{equation}
\begin{equation}\label{rhosigmamult}
    \tilde\rho_{m}(\sigma_m(x)y)=x\tilde\rho_{m}(y)\qqq x,y\in \Z[\Q/\Z]
\end{equation}
\begin{equation}\label{divisible}
    \sigma_c(\tilde\rho_b(x))=(b,c)\,\tilde\rho_{b'}(\sigma_{c'}(x))\,, \ \
    b'=b/(b,c)\,, \ \ c'=c/(b,c)\,,
\end{equation}
where $(b,c)=\text{gcd}(b,c)$.
\end{prop}
Note that taking $b=c=n$ in \eqref{divisible} gives
\begin{equation}\label{divisiblebis}
    \sigma_n(\tilde\rho_n(x))=n\,x\qqq x\in \Z[\Q/\Z].
\end{equation}
On the contrary, if we take $b=n$ and $c=m$ to be relatively prime we get
\begin{equation}\label{divisibleter}
    \sigma_n\circ\tilde\rho_m= \tilde\rho_m\circ \sigma_n.
\end{equation}
We recall from \cite{ccm} (Definition 4.7 and \S 4.2) the following facts. The integral $BC$-algebra
 is the algebra $\cH_\Z=\Z[\Q/\Z]\rtimes_{\tilde\rho}\N$ generated by the group
ring $\Z[\Q/\Z]$, and by the elements $\tilde\mu_n$ and $\mu_n^*$, with
$n\in\N$,
 which satisfy
the relations
\begin{equation}\label{presoverZ1}
\begin{array}{l}
\tilde\mu_n x \mu_n^* = \tilde\rho_n(x)\ \ \ \  \\[3mm]
\mu_n^* x = \sigma_n(x) \mu_n^*  \\[3mm]
x \tilde\mu_n = \tilde\mu_n \sigma_n(x),
\end{array}
\end{equation}
where $\tilde\rho_m$, $m\in\N$ is defined in \eqref{newrhos}, as well as the
relations
\begin{equation}\label{presoverZ2}
\begin{array}{l}
\tilde\mu_{nm}= \tilde\mu_n  \tilde\mu_m \qqq n,m\in\N\\[3mm]
\mu_{nm}^* =\mu_{n}^*\mu_{m}^* \qqq n,m\\[3mm]
\mu_n^* \tilde\mu_n  =n \\[3mm]
\tilde\mu_n\mu_m^* = \mu_m^*\tilde\mu_n \ \ \ \ (n,m) = 1.
\end{array}
\end{equation}
After tensoring by $\Q$, the Hecke
algebra $\cH_\Q=\cH_\Z\otimes_\Z\Q$ has a simpler
explicit presentation with generators $\mu_n(=\frac 1n\tilde\mu_n)$, $\mu_n^*$, $n\in \N$ and $e(r)$,
for $r\in \Q/\Z$, satisfying the relations\vspace{.05in}

$\bullet$~$\mu_n^*\mu_n =1$, ~$\forall n\in \N$,\vspace{.05in}

$\bullet$~$ \mu_m \mu_n=\mu_{mn} $,  ~$\mu^*_m \mu^*_n=\mu^*_{mn} $, $\forall m,n\in \N$,\vspace{.05in}

$\bullet$~$\mu_n\mu_m^* = \mu_m^*\mu_n$,  if $(n,m) = 1$,\vspace{.05in}

$\bullet$~$e(0)=1$, $e(r)^*=e(-r)$, and $e(r)e(s)=e(r+s)$, ~ $\forall r,s\in
\Q/\Z$,\vspace{.05in}

$\bullet$~For all $n\in \N$ and all $r\in \Q/\Z$
\begin{equation}\label{scaling}
 \mu_n\, e(r)\, \mu_n^* = \frac{1}{n} \sum_{ns=r} e(s).
\end{equation}

After tensoring by $\C$ and completion one gets a $C^*$-algebra with a natural time evolution $\sigma_t$ (\cite{BC}, \cite{CMbook} Chapter III). The extremal KMS states below critical temperature vanish on the monomials $\mu_n x\mu_m^*$ for $n\neq m$ and $x\in \Q[\Q/\Z]$ and their value on $\Q[\Q/\Z]$ is given by
\begin{equation}\label{BC-KMS1infty}
\varphi_{\beta,\rho}(e(a/b)) =
\frac{1}{\zeta(\beta)}\sum_{n=1}^{\infty}n^{-\beta} \rho
(\zeta_{a/b}^n),
\end{equation}
where  $\rho\in \hat\Z^*$ determines an embedding in $\C$ of the
 cyclotomic field $\qcy$
generated by the abstract roots of  unity.

\section{$\W_0(\bar\F_p)$ and the BC-system}\label{sectwzebc}

In \cite{Quillen} Quillen makes use of the choice of an embedding
 \begin{equation}\label{embedd}
    \sigma:\bar\F_p^\times \to \C^\times
 \end{equation}
in the study of the algebraic K-theory of the general linear group over a finite field. In this section we compare the description of the universal Witt ring $\W_0(\bar \F_p)$, endowed with the structure given by the Frobenius endomorphisms $F_n$ and the Verschiebung  maps $V_n$ with the integral BC-algebra $\cH_\Z$.

By a simple comparison process we notice that the relations \eqref{rhomult}, \eqref{rhosigmamult}, \eqref{divisible} holding on $\cH_\Z$ are the same as those fulfilled by the Frobenius endomorphisms $F_n$ and the Verschiebung maps $V_n$ on $\W_0(\bar\F_p)$. More precisely, under the correspondences $\sigma_n\to F_n$,
$\tilde\rho_n\to V_n$ the two relations of \eqref{rabi9} correspond to \eqref{rhomult}, and the first three relations of Proposition \ref{rabi10} correspond respectively to \eqref{divisiblebis}, \eqref{rhosigmamult} and \eqref{divisibleter}. These results evidently point out to the existence of a strong relation between the ($\lambda$)-ring $\W_0(\bar \F_p)$ and the group ring $\Z[\Q/\Z]$ endowed with the aforementioned operators.

Next, we compare the two groups rings: $\Z[\mup]$ and $\Z[\Q/\Z]$ which arise in the description of $\W_0(\bar\F_p)$ and in the construction of the BC-algebra respectively.  One has a surjective group homomorphism: $\Q/\Z\to\mup$ induced by the canonical factorization of the groups
\begin{equation}\label{factorg}
    \Q/\Z=\mup\times \mupin
\end{equation}
where $\mupin$ is the group of fractions whose denominator is a power of $p$. Thus one obtains a corresponding factorization of the rings
\begin{equation}\label{factoralg}
    \Z[\Q/\Z]=\Z[\mup]\otimes_\Z \Z[\mupin].
\end{equation}
By using the trivial representation of $\mupin$ (\ie the augmentation $\epsilon$ of $\Z[\mupin]$), one gets a retraction $r=id \otimes \epsilon$ producing the splitting
\begin{equation}\label{retract}
    \Z[\mup]\stackrel{j_p}{\longrightarrow} \Z[\Q/\Z]\stackrel{id \otimes \epsilon}{\longrightarrow}\Z[\mup].
\end{equation}
Notice  that  $\mup$ is preserved  by the action of the map $\gamma\mapsto n\gamma$,  $\gamma\in\mup$). This implies that the endomorphisms $\sigma_n$ acting on the BC-algebra restrict naturally to determine endomorphisms $\sigma_n: \Z[\mup]\to \Z[\mup]$.\vspace{.05in}

Let us denote by $I(p)\subset\N$ the set of integers which are prime to $p$. The following lemma describes the projection of the operators $\tilde\rho_n$ of the BC-algebra on the group ring $\Z[\mup]$
\begin{prop}\label{rrhon} Let $n=p^km$, with $m\in I(p)$. For $\gamma \in \Q/\Z$, we write modulo $1$
\begin{equation}\label{decdec}
\gamma=\frac a b+ \frac {c}{p^{s}}\,, \ b\in I(p)\,, \ a,b,c ,s\in \N.
\end{equation}
Then, with $\tilde\rho_n$ as in \eqref{newrhos} we have
\begin{equation}\label{comrhon}
   r\circ \tilde\rho_n(e(\gamma))=p^k\sum_{w=0}^{m-1} e(\frac{f + wb}{bm})
\end{equation}
where $y=\frac{f}{bm}$, $f\in \Z/bm\Z$, is the unique solution in $\Q/\Z$, with denominator prime to $p$  of the equation
\begin{equation}\label{yy}
p^ky=\frac{a}{bm}\in \Q/\Z.
\end{equation}
\end{prop}
\proof The existence and uniqueness of the decomposition \eqref{decdec} derives from the factorization \eqref{factorg}. For $d\in I(p)$, the endomorphism of $\Q/\Z$: $x\mapsto px$ restricts to an automorphism on the subgroup $G_d=\{\frac ad\in\Q/\Z|\, a\in \Z\}\subset\Q/\Z$. For $d=bm$, this fact shows the existence and uniqueness of the solution $y=\frac{f}{bm}$ of \eqref{yy}.  One has $p^ky=\frac{a}{bm}+ j$ for some integer $j\in \Z$, thus
\begin{equation*}
y=\frac{a}{bmp^k}+\frac{j}{p^k}=\frac{a}{bn}+\frac{j}{p^k},\ \ ny=\frac ab +jm.
\end{equation*}
By applying \eqref{factorg}, one also has a decomposition of the form
\begin{equation}\label{cps}
\frac{c}{np^s}=\frac{c}{mp^{s+k}}=\frac{d}{m}+\frac{e}{p^{s+k}}.
\end{equation}
One has $ny=\frac ab$ modulo $1$, $n \frac{c}{np^s}=\frac {c}{p^{s}}$,  thus
the solutions of the equation $n\gamma'=\gamma$ in $\Q/\Z$ which enter in \eqref{newrhos} are of the form
\begin{equation*}
\gamma'=y+\frac{c}{np^s}+\frac u m + \frac{v}{p^k}\,, \ \ u\in\{0,\ldots,m-1\}\,, \ v\in \{0,\ldots,p^k-1\}.
\end{equation*}
By using \eqref{cps} one derives
\begin{equation*}
\gamma'=y+\frac u m +\frac{d}{m}+ \frac{v}{p^k}+\frac{e}{p^{s+k}}\,, \ \ u\in\{0,\ldots,m-1\}\,, \ v\in \{0,\ldots,p^k-1\}.
\end{equation*}
For the projection $r(e(\gamma'))\in\Z[\mup]$ one thus gets that
\begin{equation*}
r(e(\gamma'))=e(y+\frac w m)=e(\frac{f + wb}{bm}) \,, \ \ w\in\{0,\ldots,m-1\}\
\end{equation*}
which is repeated with multiplicity $p^k$. The equation \eqref{comrhon} follows. \endproof

\begin{cor}\label{corproj}
One has
\begin{equation}\label{projproj}
   r\circ \tilde\rho_n(x)=r\circ \tilde\rho_n(r(x))\qqq x\in \Z[\Q/\Z]\,, \ n\in \N.
\end{equation}
and
\begin{equation}\label{invfr}
    r\circ \tilde\rho_{p^k}(x)=p^k\sigma_{p^k}^{-1}(r(x))\qqq x\in \Z[\Q/\Z],~k\in\N\,.
\end{equation}
\end{cor}

\proof The two statements follow from \eqref{comrhon}.\endproof

\begin{defn}\label{defnxp} For $p$ a prime number, we denote by $X_p$  the space of all injective group
homomorphisms $\sigma:\bar\F_p^\times \to \Q/\Z$.
\end{defn}
The relation between $\W_0(\bar \F_p)$ and the abelian part $\Z[\Q/\Z]$ of the integral BC-algebra $\cH_\Z$ is described by the following lemma
\begin{prop}\label{bcrelatenew} Let $\sigma\in X_p$ and let $\tilde\sigma$ be the associated ring isomorphism
\begin{equation*}\tilde\sigma: \W_0(\bar \F_p)\stackrel{\sim}{\longrightarrow}\Z[\mup]\subset \Z[\Q/\Z].\end{equation*}
Then the Frobenius $F_n$ and Verschiebung  maps $V_n$ on $\W_0(\bar \F_p)$ are obtained by restriction of the ring endomorphisms $\sigma_n$ and the maps $\tilde\rho_n$ on  $\Z[\Q/\Z]$ by the formulas
\begin{equation}\label{restrictfr}
    \tilde\sigma\circ F_n=\sigma_n\circ \tilde\sigma\,, \ \ \tilde\sigma\circ V_n=r\circ\tilde\rho_n\circ \tilde\sigma.
\end{equation}
\end{prop}
\proof In section \ref{sectwzero} we recalled (\cf \cite{Gray} for details) that the Frobenius $F_n$ on $\W_0(A)$ is  given by $F_n(E,f)=(E,f^n)$. At the level of the divisor of the eigenvalues of $f$ (it is a divisor in the virtual case), \ie at the level of the associated element in $\Z[k^\times]$, $A=k=\bar\F_p$, the Frobenius $F_n$ corresponds to the group homomorphism $g\mapsto g^n$ (\cf Proposition \ref{wofield}).  The Verschiebung  maps $V_n$ are described by the  operation \eqref{vermat} on matrices. The  maps $V_n$ are additive and hence determined by the elements $V_n([\alpha])$ where $\alpha\in k^\times$. They correspond to the $n$ eigenvalues of the following matrix
\[
V_n(\alpha) = \begin{bmatrix}0 & 0&\cdots&\cdots&\alpha\\
1&0&0&\cdots&0\\
0&1&0&\cdots&0\\
\cdots&\cdots&\cdots&\cdots&\cdots\\
0&0&0&1&0
\end{bmatrix}.
\]
Since the $n$-th power of the above matrix is the multiplication by $\alpha$, all its eigenvalues fulfill the equation $\beta^n=\alpha$. In fact the characteristic polynomial of the above matrix is $P(X)=X^n-\alpha$. Let $n=p^km$, where $m$ is prime to $p$. Since $\bar\F_p$ is a perfect field, the root $\alpha^{p^{-k}}\in \bar\F_p$ of $X^{p^k}-\alpha$ is unique and it admits $m$ distinct roots of order $m$: $\beta^m=\alpha^{p^{-k}}$, which are the $m$ roots of $P(X)$. They take the form $\xi\beta_0$, with $\xi^m=1$. Thus the corresponding divisor is
\begin{equation}\label{divi}
\delta=\sum_{ \xi^m=1}p^k [\xi\beta_0].
\end{equation}
We now compare the above description of the divisor associated to $V_n(E,f)$ with $r\circ\tilde\rho_n(e(\gamma))$, where $\gamma=\sigma(\alpha)=\frac ab\in \mup\subset\Q/\Z$. The elements $\xi\beta_0\in \bar\F_p$ are the $m$ distinct roots of the equation $X^n=\alpha$. Similarly, with the notations of \eqref{comrhon}, the elements
\begin{equation*}
\frac{f + wb}{bm}\in \mup\subset \Q/\Z, \ \ w\in\{0,\ldots,m-1\}
\end{equation*}
are the $m$ solutions in $\mup$ of the equation $nz=\gamma$. One thus gets
\begin{equation*}
\tilde\sigma(\delta)=p^k\sum_{w=0}^{m-1} e(\frac{f + wb}{bm}).
\end{equation*}
Thus \eqref{comrhon} shows that
\begin{equation*}
   \tilde\sigma(V_n([\alpha]))=\tilde\sigma(\delta)=r\circ \tilde\rho_n(e(\gamma))=
   r\circ\tilde\rho_n\circ \tilde\sigma([\alpha]).
\end{equation*}
\endproof

\begin{thm}\label{bcrelate1} Let $\sigma\in X_p$.
The following formulas define a representation $\pi_\sigma$ of the  integral BC-system $\cH_\Z$ as additive endomorphisms of
$\W_0(\bar \F_p)$
\begin{equation}\label{repi}
    \pi_\sigma(x)\xi=\tilde\sigma^{-1}(r(x))\,\xi\,, \ \ \pi_\sigma(\mu_n^*)=F_n\,, \ \ \pi_\sigma(\tilde\mu_n)=V_n
\end{equation}
for all $\xi \in \W_0(\bar \F_p)$, $x\in \Z[\Q/\Z]$ and $n\in \N$.
\end{thm}

\proof By construction $x\mapsto \tilde\sigma^{-1}(r(x))$ is a homomorphism of the group ring $\Z[\Q/\Z]$ to $\W_0(\bar \F_p)$ and hence, by composition with the left regular representation, $\pi_\sigma$ gives a representation of $\Z[\Q/\Z]$. The $F_n$ and $V_n$ are additive. It remains to check the relations \eqref{presoverZ1} and \eqref{presoverZ2}. The latter ones follow from \eqref{rabi9} for the first two, and from $(1)$ and $(3)$ of Proposition \ref{rabi10}  for the last two. To check the first relation of \eqref{presoverZ1} one needs to show that
\begin{equation}\label{toshow}
V_n \pi_\sigma(x)F_n=\pi_\sigma(\tilde\rho_n(x)).
\end{equation}
One has $\pi_\sigma(x)=\pi_\sigma(r(x))$ for all $x\in \Z[\Q/\Z]$.
Thus, by applying \eqref{projproj}, one can replace $x$ by $r(x)$ without changing both sides of the equation. Thus we can assume that $x=\tilde\sigma(z)$ for some $z\in \W_0(\bar \F_p)$. Then $\pi_\sigma(x)$ is just the multiplication by $z$. One has by \eqref{restrictfr}
\begin{equation*}
r\circ \tilde\rho_n(x)=r\circ \tilde\rho_n(\tilde\sigma(z))=\tilde\sigma\circ V_n(z).
\end{equation*}
Thus $\pi_\sigma(\tilde\rho_n(x))$ is the multiplication by $V_n(z)$ and \eqref{toshow} follows from \begin{equation*}
V_n(zF_n(\xi))=V_n(z)\xi\qqq \xi \in \W_0(\bar \F_p)
\end{equation*}
which is statement $(2)$ of Proposition \ref{rabi10}. Let us check the other two relations of \eqref{presoverZ1}. The second one means
\begin{equation*}
F_n\pi_\sigma(x)=\pi_\sigma(\sigma_n(x))F_n
\end{equation*}
and since $r\circ\sigma_n=\sigma_n\circ r$ we can assume as before that $x=\tilde\sigma(z)$, for some $z\in \W_0(\bar \F_p)$. Then  $\pi_\sigma(x)$ is  the multiplication by $z$ and, by \eqref{restrictfr}, $\pi_\sigma(\sigma_n(x))$ is the multiplication by $F_n(z)$. The required equality then follows since $F_n$ is multiplicative. The last relation of \eqref{presoverZ1} means
\begin{equation*}
\pi_\sigma(x)V_n=V_n\pi_\sigma(\sigma_n(x))
\end{equation*}
and assuming $x=\tilde\sigma(z)$ it reduces to
\begin{equation*}
zV_n(\xi)=V_n(F_n(z)\xi) \qqq \xi \in \W_0(\bar \F_p)
\end{equation*}
which in turn follows from statement $(2)$ of Proposition \ref{rabi10}.\endproof

\section{The Witt vectors functor and the truncation quotients}\label{sectbigwitt}

In this section we provide a short overview on the construction of the universal Witt scheme in the form that is most suitable to the applications contained in this paper, for more details we refer to \cite{witt69,Mu,Ca,Haz,Au,Ro,Hesselholt}. In the second part of the section we connect the universal ring $\W_0(A)$ with $\W(A)$. \vspace{.05in}

The construction of the ring of big Witt vectors (or generalized Witt vectors) is described by a covariant endofunctor $\W: \An\to \An$ in the category of commutative rings (with unit). For $A\in{\rm obj}(\An)$, and as a functor to the category of sets, one defines
\[
\W(A) = A^\N=\{(x_1,x_2,x_3,\ldots)|x_i\in A\}.
\]
To a truncation set $N\subseteq\N$ (\ie a subset of $\N$ which contains every positive divisor of each of its elements), one associates the truncated functor
\[
\W_N: \An \to \Se,\qquad \W_N(A) = A^N.
\]
 As a functor to the category of sets, $\W_N$ is left represented by the polynomial ring $R_N = \Z[x_n|n\in N]$. Then it follows that the big Witt vectors functor $\W = \W_\N$ is left represented by the symmetric algebra $\sym = \Z[x_1,x_2,x_3,\ldots]$
\begin{equation}\label{rep}
\W(A) = \Hom_{\An}(\sym,A)\qquad\forall A\in{\rm obj}(\An).
\end{equation}
As an endofunctor in the category of commutative rings $\W_N: \An \to \An$ is {\em uniquely} determined by requiring that for any commutative ring $A$ and for any $n\in N$, the following map, called the $n$-th ghost component is a ring homomorphism
\begin{equation}\label{ghost0}
gh_n: \W_N(A)\to A,\quad gh_n(x)= \sum_{d|n} d x_d^{n/d}.
\end{equation}
For $t$ a variable, the functorial bijection of sets
\begin{equation}\label{phiaa}
\varphi_A: \W(A) \to \Lambda(A)=1+tA[[t]],\qquad x = (x_n)_{n\in\N}\mapsto f_x(t) = \prod_{n\in\N}(1-x_nt^n)^{-1}
\end{equation}
 transports the ring structure from $\W(A)$ to the multiplicative abelian group $\Lambda(A)$  of power series over $A$ with constant term $1$, under the usual multiplication of power series (the power series $1$ acts as the identity element). In other words one has
 \begin{equation*}
\varphi_A(x+y)=\varphi_A(x)\varphi_A(y)\qquad\forall x,y\in \W(A).
\end{equation*}
 To make the description of the corresponding product $\star$ on $\Lambda(A)$ more explicit one introduces first the $n$-ghost components $w_n: \Lambda(A) \to A$, $n\in\N$, which  are defined by the formula
\[
w(f) = w(1+a_1t+a_2t^2+a_3t^3+\cdots) =  w_1t+w_2t^2+\cdots = t\frac d{dt}(\log(f(t))).
\]
For example, the first three ghost components are given by the universal formulas
\[
w_1(f) = a_1,\quad w_2(f)=-a_1^2+2a_2,\quad w_3(f)=a_1^3-3a_1a_2+3a_3.
\]
For products of the form $\prod_{k=1}^m(1-\xi_kt)^{-1} = 1+a_1t+a_2t^2+\cdots = f(t)$ this means that
\begin{align*}
w_1t+w_2t^2+w_3t^3+\cdots &= t\frac d{dt}(\log(f(t)))= t\frac d{dt}\sum_{k= 1}^m\log((1-\xi_kt)^{-1}) =\\ &=\sum_{i=1}^\infty(\xi_1^i+\xi_2^i+\cdots+\xi_m^i)t^i = \sum_{i=1}^\infty p_i(\xi)t^i.
\end{align*}
Thus the ghost components are given by the power sums in the $\xi_k$'s. Then, the product $\star$ on $\Lambda(A)$ is {\em uniquely} determined by requiring that these ghost components are (functorial) ring homomorphisms. In fact, distributivity and functoriality together force the multiplication of power series in $\Lambda(A)$ be expressed by the following rule
\begin{equation}\label{mult}
f(t) = \prod_i(1-\xi_it)^{-1},\quad g(t) = \prod_i(1-\eta_it)^{-1},\quad\Rightarrow\quad(f\star g)(t)=\prod_{i,j}(1-\xi_i\eta_jt)^{-1}
\end{equation}
where
\[
t\frac d{dt}(\log(\prod_{i,j}(1-\xi_i\eta_jt)^{-1}))=\sum_{n=1}^\infty p_n(\xi)p_n(\eta)t^n.
\]
It follows that multiplication according to \eqref{mult} translates into component-wise multiplication for the ghost components on $\Lambda(A)$. It is expressed by explicit polynomials with integral coefficients of the form
\begin{multline*}
(1+\sum a_nt^n)\star (1+\sum b_nt^n)=1+a_1 b_1 t+\left(a_1^2 b_1^2-a_2 b_1^2-a_1^2 b_2+2 a_2 b_2\right) t^2+\\
+(a_1^3 b_1^3-2 a_1 a_2 b_1^3+a_3 b_1^3-2 a_1^3 b_1 b_2+5 a_1 a_2 b_1 b_2-3 a_3 b_1 b_2+a_1^3 b_3-3 a_1 a_2 b_3+3 a_3 b_3) t^3+\cdots
\end{multline*}
The ghost components $gh_n(x)$ of a Witt vector $x=(x_1,x_2,x_3,\ldots)\in\W(A)$ become  the ghost components of $\varphi_A(x)$, \ie
\begin{equation}\label{ghost}
gh_n: \W(A)\to A,\quad gh_n(x) = w_n(\varphi_A(x)).
\end{equation}
It follows that the bijection $\varphi_A: \W(A) \to \Lambda(A)$ becomes a ring isomorphism.\vspace{.05in}

Note moreover that the homomorphism of abelian groups  $L: \W_0(A)\to \Lambda(A)$ of \eqref{charpol}, preserves the product, \ie
\begin{equation}\label{multi8}
    L((E,f)\otimes (F,g))=L((E,f))\star L((F,g))
\end{equation}
so that it defines an injective ring homomorphism.\vspace{.05in}

Two Witt vectors  $x,y\in\W(A)$ are added and multiplied by means of universal polynomials with
integer coefficients
\[
x+_{\W}y = (\mu_{S,1}(x,y),\mu_{S,2}(x,y),\ldots),\quad x\times_{\W}y = (\mu_{P,1}(x,y),\mu_{P,2}(x,y),\ldots).
\]
The polynomials $\mu_{S,i},\mu_{P,j}$ are recursively
computed using the ghost components by the formulas
\[
gh_n(\mu_{S,1}(x,y),\mu_{S,2}(x,y),\ldots)=gh_n(x)+gh_n(y),
\]
\[
 gh_n(\mu_{P,1}(x,y),\mu_{P,2}(x,y),\ldots) = gh_n(x)gh_n(y).
\]
Notice that the polynomials $gh_n(x)$ depend only on the $x_d$ for $d$ a divisor of $n$, hence  the $n$-th addition and multiplication polynomials $\mu_{S,n}$, $\mu_{P,n}$ are polynomials that
only involve the $x_d$ and $y_d$ with $d$ a divisor of $n$. Thus, for a truncation set $N\subseteq\N$, the polynomial ring $R_N = \Z[x_n|n\in N]$ is a sub Hopf algebra and a sub co-ring object of
$\sym$, this means that it defines a quotient functor, which coincides with
$\W_N$. This result applies in particular to the truncation set $N=\{p^n\mid n\geq 0\}$,
 where $p$ is a prime number. Thus  the p-adic Witt vectors $\W_{p^\infty}(A)$ can be interpreted as a functorial quotient of the big Witt vectors (similarly one obtains $\W_{p^n}(A)$ as the p-adic Witt vectors of length $n+1$).\vspace{.05in}

The Teichm\"uller representative is a multiplicative map which defines a section to the ghost map $gh_1$. If $N\subset\N$ is a truncation set, the Teichm\"uller representative is defined as
\begin{equation*}
[\cdot]_N: A \to \W_N(A),\quad a \mapsto [a]_N = ([a]_N)_{n\in N},\quad [a]_{N,n}=\begin{cases} a&\text{if $n=1$},\\0&\text{if $n>1$}.\end{cases}
\end{equation*}
One has $gh_n([a]_N)=a^n$ for all $n\in N$.\vspace{.05in}

On the functorial ring $\W(A)$ one can introduce several functorial operations which derive from
(the large number of) ring endomorphisms of $\sym$ and by applying the representability property \eqref{rep}. For instance, the Verschiebung (shift) additive functorial endomorphisms on $\W$ and its quotients, arise from the ring endomorphism \begin{equation*}{\bf V}_n: \sym \to \sym\,, \ \ x_i\mapsto \begin{cases}x_{i/n}&\text{if $i$ is divisible by $n$}\\0&\text{otherwise}\end{cases}\end{equation*}
which corresponds to the map $f(t)\mapsto f(t^n)$ in $\Lambda(A)$.\vspace{.05in}

For $N\subset\N$ a truncation set, the shift is the additive map given by
\[
V_n: \W_{N/n}(A) \to \W_N(A),\quad V_n((a_d|d\in N/n))=(a_m'|m\in N);~ a_m'=\begin{cases} a_d&\text{if $m=nd$}\\0&\text{otherwise}\end{cases}
\]
where $N/n=\{d\in\N|nd\in N\}$. This means that the composite with the ghost components is given by
\begin{equation}\label{ghost12}
   gh_mV_n = \begin{cases}n gh_{m/n}&\text{if $n$ divides $m$}\\0&\text{otherwise.}\end{cases}
\end{equation}\vspace{.05in}

The n-th Frobenius is the (unique) natural ring homomorphism
\[
F_n: \W_N(A) \to \W_{N/n}(A)
\]
which is defined on the ghost components by the formula $gh_rF_n=gh_{rn}$. Thus by definition the n-th Frobenius map makes the following diagram commute
\begin{equation*}
\begin{CD}
\W_N(A)@>gh>> A^N\\
@VF_nVV  @VVF_n^{gh}V\\
\W_{N/n}(A)@>gh>> A^{N/n}
\end{CD}
\end{equation*}
where $F_n^{gh}$ takes a sequence $(a_m|m\in N)$ to the sequence whose d-th component is $a_{dn}$.
At the level of the components $x_j$ of a Witt vector $x\in\W_N(A)$, the Frobenius $F_n$ is given by polynomials with integral coefficients. For instance, the following are the first $5$ components of $F_3(x)$
\begin{eqnarray}
F_3(x)_1&=& x_1^3+3 x_3 \nonumber\\
F_3(x)_2&=& x_2^3-3 x_1^3 x_3-3 x_3^2+3 x_6 \nonumber\\
F_3(x)_3&=&  -3 x_1^6 x_3-9 x_1^3 x_3^2-8 x_3^3+3 x_9 \nonumber\\
F_3(x)_4&=&  -3 x_1^9 x_3+3 x_1^3 x_2^3 x_3-18 x_1^6 x_3^2+3 x_2^3 x_3^2-36 x_1^3 x_3^3\nonumber\\
&& -24 x_3^4+x_4^3-3 x_2^3 x_6+9 x_1^3 x_3 x_6+9 x_3^2 x_6-3 x_6^2+3 x_{12} \nonumber\\
F_3(x)_5&=&  -3 x_1^{12} x_3-18 x_1^9 x_3^2-54 x_1^6 x_3^3-81 x_1^3 x_3^4-48 x_3^5+x_5^3+3 x_{15}.\nonumber
\end{eqnarray}
Note that when $p$ is a rational prime one has (\cf \cite{Rabi} Proposition 5.12)
\begin{equation}\label{frobppower}
    F_p(x)_m\equiv x_m^p~  ({\rm mod}~pA).
\end{equation}
One also has (\cf \cite{Rabi} Proposition 5.9)
\begin{equation}\label{rabi9bis}
    V_{nm}=V_n\circ V_m=V_m\circ V_n\,, \ \ F_{nm}=F_n\circ F_m=F_m\circ F_n
\end{equation}
where for the maps $F_n$ one assumes $nN\subset N$ and $mN\subset N$.

Proposition \ref{rabi10}  extends without change, (\cf \cite{Rabi} Proposition 5.10).
\begin{prop} \label{rabi10bis}  Let $N\subset\N$ be a truncation set, and $n\in N$ with $nN\subset N$. Let $A$ be a commutative ring and $x,y\in \W_N(A)$. Then\vspace{.05in}

$(1)$~$F_n\circ V_n(x)=nx$.\vspace{.05in}

$(2)$~$V_n(F_n(x)y)=xV_n(y)$.\vspace{.05in}

$(3)$~If $m$ is prime to $n$, one has $V_m\circ F_n=F_n\circ V_m$.\vspace{.05in}

$(4)$~One has $V_n(x)V_n(y)=nV_n(x y)$.
\end{prop}
\proof We refer to \cite{Rabi} Proposition 5.10.
The statement $(4)$ differs slightly from this reference, it can be checked directly using Proposition \ref{rabi10}. It implies that when $n$ is invertible in $\W_N(A)$ then
$\frac 1n V_n$ defines a ring endomorphism.\endproof

It is important to see how the description of the universal ring  $\W_0(A)$ fits with the definition of $\W(A)$.  There is a canonical ring monomorphism   $\W_0(A)\hookrightarrow \W(A)$ which is given as the composite of  the injective ring homomorphism $L: \W_0(A)\to \Lambda(A)$ as in \eqref{charpol} and of the  ring isomorphism  $\varphi_A^{-1}:\Lambda(A)\stackrel{\sim}{\to} \W(A)$ (\cf\eqref{phiaa})
\begin{equation}\label{natmap}
 \W_0(A)\to \Lambda(A)\simeq\W(A),\qquad   (E,f)\mapsto \det(1-tM(f))^{-1}.
\end{equation}
In the case $A=\bar \F_p$ the characteristic polynomial $\det(1-tM(f))=\det(1-tf)$ factorizes as a product of terms $(1-t\alpha_j)$ of degree one,  where the $\alpha_j\in \bar \F_p$ are the eigenvalues of $f$ (\cf \eqref{Ldelta}).
\begin{lem}\label{embedlem}
Let $[\cdot]: \bar \F_p\to \W(\bar \F_p),~ x\mapsto \tau(x):=[x]$ be the \te lift and let $\delta:\W_0(\bar\F_p)\to \Z[\bar\F_p^\times]$ be the isomorphism of \eqref{eigeniso}. Then the canonical map \eqref{natmap}   is given explicitly as
\begin{equation}\label{map2}
j:\W_0(\bar \F_p)\to\W(\bar \F_p),\quad j\circ \delta^{-1} : \Z[ \bar \F_p^\times]\ni \sum n_j \alpha_j\mapsto \sum n_j \tau(\alpha_j)\in \W(\bar \F_p).
\end{equation}
\end{lem}
This Lemma together with Theorem \ref{almthm} shows that the subring $\W_0(\bar \F_p)\subset\W(\bar\F_p)$ is just  the group ring $\Z[ \bar \F_p^\times]$ and is freely generated over $\Z$ by the \te lifts.

\section{The $p$-adic representations $\pi_\sigma$ of the BC-system}\label{sectcompl}

  In this section we shall implement the results of \cite{Ca,Ro,Au} to describe the ring $\W(\bar\F_p)$, then using the embedding with dense image $\W_0(\bar\F_p)\hookrightarrow \W(\bar\F_p)$, we will extend the representation $\pi_\sigma$ of $\cH_\Z$ on $\W_0(\bar\F_p)$ (Theorem \ref{bcrelate1}) to a representation of the integral BC-system on $\W(\bar\F_p)$. Such representation is the $p$-adic analogue of the irreducible complex representation \eqref{repell2N}.\vspace{.05in}

 We begin by recalling the definition of the isomorphism
 \begin{equation}\label{robisobis}
   \W(\bar\F_p)\simeq\W_{p^\infty}(\bar\F_p)^{I(p)}
\end{equation}
where $I(p)\subset \N$ is the set
of positive integers which are prime to $p$ and $p^\infty$ is the set of integer powers of $p$.
 At the  conceptual level, this isomorphism is a special case of the general functorial isomorphism holding for any commutative ring $A$ with unit
 (\cite{Au} Theorem 1)
\begin{equation}\label{auer}
   \W(A)=\W_{I(p)}(\W_{p^\infty}(A)).
\end{equation}
When
 $A$ is an $\F_p$-algebra, every element of $I(p)$ is invertible in $B=\W_{p^\infty}(A)$, thus one derives a canonical isomorphism $\W_{I(p)}(B)\simeq B^{I(p)}$ which is defined in terms of the ghost components. Let $\Z_{(p)}$ be the ring $\Z$ localized at the prime ideal $p\Z$ so that every element of $I(p)$ is invertible in $\Z_{(p)}$.
A central role,  in the ring  $\Lambda(\Z_{(p)})$,  is played by the Artin-Hasse exponential, this is the power series
\begin{equation}\label{artha}
  E_p(t)=  {\rm hexp}(t)={\rm exp}(t+\frac{t^p}{p}+\frac{t^{p^2}}{p^2}+ \cdots)\in\Lambda(\Z_{(p)}).
\end{equation}
The following properties are well known (\cf \cite{Au,Ro})
\begin{prop}\label{proj}

$(1)$~$E_p(t)$ is an idempotent of $\Lambda(\Z_{(p)})$.\vspace{.05in}

$(2)$~For $n\in I(p)$, the series $E_p(n)(t):=\frac 1nV_n(E_p)(t)\in \Lambda(\Z_{(p)})$  determine an idempotent. As $n$ varies in $I(p)$, the $E_p(n)$ form a partition of unity by idempotents.\vspace{.05in}

$(3)$~ For $n\notin p^\N$, $F_n(E_p)(t)=1(=0_\Lambda)$ and $F_{p^k}(E_p)(t)=E_p(t)$, $\forall k\in\N$.
\end{prop}
To check $(1)$ directly, one shows that there exists a unique sequence $(x_n)_{n\in \N}\in \W(\Z_{(p)})$ such that\vspace{.05in}

$\bullet$~$x_1=1$\vspace{.05in}

$\bullet$~$x_{p^k}=0$ for all $k>0$\vspace{.05in}

$\bullet$~$F_m(x)_{p^k}=0$ for all $m\in I(p)$ and $k\geq 0$.\vspace{.05in}

This follows by noticing that the coefficient of $x_{mp^k}$ in $F_m(x)_{p^k}$ is $m\in I(p)$ which is invertible in $\Z_{(p)}$, so that one determines the $x_n$ inductively. One then checks that the ghost components of $(x_n)_{n\in \N}\in \W(\Z_{(p)})$ are the same as those of $E_p(t)$, \ie $gh_n (x)$ is equal to $1$ if $n\in p^\N$ and is zero otherwise.

Note that any $n\in I(p)$ is invertible in $\Lambda(\Z_{(p)})$. Division by $n$ corresponds to the extraction of the $n$-th root of the power series $f(t)=1+g(t)$. Formally, this is given by the binomial formula
\begin{equation}\label{rootn}
  f^{\frac 1n}= (1+g)^{\frac 1n}=1+\frac 1n g+\ldots +\frac{\frac 1n(\frac 1n-1)\cdots (\frac 1n-k+1)}{k!}g^k+\ldots
\end{equation}
The $p$-adic valuation of the rational coefficient of $g^k$ is positive because $\frac 1n\in\Z_p$, thus this coefficient can be approximated arbitrarily by a binomial coefficient.
It follows from Proposition \ref{rabi10}, $(4)$ that $\frac 1nV_n$ is an endomorphism of $\Lambda(\bar\F_p)$ and also a right inverse of $F_n$.

One easily derives from \cite{Ca,Au,Ro} the following result
\begin{prop}\label{proj1} Let $A$ be an $\F_p$-algebra.

$(a)$~The map
\begin{equation}\label{winfiso}
   \psi_A: \winf(A) \to \Lambda(A)_{E_p},\quad x = (x_{p^n})_{n\in\N},\quad \psi_A(x)(t):= h_x(t)=\prod_\N E_p(x_{p^n}t^{p^n})
\end{equation}
 is an isomorphism onto the reduced ring $\Lambda(A)_{E_p}=\{x\in \Lambda(A)\mid x\star E_p=x\}$.\vspace{.05in}

$(b)$~For $n\in I(p)$, the composite $\psi_A^{-1}\circ F_n$ is an isomorphism of the reduced algebra
$\Lambda(A)_{E_p(n)}$ with $\winf(A)$.\vspace{.05in}

$(c)$~The composite
\begin{equation}\label{isowinf}
    \theta_A(x)=(\theta_A(x))_n=\psi_A^{-1}\circ F_n(x\star E_p(n)), \ \ n\in I(p),~x\in \Lambda(A)
\end{equation}
 is a canonical isomorphism $\theta_A:\Lambda(A)\to\winf(A)^{I(p)}=\W(A)$.\vspace{.05in}

$(d)$~The composite isomorphism $\Theta_A:=\theta_A\circ\varphi_A: \W(A)\to \winf(A)^{I(p)}$ is given explicitly on the components by
\begin{equation}\label{thetacomp}
    ( \Theta_A(x)_n)_{p^k}=F_n(x)_{p^k}\qqq x\in \W(A) \qqq n\in I(p).
\end{equation}
\end{prop}
\proof The first three statements follow from \cite{Ca} \S 3.b, \cite{Au}, Thm. 1 and Prop. 1, \cite{Ro} Thm. 9.15. We prove $(d)$.
Since the Frobenius $F_n$ is an endomorphism and $F_n(x\star E_p(n))=F_n(x)\star E_p$, one can rewrite
\eqref{isowinf} as
\begin{equation}\label{isowinfbis}
   (\theta_A(x))_n=\psi_A^{-1}(E_p\star F_n(x))\qqq n\in I(p).
\end{equation}
Thus, to show \eqref{thetacomp} it is enough to prove it for $n=1$. One needs to check that
for all $x\in \W(A)$, one has
\begin{equation*}
E_p\star \varphi_A(x)=\prod_\N E_p(x_{p^n}t^{p^n}).
\end{equation*}
Indeed, this follows from distributivity and the identity
\begin{equation*}
E_p\star (1-xt^n)^{-1}=\left\{
                         \begin{array}{ll}
                           1 & \hbox{if}\  n\notin p^\N\\
                           E_p(x t^{p^k}) & \hbox{if}\  n=p^k.
                         \end{array}
                       \right.
\end{equation*}
The above identity can be checked directly knowing that $(1-xt^n)^{-1}=V_n(\tau(x))$ and by applying the equality
\begin{equation*}
E_p\star (1-xt^n)^{-1}=E_p\star V_n(\tau(x))=V_n(F_n(E_p)\star\tau(x))
\end{equation*}
together with Proposition \ref{proj} $(3)$ and the equality
\begin{equation*}
\tau(y)\star f(t)=(1-yt)^{-1}\star f(t)=f(yt)
\end{equation*}
which holds
for any element $f(t)\in \Lambda(A)$.
In particular, for the \te lift $\tau(y)=[y]$ of an element  $y\in A$ one gets
\begin{equation}\label{isote}
    \theta_A(\tau(y))_n=\tau(y^n)\qqq n\in I(p)
\end{equation}
where, on the right hand side, $\tau$  denotes the original \te lift $\tau:  A\to \winf(A)$. Indeed one has
$F_n(\tau(y))=\tau(y^n)$.\endproof

\begin{cor}\label{fixedfr} Let $A$ be an $\F_p$-algebra. Then, the common fixed points of the endomorphisms $F_n: \W(A)\to \W(A)$ for $n\in I(p)$, are the elements of the form
\begin{equation}\label{scalmul}
    L(\lambda)=\sum_{m\in I(p)} \frac 1mV_m(E_p\star\lambda)\,, \ \ \lambda \in \winf(A).
\end{equation}
One also has
\begin{equation*}
\varphi_{A}(L(\lambda))=\prod_{n\in I(p)} h_\lambda(t^n)^{\frac 1n}.
\end{equation*}
\end{cor}
\proof Let $x\in \W(A)$ with $F_n(x)=x$ for all $n\in I(p)$. Then, it follows from \eqref{thetacomp} and \eqref{isowinf} that all the components $(\theta_A(x))_n$ are equal, so that for some $\lambda \in \winf(A)$ one has
\begin{equation*}
E_p(n)\star x=\frac 1n V_n(E_p\star\lambda)
\end{equation*}
and $x$ is of the required form. Conversely, by Proposition \ref{proj} $(3)$, one has $F_a(E_p)=0_\Lambda$ for all $a\in I(p),\, a\neq 1$. Thus when one applies $F_k$ to $\frac 1n V_n(E_p\star\lambda)$, one gets $1(=0_\Lambda)$ unless $k|n$ using Proposition \ref{rabi10bis} $(2)$, $(3)$. When $k|n$ one obtains $\frac 1m V_m(E_p\star\lambda)$,
with $m=n/k$. Thus the elements of the form \eqref{scalmul} are fixed under all $F_k$. \endproof

We now apply these results to the case $A=\bar\F_p$. We identify $\winf(\bar\F_p)$ with a subring of $\C_p$ (the $p$-adic completion of an algebraic closure of $\Q_p$). Let $\cun\subset \C_p$ be the completion of the maximal unramified extension of $\Q_p$. Then one knows that $\winf(\bar\F_p)=\zun$ is the ring of integers of $\cun$. With $\Theta$ the isomorphism of \eqref{thetacomp}, we have
\begin{equation}\label{thetacomp1}
    \Theta: \W(\bar\F_p)\stackrel{\sim}{\to} (\zun)^{I(p)}\,, \ \ ( \Theta(x)_n)_{p^k}=F_n(x)_{p^k}, \quad \forall n\in I(p),~\forall x\in\W(\bar\F_p).
\end{equation}
Thus $\Theta$ makes $\W(\bar\F_p)$ a module over $\zun$.

To the Frobenius automorphism of $\bar\F_p$ corresponds, by functoriality, a canonical automorphism $\fr$ of $\zun$ which extends to a continuous automorphism
\begin{equation}\label{frobb}
    \fr\in \Aut(\cun).
\end{equation}
We can now describe the $p$-adic analogues of the complex irreducible representations of the BC-system (\cf\eqref{repell2N}). We recall that $X_p$ denotes the space of all injective group
homomorphisms $\sigma:\bar\F_p^\times \to \C^\times$. The  choice of $\sigma\in X_p$ determines an embedding $\rho: \qcp\to \C_p$ of the
 cyclotomic field
generated by the abstract roots of  unity of order prime to $p$ inside $\C_p$.\vspace{.05in}

In the following we shall use the simplified notation $\zuni=\zun$. For $m\in I(p)$, we let $\epsilon_m$  be the vector in $\W(\bar\F_p)$ with only one non-zero component: $\epsilon_m(m)=1$.

\begin{thm}\label{bcrelate2} Let $\sigma\in X_p$. The representation $\pi_\sigma$ as in Theorem \ref{bcrelate1} extends by continuity to a representation of the  integral BC-algebra $\cH_\Z$ on
$\W(\bar \F_p)$. For $n\in I(p)$ and for $x\in \Z[\Q/\Z]$,  $\pi_\sigma (\mu_n)$, $\pi_\sigma(x)$ and $\pi_\sigma(\mu_n^*)$  are   $\zuni$-linear operators on $\W(\bar\F_p)$
\begin{align}\label{reppad}
\pi_\sigma (\mu_n)\epsilon_m = \epsilon_{nm}, \ \ \ \
\pi_\sigma(e(a/b))&\epsilon_m=\rho(\zeta_{a/b}^m)\epsilon_m \qqq a \in \Z,\ \forall b, m\in I(p)
\end{align}
\begin{equation}\label{munstar}
    \pi_\sigma(\mu_n^*)\epsilon_k=\left\{
                                      \begin{array}{ll}
                                        0 & \hbox{if}\, k\notin n\N\,  \\
                                        \epsilon_{k/n} & \hbox{if}\, k\in n\N\,.
                                      \end{array}
                                    \right.
\end{equation}
 One has $\pi_\sigma(x)=\pi_\sigma(r(x))$ for all $x\in \Z[\Q/\Z]$ ($r: \Z[\Q/\Z]\to \Z[\mup]$ the retraction as in \eqref{retract}) and
\begin{equation}\label{geomfrob}
   \pi_\sigma (\mu_p)=\fr^{-1},\quad \pi_\sigma(\mu_p^*) = \fr
\end{equation}
where $\fr$ is the Frobenius automorphism, acting componentwise as a skew-linear operator.
\end{thm}
\proof Theorem \ref{bcrelate1} and the density of $\W_0(\bar\F_p)$ in $\W(\bar\F_p)$ (\cf\eg\cite{Haz1},~1.8) show that $\pi_\sigma$
extends by continuity to a representation of the  integral BC-algebra $\cH_\Z$ on
$\W(\bar \F_p)$. In view of the invertibility of the elements $n\in I(p)$ in $\W(\bar\F_p)$,  the description of the representation $\pi_\sigma$ is simplified by using the elements $\mu_n=\frac 1n\tilde\mu_n$,   to stress the analogy with the complex case. It follows from Corollary \ref{fixedfr}  that the subring $\zuni$ of $\W(\bar \F_p)$ is the fixed subring for the action of  the operators $F_n$, for all $n\in I(p)$. For $n\in I(p)$, the operators $F_n$ are $\zuni$-linear likewise the $V_n$ (\cf Proposition \ref{rabi10}, $(2)$) which correspond to the $\tilde\mu_n$ by means of the representation $\pi_\sigma$. Thus we obtain the first equality in \eqref{reppad}. The operators $\pi_\sigma(e(a/b))$ are the multiplication operators (\cf Corollary~\ref{noncan}) by $\tau(e(a/b))$, thus they are $\zuni$-linear and the second equation in \eqref{reppad} follows from \eqref{isote}. By applying \eqref{repi} one has $\pi_\sigma(\mu_n^*)=F_n$ for all $n$. Taking $n=p$, one gets that $\pi_\sigma(\mu_p^*)=F_p$ which coincides with $\fr$ acting componentwise, as it follows from the commutation $F_p\circ \frac 1nV_n=\frac 1nV_n\circ F_p$ for $n\in I(p)$ and \eqref{frobppower}.
Since $\mu_p^*\mu_p=1$ and $\fr$ is invertible one gets \eqref{geomfrob}.
\endproof

\begin{defn}\label{theideal} We denote by $\cJ_p\subset \cH_\Z$ the two sided ideal generated by the elements
\begin{equation}\label{ele1}
    1-e(p^{-k})\qqq k\in \N.
\end{equation}
\end{defn}

\begin{prop}\label{ideal}
One has
$
    \cJ_p= \Ker\, \pi_\sigma
$
(\cf\eqref{reppad}) and the intersection $\Z[\Q/\Z]\cap \cJ_p$ is the ideal $\cJ_p^0$ of $\Z[\Q/\Z]$ generated
by the elements as in \eqref{ele1}.

The sequence of commutative algebras
\begin{equation}\label{qual}
 0\to \cJ_p^0\to  \Z[\Q/\Z]\stackrel{r}{\longrightarrow}\Z[\mup]\to 0
\end{equation}
is exact.
\end{prop}
\proof
Let $r=id_{\Z[\mup]} \otimes \epsilon:\Z[\Q/\Z]\to \Z[\mup]$ be the retraction map introduced in \eqref{retract}. By construction, one has
\begin{equation}\label{cjpo}
\cJ_p^0=\Ker(r).
\end{equation}
Since $\pi_\sigma(e(a/b))$ only depends upon $r(e(a/b))$ it follows that $ \cJ_p\subset \Ker\, \pi_\sigma$. One knows (\cf~\cite{ccm}, Lemma 4.8) that any element of the algebra $\cH_\Z$ can be written as a finite sum of monomials of the form
\begin{equation}\label{hzelem}
   \sum_{\{a,b\}\in\N^2\atop (a,b)=1}\tilde\mu_a\, x_{\{a,b\}}\,\mu_b^*,\qquad x_{\{a,b\}}\in \Z[\Q/\Z].
\end{equation}
We show that for any finite sum $X$ as in \eqref{hzelem} we have
\begin{equation}\label{cjo}
    \pi_\sigma(X)=0\implies x_{\{a,b\}}\in \cJ_p^0.
\end{equation}
It is enough to prove that $r(x_{\{a,b\}})=0$, $\forall a,b\in\N$ and since $\Z[\mup]$ is torsion free it suffices to show that $a r(x_{\{a,b\}})=0$, $\forall a,b\in\N$. Define $y: \Q_+^*\to \Z[\mup]$, $y(\frac ab):=a\, r(x_{\{a,b\}})$, then $y$ has finite support.
For any group homomorphism
\begin{equation*}
   \chi:\mup\to \cO^\times
\end{equation*}
there is a unique  ring homomorphism $h_\chi$ with
\begin{equation}\label{ringhomo}
    h_\chi:\Z[\mup]\to \cO,\qquad h_\chi(e(\gamma))=\chi(\gamma)\qqq \gamma\in \mup.
\end{equation}
This applies in particular, for any integer $j$, to $\chi=\rho^j$ where we view $\rho: \qcp\to \C_p$
as a group homomorphism $\rho: \mup\to \cO^\times$. One has
\begin{equation}\label{basechar}
    \bigcap_{j\in\Z}\,\Ker\, h_{\rho^j}=\{0\}
\end{equation}
since  an injective character of a finite cyclic group generates the dual group.
Let $n,m\in I(p)$ be relatively prime. Then one has for any $j\in I(p)$ and $z\in \cO$
\begin{equation}\label{compnm}
    (\pi_\sigma(X)z\epsilon_{jm})_{jn}=\sum_{k\in\Z}h_{\rho^j}(y(p^{-k}\frac nm))\,\fr^k(z).
\end{equation}
Thus if $\pi_\sigma(X)=0$ one has for all $j$ and $m,n$ as above
\begin{equation}\label{vanis}
   \sum_{k\in\Z}h_{\rho^j}(y(p^{-k}\frac nm))\,\fr^k(z)=0\qquad \forall z\in \cO.
\end{equation}
For $z$ a root of unity one has $\fr^k(z)=z^{p^k}$, thus the  polynomial
\begin{equation*}
    \sum_{k\in\Z}h_{\rho^j}(y(p^{-k}\frac nm))\, Z^{p^{k+n}}
\end{equation*}
vanishes, for $n$ large enough, on all roots of unity thus it is identically zero, hence all its coefficients must vanish \ie
\begin{equation}\label{vanis1}
  h_{\rho^j}(y(p^{-k}\frac nm))=0\qqq k\in\Z, \ j\in \N.
\end{equation}
It then follows from \eqref{basechar} that $y(\frac ab)=y(p^{-k}\frac nm)=0$, hence \eqref{cjo} holds and the proof that any element of $\Ker\, \pi_\sigma$ is in $\cJ_p$ is  complete. Finally, if $x\in \Z[\Q/\Z]$ belongs to  $\Ker\, \pi_\sigma$ one has $x\in \cJ_p^0$ by \eqref{cjo} and thus the intersection $\Z[\Q/\Z]\cap \cJ_p$ is the ideal $\cJ_p^0$.\endproof

\begin{defn}\label{defnhp}
We denote by $\cH_\Z^{(p)}$ the quotient by $\cJ_p$ of the subalgebra of $\cH_\Z$ generated by
$
    \Z[\Q/\Z]$, $\tilde\mu_n$, $\mu_n^*$, for $n\in I(p)
$.
\end{defn}
The algebra $\cH_\Z^{(p)}$ is generated by $\Z[\mup]$ the operators $\tilde\mu_n$ and $\mu_n^*$, for $n\in I(p)$ and its presentation is similar to the presentation of $\cH_\Z$. The relations are
\begin{equation}\label{presoverZ2bis}
\begin{array}{l}
\tilde\mu_{nm}= \tilde\mu_n  \tilde\mu_m\,, \ \mu_{nm}^* =\mu_{n}^*\mu_{m}^* \qqq n,m\in I(p)\\[3mm]
\mu_n^* \tilde\mu_n  =n\qqq n\in I(p) \\[3mm]
\tilde\mu_n\mu_m^* = \mu_m^*\tilde\mu_n \qqq n,m\in I(p)\quad\text{with}~ (n,m) = 1
\end{array}
\end{equation}
as well as the
relations
\begin{equation}\label{presoverZ1bis}
\begin{array}{l}
\tilde\mu_n x \mu_n^* = \tilde\rho_n(x)\,, \
\mu_n^* x = \sigma_n(x) \mu_n^*  \,, \
x \tilde\mu_n = \tilde\mu_n \sigma_n(x)
\end{array}
\end{equation}
where $\tilde\rho_n$, $n\in I(p)$ is defined by
\begin{equation}\label{tr}
    \tilde\rho_n(e(\gamma))=
\sum_{n\gamma'=\gamma}e(\gamma')\qqq \gamma\in \mup.
\end{equation}
Given an algebra $\cA$, an automorphism $\theta\in \Aut(\cA)$ and an integer $p$ we let $\cA\rtimes_{\theta,\,p}\Z$ be the subalgebra of the algebraic cross product   $\{\displaystyle{\sum_{n\in\Z}} a_n V^n\mid a_n\in \cA\}$ determined by the condition
\begin{equation}\label{crosp}
    a_{-n}\in p^n\cA\qqq n\in \N.
\end{equation}
If we let $V=U^*$ and $pV^{-1}= \tilde U$, then it is easy to see that $\cA\rtimes_{\theta,\,p}\Z$ is generated by $\cA,\tilde U,U^*$ with the relations
\begin{equation}\label{crospr}
 U^*\tilde U=p,\  \tilde U x U^* = p\theta^{-1}(x)\,, \
U^* x = \theta(x) U^*  \,, \
x \tilde U = \tilde U \theta(x) \qqq x\in \cA.
\end{equation}

\begin{prop} \label{crossfrob} There exists a unique automorphism $\fr\in \Aut(\cH_\Z^{(p)})$ such that
\begin{equation}\label{frfr}
    \fr(e(\gamma))=e(\gamma)^p\qqq \gamma\in \mup\,,\ \fr(\tilde\mu_n)=\tilde\mu_n\,, \
    \fr(\mu_n^*)=\mu_n^*\qqq n\in I(p).
\end{equation}
One derives an isomorphism
\begin{equation}\label{isocros}
    \cH_\Z/\cJ_p=\cH_\Z^{(p)}\rtimes_{\fr,\,p}\Z.
\end{equation}
\end{prop}
\proof The map $\gamma\to p\gamma$ defines an automorphism of $\mup$. Its linearization $\fr$ acts on  $\Z[\mup]$ and commutes with the endomorphisms $\sigma_n$ and $\tilde \rho_n$. In fact by applying the isomorphism of Proposition \ref{bcrelatenew}, $\fr$ corresponds to the Frobenius automorphism of $\bar\F_p$. Thus it extends to an automorphism $\fr\in \Aut(\cH_\Z^{(p)})$.

The second statement follows by comparing the presentation of $\cH_\Z/\cJ_p$ with that of the crossed product $\cH_\Z^{(p)}\rtimes_{\fr,\,p}\Z$ as in \eqref{crospr}. \endproof

\begin{prop} \label{bcrelate3}Let $\sigma\in X_p$.\vspace{.05in}

$(1)$~The restriction $\pi_\sigma|_{\cH_\Z^{(p)}}$ of the representation $\pi_\sigma$ (as in Theorem~\ref{bcrelate2}) to $\cH_\Z^{(p)}$ is $\zuni$-linear and indecomposable over $\zuni$.\vspace{.05in}

$(2)$~The representations $\pi_\sigma|_{\cH_\Z^{(p)}}$ are pairwise inequivalent.\vspace{.05in}

$(3)$~The representation $\pi_\sigma$ is linear and indecomposable over $\Z_p$.\vspace{.05in}

$(4)$~Two representations $\pi_\sigma$ and $\pi_{\sigma'}$ are equivalent over $\Z_p$ if and only if there exists $\alpha\in \Aut(\bar\F_p)$ such that $\sigma'=\sigma\circ \alpha$.

\end{prop}

\proof
$(1)$~The $\zuni$-linearity property is checked directly on the generators using Theorem \ref{bcrelate2}.
 It follows from \eqref{reppad} that the vector $\epsilon_1$ is cyclic for $\cH_\Z^{(p)}$, \ie $\pi_\sigma(\cH_\Z^{(p)})\,\epsilon_1$ is dense in $\W(\bar \F_p)=\zuni^{I(p)}$. One has
 \begin{equation}\label{subspace}
\zuni \epsilon_1=\{\xi\in \W(\bar \F_p)\mid \pi_\sigma(\mu_n^*)(\xi)=0\qqq n\neq 1,\  n\in I(p)\}.
	 \end{equation}
For any $\zuni$-linear continuous operator $T$ in the commutant of $\cH_\Z^{(p)}$ one has $\pi_\sigma(\mu_n^*) T\epsilon_1=T\pi_\sigma( \mu_n^*) \epsilon_1 =0$, $\forall n>1$ and by \eqref{subspace} there exists $\lambda\in \zuni$ such that $T\epsilon_1=\lambda\epsilon_1$. Thus  since $\epsilon_1$ is cyclic,  $T$ is given by the module action of $\lambda\in \zuni$.

$(2)$~By \eqref{reppad}, the action of $\pi_\sigma(e(\gamma))$ for $\gamma\in \mup$ on the subspace \eqref{subspace} is given by the multiplication by $\rho(\gamma)\in \zuni$. Thus $\rho$ is an invariant of the representation.

$(3)$~Any element of the commutant of the action of $\cH_\Z$ is
 given by the module action of $\lambda\in \zuni$, where $\lambda$ is fixed  for the action of the Frobenius on $\zuni$, \ie  $\lambda\in \Z_p$. This shows that $\pi_\sigma$ is indecomposable.

 $(4)$~We show first that if there exists $\alpha\in \Aut(\bar\F_p)$ such that $\sigma'=\sigma\circ \alpha^{-1}$, the representations $\pi_\sigma$ and $\pi_{\sigma'}$ are equivalent over $\Z_p$. Let $\tilde\alpha=\winf(\alpha)\in \Aut(\winf(\bar\F_p))=\Aut(\zuni)$
 and define $U:\zuni^{I(p)}\to \zuni^{I(p)}\,, \ \ (U\xi)_n=\tilde\alpha(\xi_n)\qqq n\in I(p)$.
 One has $U\epsilon_n=\epsilon_n$ for all $n\in I(p)$ and if $T$ is an $\zuni$-linear operator so is $UTU^{-1}$. It thus follows from \eqref{reppad} and \eqref{munstar} that $U\pi_\sigma (\mu_n)U^{-1}=\pi_\sigma (\mu_n)$ and $U\pi_\sigma (\mu_n^*)U^{-1}=\pi_\sigma (\mu_n^*)$. For $x\in \Z[\Q/\Z]$, $U\pi_\sigma (x)U^{-1}$ only depends on $r(x)$ and for $x=e(a/b)$, $b\in I(p)$, one has $U\pi_\sigma(e(a/b))U^{-1}\epsilon_m=\tilde\alpha(\rho(\zeta_{a/b}^m))\epsilon_m
    =\rho'(\zeta_{a/b}^m)\epsilon_m=\pi_{\sigma'}(e(a/b))$.

 Moreover since $\tilde\alpha$ commutes with $\fr$, it follows from \eqref{geomfrob}  that $U\pi_\sigma (\mu_p)U^{-1}=\pi_\sigma (\mu_p)$ and $U\pi_\sigma (\mu_p^*)U^{-1}=\pi_\sigma (\mu_p^*)$. Thus one gets the required equivalence.

 Conversely, assume that two representations $\pi_\sigma$ and $\pi_{\sigma'}$ are equivalent over $\Z_p$. By \eqref{subspace} the $\Z_p$-linear representation $\pi_\sigma$ (and similarly $\pi_{\sigma'}$) determines uniquely the following representation of $\Z[\mup]$ in the $\Z_p$-module $\zuni$
 \begin{equation}\label{reprep}
   \beta_\sigma(e(a/b))\xi=\rho(\zeta_{a/b})\xi\qqq \xi \in \zuni.
 \end{equation}
In turns this determines an extension of the $p$-adic valuation to the subfield $\qcp\subset \qcy$  generated over $\Q$ by $\emup$. Indeed the formula
\begin{equation}\label{valext}
    {\rm val}(x)=\inf\{k\geq 0\mid  \beta_\sigma(x)\zuni\subset p^k\zuni \}\qqq x\in \Z[\mup]
\end{equation}
only depends on the class of $x$ in $\qcp$ and extends uniquely to a valuation on $\qcp$. The conclusion then follows from Proposition \ref{try}.
 \endproof

\section{The KMS theory of the BC-system at a prime $p$}\label{sectL}

In \cite{BC} it was shown that the extremal, complex KMS states below critical temperature of the BC-system (\cf\eqref{BC-KMS1infty}) are of the form
\begin{equation}\label{BC-KMS1inftyfull}
\varphi_{\beta,\rho}(X) =\frac{\Tr(\pi_\rho(X)e^{-\beta
H})}{\Tr(e^{-\beta H})} \qqq X\in \cH_\Z
\end{equation}
where $H$ is the Hamiltonian operator of multiplication by $\log n$ in the canonical basis $\epsilon_n$ of the
Hilbert space $\ell^2(\N)$ and $\pi_\rho$ is the irreducible representation of the algebra $\cH_\Q$  given by
\begin{equation}\label{repell2N}
\pi_\rho (\mu_n)\epsilon_m = \epsilon_{nm}, \ \ \pi_\rho (\mu_n^*)=\pi_\rho (\mu_n)^*\,,  \ \
\pi_\rho(e(a/b))\epsilon_m=\rho(\zeta_{a/b}^m)\epsilon_m,
\end{equation}
where $\rho\in \hat\Z^*$ determines an embedding in $\C$ of the
 cyclotomic field $\qcy$
generated by the abstract roots of  unity. Thus the extremal KMS states $\varphi_{\beta,\rho}$ are directly computable using the representation $\pi_\rho$ and the explicit description of the Hamiltonian.

      In section~\ref{sectcompl}, we have described the $p$-adic analogue of the representation $\pi_\rho$. In this section, our goal is to obtain the $p$-adic analogue of the KMS states $\varphi_{\beta,\rho}$. The guiding equation is provided by the general algebraic formulation of the KMS condition which is described by the equality
\begin{equation}\label{kmscnt}
   \varphi(x\sigma(y))=\varphi(y\,x)\qqq x,y\in \cA,
\end{equation}
where $\varphi$ is a linear form on an algebra $\cA$ endowed with an automorphism $\sigma \in \Aut(\cA)$.
In our case the algebra is
\begin{equation}\label{cA}
    \cA=\cH_{\C_p}^{(p)}=\cH_\Z^{(p)}\otimes_\Z \C_p.
\end{equation}
In \S \ref{sectkmsauto} we introduce, using the Iwasawa logarithm as a substitute for the above complex Hamiltonian $H$, the automorphisms $\sigma^{(\beta)} \in \Aut(\cH_{\C_p}^{(p)})$. These  automorphisms are defined for $\beta$ in the ``extended $s$-disk" $D_p$ (\cf~\eqref{domain} below).
In \S \ref{sectexten} we shall show how to extend their definition from the domain $D_p$ to a covering $M$ of $\C_p$. The construction of the KMS states is based on the classical construction of the $p$-adic L-functions and $p$-adic polylogarithm and many properties that we obtain rely on the simplifications which occur when $\beta=1-k(p-1)$ ($k\in\Z$). In \S \ref{sectcycloid} we prove the identities in the cyclotomic field, involving Bernoulli polynomials,  which are behind the verification of the KMS condition. In \S \ref{sectlinform} we provide the construction of the linear forms $\varphi_{\beta,\rho}$ using some of the results from \cite{Wa} (\cf Chapter V). In \S \ref{sectkmscond} we prove that the functionals $\varphi_{\beta,\rho}$ fulfill the KMS condition with respect to the automorphism $\sigma^{(\beta)} \in \Aut(\cH_{\C_p}^{(p)})$. Unlike the complex case, this construction exhibits the (new) phenomenon of the invariance of the linear forms $\varphi_{\beta,\rho}$ under the symmetry of $\cH_{\C_p}^{(p)}$ given by the automorphism $e(\gamma)\mapsto e(-\gamma)$.\vspace{.05in}

Throughout this section we fix a finite, rational prime $p$ and an algebraic closure $\bar\Q_p$ whose completion is denoted $\C_p$. We also use the following notation
\begin{equation}\label{notationq}
q=4,\ \text{if}\ p=2,\qquad
    q=p,\ \text{if}\ p\neq 2\,.
\end{equation}
and
\begin{equation}\label{phiq}
    \varphi(q)=2,\ \text{if}\ p=2,\qquad \varphi(q)=p-1,\ \text{if}\ p\neq 2\,.
\end{equation}

We consider the ``extended $s$-disk"
\begin{equation}\label{domain}
  D_p:=  \{\beta\in \C_p\mid |\beta|_p<qp^{-1/(p-1)}>1\}\,, \ \
\end{equation}
and first develop the theory for $\beta\in D_p$. In \S \ref{sectexten} we shall explain how the Iwasawa construction of $p$-adic $L$-functions allows one to extend the whole theory from the domain $D_p$ to the covering  of $\C_p$ given by the multiplicative group $M$ which is the open disk of radius one and center $1$ in $\C_p$.

\subsection{The automorphisms $\sigma^{(\beta)} \in \Aut(\cH_{\C_p}^{(p)})$}\label{sectkmsauto}

Let $\Z_{(p)}^\times\subset \Q^\times$ be the multiplicative  group of rational fractions whose numerator and denominator are prime to $p$.
\begin{lem}\label{uniqueext} Let $r\in \Z_{(p)}^\times$. There exists a unique analytic function
\begin{equation}\label{rbeta}
     D_p\to \C_p,\quad \beta\mapsto r^{(\beta)}
\end{equation}
such that
\begin{equation}\label{interpol}
   r^{(\beta)}=r^\beta\qqq \beta=1-k \varphi(q).
\end{equation}
\end{lem}
\proof We recall that the Iwasawa logarithm $\log_p$ is the unique extension of the function defined in the open unit  disk centered at $1$ by
\begin{equation}\label{iwa}
   -\log_p(1-x)=\sum_{n=1}^\infty\frac{x^n}{n}\qqq x\in\C_p,\ |x|_p<1
\end{equation}
to a map $\log_p: \C_p^\times\to \C_p$ such that
\begin{equation}\label{iwa1}
   \log_p(xy)=\log_p(x)+\log_p(y)\qqq x,y\in \C_p
\,, \ \ \log_p(p)=0\,.
\end{equation}
One has $\log_p(-1)=0$ since $-1$ is a root of unity, and
\begin{equation}\label{logpint}
    |\log_p(r)|_p\leq q^{-1}\qqq r\in \Z_{(p)}^\times.
\end{equation}
Moreover the exponential function is defined by the series
\begin{equation}\label{expo}
   {\rm exp}(x)=\sum_{n=1}^\infty\frac{x^n}{n!}\qqq x\in \C_p\,, \ |x|_p<r_p=p^{-\frac{1}{p-1}}.
\end{equation}
We define
\begin{equation}\label{defnrbeta}
   r^{(\beta)}:=r\,{\rm exp}((\beta-1)\log_p(r))\qqq \beta\in D_p.
\end{equation}
   This is a well-defined, analytic function of $\beta\in D_p$ since $\beta-1\in D_p$ and thus $|(\beta-1)\log_p(r)|_p<r_p$ by \eqref{logpint}. We show that \eqref{interpol} holds. This follows from the equality
\begin{equation*}
    {\rm exp}(k \varphi(q)\log_p(r))=r^{k \varphi(q)}\qqq r\in \Z_{(p)}^\times, \ k\in \Z
\end{equation*}
which holds for $r=-1$ since $\varphi(q)$ is even. In general, \eqref{interpol} follows from the formula
\begin{equation}\label{agree}
    {\rm exp}(n \log_p(a))=a^n\qqq a\in \Z_p^*\,, \ n\in \varphi(q)\Z
\end{equation}
as shown in \cite{Wa} (Chapter~5, p. 52), where the notation
\begin{equation}\label{langle}
\langle a\rangle=
{\rm exp}(\log_p(a))
\end{equation}
is introduced. The uniqueness follows from the discreteness of the
set of zeros of analytic functions. \endproof

\begin{lem}\label{multext} Let $\beta \in D_p$, then
\begin{equation}\label{rbetamul}
   \Z_{(p)}^\times\ni r\mapsto  r^{(\beta)}\in \C_p^\times
\end{equation}
is a group homomorphism. Moreover, for $r\in \Z_{(p)}^\times$
\begin{equation}\label{mulmul}
    r^{(\beta_1)}r^{(\beta_2)}=r^{(\beta_1+\beta_2)}r^{(0)}\qqq \beta_j\in D_p
\end{equation}
\end{lem}
\proof This follows from \eqref{iwa1} and the equality (\cf \cite{Rob})
\begin{equation}\label{expadd}
    {\rm exp}(x_1+x_2)={\rm exp}(x_1){\rm exp}(x_2)\qqq x_j, \ |x_j|_p<r_p.
\end{equation}
\endproof
The standard notation for $r^{(0)}$ is $\omega(r)$: it is the unique $\varphi(q)$ root of unity which is congruent to $r$ modulo $q$. In particular one has
\begin{equation}\label{rootu}
    (r^{(0)})^{\varphi(q)}=1\qqq r\in \Z_{(p)}^\times.
\end{equation}
\begin{prop}\label{kmsauto}
$(1)$~For $\beta\in D_p$ there exists a unique automorphism $\sigma^{(\beta)} \in \Aut(\cH_{\C_p}^{(p)})$ such that
\begin{equation}\label{autocond}
   \sigma^{(\beta)} (\tilde\mu_a e(\gamma)\mu_b^*) =\left(\frac ba\right)^{(\beta)} \tilde\mu_a e(\gamma)\mu_b^*\qqq a,b\in I(p),\ \gamma\in \mup.
\end{equation}
$(2)$~One has
\begin{equation}\label{autoprop}
    \sigma^{(\beta_1)}\circ \sigma^{(\beta_2)}=\sigma^{(\beta_1+\beta_2)}\circ \sigma^{(0)}\qqq \beta_j\in D_p
\end{equation}
and $\sigma^{(0)}$ is an automorphism of order $\varphi(q)$.
\end{prop}
\proof It suffices to check that $\sigma^{(\beta)}$ preserves the presentation given by the relations \eqref{presoverZ2bis} and \eqref{presoverZ1bis}. This follows from the multiplicativity shown in Lemma \ref{multext}.
Similarly \eqref{autoprop} follows from \eqref{mulmul}. The last statement follows from \eqref{rootu}.\endproof

\subsection{Cyclotomic identities for the polylogarithm}\label{sectcycloid}

We recall that the Bernoulli polynomials $B_n(u)$ are defined inductively as follows
\begin{equation*}
    B_0(x)=1\,, \ \ B'_n(x)=nB_{n-1}(x)\,, \ \ \int_0^1B_n(x)dx=0.
\end{equation*}
Equivalently, these polynomials can be introduced using the generating function
\begin{equation}\label{berpo}
   F(u,t)=\frac{t e^{ut}}{e^t-1}=\sum_{n=0}^\infty B_n(u)\frac{t^n}{n!}.
\end{equation}
The first few are
\begin{eqnarray}
  B_0(u)&=&1 \nonumber\\
 B_1(u)&=&-\frac{1}{2}+u \nonumber\\
 B_2(u)&=&\frac{1}{6}-u+u^2 \nonumber\\
 B_3(u)&=&\frac{u}{2}-\frac{3 u^2}{2}+u^3\nonumber \\
 B_4(u)&=&-\frac{1}{30}+u^2-2 u^3+u^4 \nonumber\\
 B_5(u)&=&-\frac{u}{6}+\frac{5 u^3}{3}-\frac{5 u^4}{2}+u^5.\nonumber
\end{eqnarray}
These polynomials fulfill the equation $B_n(1-u)=(-1)^n B_n(u)$.  The Bernoulli numbers are $B_n=B_n(0)$.
Using \eqref{berpo}, one checks the identity (\cf~\cite{Wa},  Chapter~4, Proposition 4.1)
\begin{equation}\label{sumber}
    g^{n-1}\sum_{j=0}^{g-1}B_n(\frac{x+j }{g})=B_n(x)\,.
\end{equation}
We also introduce inductively the rational fractions $\ell_{\beta}(z)$ for $\beta\in -\N$, as follows
\begin{equation}\label{multilog1}
z\partial_z\ell_{\beta}(z)=\ell_{\beta-1}(z)\,, \ \ell_{0}(z)=\frac{z}{1-z}.
\end{equation}
 For $\alpha\in \Q/\Z$ we denote by $\zeta_\alpha\in \qcy$ the class of $e(\alpha)\in \Q[\Q/\Z]$ modulo the cyclotomic ideal (\cf Definition \ref{acf}). It is a root of unity whose order is the denominator of $\alpha$.

\begin{lem}\label{indeplem}
Let $n >1$, $a,b\in \N$. Then
\begin{equation}\label{alternat}
  b^{n-1}  \sum_{j=0}^{b-1} \zeta_{a/b}^j B_n(\frac jb)=\left\{
                                                          \begin{array}{ll}
                                                            -n\ell_{1-n}(\zeta_{a/b}), & \hbox{if}\ \zeta_{a/b}\neq 1 \\
                                                            B_n, & \hbox{if}\ \zeta_{a/b}= 1.
                                                          \end{array}
                                                        \right.
\end{equation}
\end{lem}
\proof The equality \eqref{alternat} for $\zeta_{a/b}= 1$ follows from \eqref{sumber}. Thus we can
assume that $z=\zeta_{a/b}\neq 1$. The Taylor expansion at $t=0$ of $(ze^t-1)^{-1}$ is given by
\begin{equation}\label{Tayl}
    (ze^t-1)^{-1}=(z-1)^{-1}-\sum_{n=1}^\infty \ell_{-n}(z)\frac{t^n}{n!}
\end{equation}
since $(z-1)^{-1}=-1-\ell_0(z)$ and $\partial_t$ agrees with $z\partial_z$.
Then for $b\in \N$ and $t$ such that $ze^{\frac tb }\neq 1$ one has
\[
\sum_{j=0}^{b-1} z^je^{\frac jb t}=\frac{z^be^t-1}{ze^{\frac tb }-1}.
\]
Since $z^b=1$, one derives
\[
\sum_{n=0}^\infty \left(\sum_{j=0}^{b-1} z^j B_n(\frac jb)\right)\frac{t^n}{n!}=\sum_{j=0}^{b-1} z^jF(\frac jb, t)=\frac{t }{ze^{\frac tb }-1}.
\]
Since $z\neq 1$, taking the Taylor expansion at $t=0$ using \eqref{Tayl}, gives the equality
\begin{equation}\label{alln}
   \sum_{j=0}^{b-1} z^j B_n(\frac jb)=-\frac{n}{b^{n-1}}\ell_{1-n}(z),\qquad\forall n>1\,.
\end{equation}
\endproof
\begin{prop}\label{alternat1} Let $n>1$, $a,b\in \N$.

$(1)$~The following sum only depends upon $n$ and  $\frac ab\in \Q/\Z$
\begin{equation}\label{yab}
    Y_n(a/b)=f^{n-1}  \sum_{j=0}^{f-1} \zeta_{a/b}^j B_n(\frac jf)\qqq f\in b\N,\ f\neq 0.
\end{equation}

$(2)$~One has
\begin{equation}\label{yab1}
\frac 1b \sum_{a=0}^{b-1}Y_n(a/b)=b^{n-1}B_n=b^{n-1}Y_n(0).
\end{equation}

$(3)$~For $g\geq 1$, $x^g\neq 1$ one has
\begin{equation}\label{dupli}
    \frac 1g\sum_{j=0}^{g-1} \ell_{1-n}(\zeta_{j/g}\, x)=g^{n-1}\ell_{1-n}(x^g).
\end{equation}
\end{prop}
\proof $(1)$~Follows from \eqref{alternat}. To obtain $(2)$, note that
\begin{equation*}
    \frac 1b \sum_{a=0}^{b-1}\zeta_{a/b}^j=0 \qqq j\neq 0 \ (b)\,, \ \ \frac 1b \sum_{a=0}^{b-1}\zeta_{a/b}^j=1 \qqq j= 0 \ (b).
\end{equation*}
$(3)$~One checks \eqref{dupli} as an identity between rational fractions by induction on $n\in \N$. It holds for $n=1$ by applying the operation $-z\partial_z\log()$ to both sides of the identity
\begin{equation*}
    \prod_{j=0}^{g-1}(1-\zeta_{j/g}z)=1-z^g.
\end{equation*}
To obtain \eqref{dupli} for $n$ assuming it for $n-1$ one applies the operation $z\partial_z$ to both sides of the identity for $n-1$.
\endproof
Combining \eqref{dupli} with \eqref{alternat} we obtain, using \eqref{yab1} when $\alpha\in \Z$
\begin{equation}\label{dupliY}
    \frac 1b\sum_{j=0}^{b-1}Y_n(\frac{\alpha+j}{b})=b^{n-1}Y_n(\alpha)\qqq \alpha \in \Q/\Z
\end{equation}

\subsection{The linear forms $\varphi_{\beta,\rho}$}\label{sectlinform}

In this section we shall provide a meaning to expressions of the form
\begin{equation}\label{padicsum}
   Z_\rho(\frac ab,\beta)= \sum_{m\in I(p)}\rho(\zeta_{a/b}^m)\, m^{-\beta},\qquad\beta\in D_p
\end{equation}
where $\frac ab\in\Q$, $b\in I(p)$ is an integer prime to $p$ and $\rho: \qcy\hookrightarrow\C_p$. Note that as a function of $m\in I(p)$,  $\rho(\zeta_{a/b}^m)$ only depends on the residue  of $m$ modulo $b$. We let $f=bp$ and decompose the sum \eqref{padicsum} according to the residue $\alpha$ of $m$ modulo $f$. One has $\Z/f\Z=\Z/p\Z\times \Z/b\Z$. The elements of $I(p)$ are characterized by the fact that their residues mod. $f$ are given by pairs $\alpha=(\alpha_1,\alpha_2)\in\Z/f\Z$, with $\alpha_1\neq 0$. For  $\alpha\in (\Z/p\Z)^\times\times \Z/b\Z$, we let $\tilde\alpha\in \N$ be the smallest integer with residue modulo $f$ equal to $\alpha$.  Then, the sum \eqref{padicsum} can be written as
\begin{equation}\label{padicsum1}
   Z_\rho(\frac ab,\beta)= \sum_{\alpha}\rho(\zeta_{a/b}^\alpha)\, \sum_{n\geq 0} (\tilde\alpha+fn)^{-\beta},\qquad\beta\in D_p.
\end{equation}
Notice that the first sum (over $\alpha$) in \eqref{padicsum1} only involves finitely many terms. Each infinite sum in \eqref{padicsum1} is of the form (with $z=\tilde\alpha/f$)
\begin{equation}\label{padicsum1.5}
\sum_{n\in \N} (\tilde\alpha+fn)^{-\beta}=f^{-\beta}\,\sum_{n\in \N} (z+n)^{-\beta},\qquad\beta\in D_p
\end{equation}
and it is well known that this expression retains a meaning in the $p$-adic context (\cf \cite{Wa} Chapter V). More precisely, the asymptotic expansion in the complex case, for $z\to \infty$ (this process goes back to Euler's computation of $\sum_1^\infty n^{-2}$)
\begin{equation*}
   \sum_{n=0}^\infty (z+n)^{-\beta}\sim \frac{z^{1-\beta}}{\beta-1}\sum_0^\infty \binom{1-\beta}{j}B_j z^{-j}
\end{equation*}
motivates the following precise formula, where we prefer to leave some freedom in the choice of the multiple $f$ of $bq$.

\begin{lem}\label{indlem} With $q$ as in \eqref{domain}, and $f\in \N$, $f\neq 0$, a multiple of $bq$, the expression
\begin{equation}\label{padidef1}
    Z_\rho(\frac ab,\beta,f):=\frac{1}{f}\sum_{1\leq c<f\atop c\notin p\N}\rho(\zeta_{a/b}^c)\,
    \frac{\langle c\rangle^{1-\beta}}{\beta-1}\sum_{j=0}^\infty \binom{1-\beta}{j}\left(\frac{f}{c}\right)^{j}B_j,\qquad \beta\in D_p,
\end{equation}
 defines a meromorphic function of $\beta\in D_p$ with a single pole at $\beta=1$.
\end{lem}
\proof It follows from \cite{Wa} (Proposition 5.8) and the inequality (\cf \cite{Wa} Theorem 5.10)
\begin{equation*}
    |\left(\frac{f}{c}\right)^{j}B_j|_p\leq p|f|_p^j
\end{equation*}
that the series
\begin{equation*}
\sum_{j=0}^\infty \binom{1-\beta}{j}\left(\frac{f}{c}\right)^{j}B_j
\end{equation*}
converges for $|\beta|_p<|f|_p^{-1}p^{-\frac{1}{p-1}}\geq qp^{-\frac{1}{p-1}}>1$.\endproof

\begin{lem} \label{indep}
For $\beta$ a negative odd integer of the form $\beta=1-m=1-k\varphi(q)$, and $f\in \N$, $f\neq 0$, $f$ a multiple of $bq$, one has, with $Y_m$ defined by \eqref{yab}
\begin{equation}\label{simpl1}
   Z_\rho(\frac ab,\beta,f)=-\frac{1}{m}\rho\left(Y_m\left(\frac ab\right)-p^{m-1}Y_m\left(\frac {pa}{b}\right)\right)\,.
\end{equation}
\end{lem}
\proof
For $1\leq c < f$, $c\notin p\N$, one has $\langle c\rangle^{1-\beta}=c^{m}$. The binomial coefficients $\binom{1-\beta}{j}$ in \eqref{padidef1} all vanish for $j> m$ and the sum defining $Z(\frac ab,\beta,f)$  is therefore finite. One has
\[
\frac{\langle c\rangle^{1-\beta}}{\beta-1}\sum_{j=0}^\infty \binom{1-\beta}{j}\left(\frac{f}{c}\right)^{j}B_j=
-\frac{c^{m}}{m}\sum_{j=0}^{m} \binom{m}{j}\left(\frac{f}{c}\right)^{j}B_j.
\]
Moreover for any integer $m>0$, the Bernoulli polynomials fulfill the equation
\[
\sum_{j=0}^{m} \binom{m}{j}z^{-j}B_j=z^{-m}B_m(z)\,.
\]
For $1\leq c < f$, $c\notin p\N$, one thus gets, taking $z=\frac cf$
\begin{equation}\label{simpl}
  \frac{1}{f}\frac{\langle c\rangle^{1-\beta}}{\beta-1}\sum_{j=0}^\infty \binom{1-\beta}{j}\left(\frac{f}{c}\right)^{j}B_j=
   -\frac{f^{m-1}}{m}B_m\left(\frac{c}{f}\right)\,.
\end{equation}
One defines for any $c\in\N$
\begin{equation}\label{tcdef}
   T(c):= -\frac{f^{m-1}}{m}B_m\left(\frac{c}{f}\right).
\end{equation}
One has
\[
 Z_\rho(\frac ab,\beta,f)=\sum_{1\leq c<f\atop c\notin p\N}T(c)\rho(\zeta_{a/b}^c)=\sum_{0\leq c < f}T(c)\rho(\zeta_{a/b}^c)-\sum_{c=jp\atop 0\leq j < f/p}T(c)\rho(\zeta_{a/b}^c).
 \]
Since $b$ divides $f$, one derives
  \[
 \sum_{0\leq c< f}T(c)\rho(\zeta_{a/b}^c)=-\frac{f^{m-1}}{m}\sum_{0\leq c< f} \rho(\zeta_{a/b}^c) B_m\left(\frac{c}{f}\right)=-\frac 1m\rho(Y_m(\zeta_{a/b}))
 \]
 while, since $b$ divides $f/p=f'$ one gets
 \[
\sum_{c=jp\atop 0\leq j < f/p}T(c)\rho(\zeta_{a/b}^c)= -\frac{f^{m-1}}{m}\sum_{0\leq j< f/p} \rho(\zeta_{a/b}^{jp}) B_m\left(\frac{j}{f'}\right)=-\frac{p^{m-1}}{m}\rho(Y_m(\zeta_{a/b}^p)).
 \]
The equality \eqref{simpl1} follows.\endproof
\begin{cor}\label{indcoro} The function
\begin{equation}\label{padidef11}
    Z_\rho(\frac ab,\beta):=Z_\rho(\frac ab,\beta,f)
\end{equation}
 is independent of the choice of $f\in bq\N$, $f\neq 0$.
\end{cor}
\proof For two choices $f,f'$ the analytic function of $\beta\in D_p$
\begin{equation*}
    (\beta-1)(Z_\rho(\frac ab,\beta,f)-Z_\rho(\frac ab,\beta,f'))
\end{equation*}
vanishes at all negative integers $1-k\varphi(q)$ by  the equality \eqref{simpl1},  thus it is identically $0$.\endproof
\begin{defn}\label{defz} The following equation defines a linear form $\varphi_{\beta,\rho}$ on $\cH_\Z^{(p)}$ for any $\beta\in D_p$
\begin{equation}\label{varphidef1}
  \varphi_{\beta,\rho}(\tilde\mu_n e(\frac ab)\mu_m^*)=\left\{
                                                  \begin{array}{ll}
                                                    Z_\rho(\frac ab,\beta) & \hbox{if}\ n=m=1\,, \\
                                                    0 & \hbox{otherwise,}
                                                  \end{array}
                                                \right.
\end{equation}
for $n,m\in I(p)$ relatively prime.
\end{defn}
The next lemma will play an important role in the proof (\cf next section) that $\varphi_{\beta,\rho}$ fulfills the KMS condition.
\begin{lem}\label{tilderhon}
For any $n\in I(p)$ and $\beta\in D_p$, $\beta\neq 1$, one has
\begin{equation}\label{varphiprop}
    \varphi_{\beta,\rho}(\tilde\rho_n(X))=\langle n\rangle^{1-\beta}\varphi_{\beta,\rho}(X)\qqq X\in \Z[\mup]
\end{equation}
(\cf~\eqref{tr} for the definition of $\tilde\rho_n$).
\end{lem}
\proof After multiplication by $\beta-1$, both sides of \eqref{varphiprop}   are analytic functions of $\beta\in D_p$. Thus it is enough to show that \eqref{varphiprop} holds for $\beta=1-k\varphi(q)=1-m$. In this case one has $\langle n\rangle^{1-\beta}=n^{m}$ and, from \eqref{simpl1} one gets
\begin{equation*}
  \varphi_{\beta,\rho}(e(\gamma))=-\frac{1}{m}\rho\left(Y_m(\gamma)-p^{m-1}Y_m(p\gamma)\right)\qqq \gamma \in \mup.
\end{equation*}
To prove the equality \eqref{varphiprop} we can assume that $X=e(\alpha)$ for $\alpha\in \mup$. One has
\begin{equation*}
    \tilde\rho_n(X)=\sum_{j=0}^{n-1} e(\frac{\alpha + j}{n})
\end{equation*}
so that
\begin{equation*}
     \varphi_{\beta,\rho}(\tilde\rho_n(X))=-\frac{1}{m}
\sum_{j=0}^{n-1}\rho\left(Y_m(\frac{\alpha + j}{n})-p^{m-1}Y_m(p\frac{\alpha + j}{n})\right).
\end{equation*}
Then \eqref{varphiprop}
 follows from  \eqref{dupliY}. Since $p$ is prime to $n$ and the rational numbers $p\frac{\alpha + j}{n}\in \Q/\Z$ form the same subset as the set made by the
$\frac{p\alpha + j}{n}$, we derive
\begin{equation*}
   \sum_{j=0}^{n-1}Y_m(\frac{\alpha + j}{n})=n^m Y_m(\alpha)\,, \ \  \sum_{j=0}^{n-1}Y_m(p\frac{\alpha + j}{n})=n^m Y_m(p\alpha).
\end{equation*}
\endproof

\subsection{The KMS$_\beta$ condition}\label{sectkmscond}

The main result of this section is the following

\begin{thm}\label{thmkms} For any $\beta\in D_p$, $\beta\neq 1$ and $\rho: \qcp\hookrightarrow\C_p$, the linear form $\varphi_{\beta,\rho}$ fulfills the KMS$_{\beta}$ condition:
\begin{equation}\label{kmscnt1}
  \varphi_{\beta,\rho}(x\sigma^{(\beta)}(y))=\varphi_{\beta,\rho}(y\,x)\qqq x,y\in \cH_{\C_p}^{(p)}.
\end{equation}
Moreover the partition function is the $p$-adic $L$-function
\begin{equation}\label{partnonzero}
  Z(\beta):=  \varphi_{\beta,\rho}(1)=L_p(\beta,1)
\end{equation}
which does not vanish for $\beta\in D_p$.
\end{thm}
\proof We fix $x,y\in \cH_{\C_p}^{(p)}$, then after multiplication by $\beta-1$, both sides of \eqref{kmscnt1}  are analytic functions of $\beta\in D_p$. We first assume that $\beta\neq 1$; we shall consider the case $\beta=1$ separately later.  Since any element of the algebra $\cH_{\C_p}^{(p)}$ can be written as a finite linear combination of $\tilde\mu_n X\mu_m^*$, for $X\in \Z[\Q/\Z]$, we may assume that
\begin{equation*}
    x=\tilde\mu_n X\mu_m^*\,, \ \ y=\tilde\mu_s Y\mu_t^*
\end{equation*}
where $n,m\in I(p)$, $(n,m)=1$, $s,t\in I(p)$, $(s,t)=1$ and $X,Y\in \Z[\Q/\Z]$. Then, we use the presentation of $\cH_{\C_p}^{(p)}$ to compute
$
xy=\tilde\mu_n\,X \, \mu^*_m\,\tilde\mu_s\,Y \, \mu^*_t\,.
$
Let $u$ be the gcd of $m=um'$ and $s=us'$. One has
\begin{equation*}
\mu^*_m\,\tilde\mu_s=\mu^*_{m'}\,\mu^*_u\,\tilde\mu_u\tilde\mu_{s'}=u
\mu^*_{m'}\,\tilde\mu_{s'}=u\,\tilde\mu_{s'}\mu^*_{m'}
\end{equation*}
\begin{equation*}
\tilde\mu_n\,X \, \mu^*_m\,\tilde\mu_s\,Y \, \mu^*_t=u\,\tilde\mu_n\,X \,
\tilde\mu_{s'}\mu^*_{m'}\,Y \,
\mu^*_t =\,u\,\tilde\mu_n\,\tilde\mu_{s'}\sigma_{s'}(X)\sigma_{m'}(Y)
\mu^*_{m'}\, \mu^*_t
\end{equation*}
Let $v$ be the  gcd of $ns'=vw$ and $m't=vz$. One has
$
\tilde\mu_n\,\tilde\mu_{s'}=\tilde\mu_w\,\tilde\mu_{v}$, $\mu^*_{m'}\,
\mu^*_t=\mu^*_{v}\, \mu^*_z
$
\begin{equation*}
\tilde\mu_n\,X \, \mu^*_m\,\tilde\mu_s\,Y \, \mu^*_t=u\,\tilde\mu_w\,\tilde\mu_{v}\,\sigma_{s'}(X)\sigma_{m'}(y)\,\mu^*_{v}\,
\mu^*_z= u\,\tilde\mu_w\,
\tilde\rho_v(\sigma_{s'}(X)\sigma_{m'}(Y))\,\mu^*_z\,.
\end{equation*}
We obtain
\begin{equation}\label{fineform}
\tilde\mu_n\,X \, \mu^*_m\,\tilde\mu_s\,Y \, \mu^*_t= u\,\tilde\mu_w\,
\tilde\rho_v(\sigma_{s'}(X)\sigma_{m'}(Y))\,\mu^*_z\,, \ \ \frac wz=\frac nm\frac st.
\end{equation}
It follows  that unless $s=m$ and $t=n$ one has $\frac wz\neq 1$ and
\begin{equation}\label{prodkm}
  \varphi_{\beta,\rho}(x\sigma^{(\beta)}(y))=\varphi_{\beta,\rho}(y\,x)=0.
\end{equation}
Thus we can assume that $s=m$ and $t=n$. Then we have
\begin{equation*}
   x\sigma^{(\beta)}(y)=m\left(\frac nm\right)^{(\beta)}\tilde\mu_n\,XY\, \mu^*_n=
m\left(\frac nm\right)^{(\beta)}\tilde\rho_n(XY)
\end{equation*}
so that, by \eqref{varphiprop} one derives
\begin{equation}\label{sigmasig}
   \varphi_{\beta,\rho}(x\sigma^{(\beta)}(y))=m\left(\frac nm\right)^{(\beta)}
\langle n\rangle^{1-\beta} \varphi_{\beta,\rho}(XY).
\end{equation}
Similarly one has, by applying again \eqref{varphiprop}
\begin{equation*}
    yx=n\tilde\mu_m YX \mu_m^*\,, \ \ \varphi_{\beta,\rho}(y\,x)=n \langle m\rangle^{1-\beta} \varphi_{\beta,\rho}(XY).
\end{equation*}
Thus \eqref{kmscnt1} follows from the equality
\begin{equation*}
    m\left(\frac nm\right)^{(\beta)}
\langle n\rangle^{1-\beta}=n \langle m\rangle^{1-\beta}
\end{equation*}
which in turn derives from \eqref{defnrbeta} and \eqref{langle}.

Now, we turn to the normalization factor (\ie partition function) in \eqref{varphidef1} which is given by
\begin{equation}\label{partfun}
    Z(\beta):=\varphi_{\beta,\rho}(1)=\frac{1}{q}\sum_{1\leq c<q\atop c\notin p\N}
    \frac{\langle c\rangle^{1-\beta}}{\beta-1}\sum_{j=0}^\infty \binom{1-\beta}{j}\left(\frac{q}{c}\right)^{j}B_j\,.
\end{equation}
This is the  $p$-adic $L$-function for the character $\chi=1$ (\cf \cite{Wa}, Chapter 5, Theorem 5.11)
\begin{equation}\label{partfun1}
Z(\beta)=L_p(\beta,1).
\end{equation}
Moreover, notice that the Iwasawa construction of $L$-functions (\cf \cite{Wa}, Chapter 7, Theorem 7.10) yields a formal power series
    $\frac 12 g(T)\in \Z_p[[T]]^\times$
such that (with $q$ as in \eqref{notationq}) the following equality holds
\begin{equation}\label{iwas}
    L_p(\beta,1)=g\left((1+q)^\beta-1\right)/\left(1-(1+q)^{1-\beta}\right)\qqq \beta \in D_p.
\end{equation}
Since $\frac 12 g(T)\in \Z_p[[T]]^\times $ is invertible (\cf \cite{Wa} Lemma 7.12), this gives the required result.
\endproof
Note that $Z(\beta)$ has a pole at $\beta=1$, with residue given by
\[
\frac{1}{q}\sum_{1\leq c<q\atop c\notin p\N}1=\frac{\varphi(q)}{q}=
\frac{p-1}{p}\,.
\]
\begin{prop}\label{critval} When $\beta\to 1$ one has
\begin{equation}\label{limval}
   \lim_{\beta\to 1}Z(\beta)^{-1} Z_\rho(\frac ab,\beta)=
   \left\{
     \begin{array}{ll}
       1 & \hbox{if} \ \frac ab\in \Z\\
       0 & \hbox{otherwise.}
     \end{array}
   \right.
\end{equation}
\end{prop}
\proof Assume first that $\frac ab\notin \Z$. Then $\xi=\rho(\zeta_{a/b})$ is a non-trivial root of unity, whose order $m>1$ divides $b$ which is prime to $p$ and hence prime to $q$. Thus using the decomposition $\Z/bq\Z=\Z/q\Z\times \Z/b\Z$ we get
\[
\sum_{1\leq c<bq\atop c\notin p\N}\xi^c=\varphi(q)\sum_{n\in \Z/b\Z}\xi^n=0.
\]
If $\frac ab\in \Z$ the result follows from the above discussion.\endproof

Notice in particular that the limit of the functional values $ Z(\beta)^{-1} Z(\frac ab,\beta)$ as $\beta\to 1$  is independent of values of $\rho$ (\ie independent of the choice of  $\sigma\in X_p$). In the complex case, the functional values for $\beta>1$, are given by the formula \eqref{BC-KMS1infty}.
In that case, we shall now check directly that for $\beta\in \C$, $\Re(\beta)>1$, the functional values determine $\rho: \qcy \to \C$ as an  embedding  of the
 abstract cyclotomic field $\qcy$ in $\C$.

\begin{lem}\label{dense} $(1)$~Let $\lambda\in \hat\Z^*$, $\lambda\neq\pm 1$. Then the graph of the multiplication by $\lambda$ in $\Q/\Z$ is a dense subset of $\R/\Z\times \R/\Z$.

$(2)$~Let $\theta\in \Aut(\mup)$. Assume that $\theta \notin\{ \pm  p^\Z\}$. Then the graph of $\theta$ is dense in $\R/\Z\times \R/\Z$.
\end{lem}
\proof $(1)$~The set $G=\{(\alpha,\lambda\alpha)\mid \alpha\in \Q/\Z\}$ is a subgroup of $\R/\Z\times \R/\Z$ and so is its closure $\bar G$. If $G$ were not dense, then there would exist a non-trivial character $\chi$ of the compact group $\R/\Z\times \R/\Z$ whose kernel contains $\bar G$. Thus there would exist a non-zero pair $(n,m)\in \Z^2$ such that $n\alpha+m\lambda\alpha\in \Z$, for all $\alpha\in \Q/\Z$. This would imply that the multiplication by $\lambda\in \hat\Z^*$ in the group $\Q/\Z=\A_{\Q,\,f}/\hat\Z$ ($\A_{\Q,\,f}$ are the finite ad\`eles) ought to fulfill $n\alpha+m\lambda\alpha\in \hat\Z\qqq \alpha\in \A_{\Q,\,f}$. This implies $(n+m\lambda_p)\alpha\in \Z_p \qqq \alpha\in\Q_p$ and hence $n+m\lambda_p=0$ for all primes $p$.
If $n/m\notin\{\pm 1\}$, this contradicts the fact that $\lambda\in \hat\Z^*$ \ie $\lambda_p\in \hat\Z_p^*$ for all $p$.

$(2)$~By Lemma \ref{Gaut} the  group $G_p=\prod_{\ell\neq p}\Z_\ell^*$ is the group of automorphisms of the
group $\mup$ viewed as the additive group $\Gamma=\prod^{res}_{\ell\neq p}\Q_\ell/\Z_\ell$. Let $\lambda\in G_p$ represent $\theta\in \Aut(\mup)$.  Then the same proof as in $(1)$ shows that if the graph of $\theta$ is not dense, there exists a non-zero pair $(n,m)\in \Z^2$ such that
$n+m\lambda_\ell=0$ for all primes $\ell\neq p$. It follows that $-n/m\in \{ \pm  p^\Z\}$ and $\theta \in\{ \pm  p^\Z\}$.\endproof

From Lemma \ref{dense} we derive that, if $f:\{z\in\C\mid|z|=1\}\to \C$ is a continuous non-constant function, and $\rho_j: \qcy \to \C$, are injective, an equality of the form
\begin{equation}\label{equf}
    f(\rho_1
(\zeta_{a/b}))=f(\rho_2
(\zeta_{a/b}))\qqq a/b\in \Q/\Z
\end{equation}
necessarily implies that $\rho_2=\rho_1$ or $\rho_2=\bar\rho_1$. In the latter case one also gets
\[
f(\bar z)=f(z)\qqq z,\ |z|=1.
\]
By uniqueness of the Fourier decomposition however, this case cannot occur if $f(z)=\sum_{n=1}^{\infty}n^{-\beta} z^n$, for $\Re\mathfrak e(\beta)>1$.\vspace{.05in}

Next, we fix an integer  $\beta=1-m=1-k\varphi(q)$, $k>0$, and we  investigate the dependence on $\rho$ in the expressions \eqref{simpl1}.\vspace{.05in}

 For a chosen pair of embeddings $\rho,\rho'$,  assume that $Z_\rho(\frac ab,\beta)=Z_{\rho'}(\frac ab,\beta)$ holds for all $a/b\in \mup$, \ie the equality holds for all fractions with denominator $b$ prime to $p$. It follows from \eqref{simpl1} that one has (with $\fr$ the Frobenius automorphism of $\cun$)
\begin{equation}\label{indp}
    (1-p^{m-1}\fr)^{-1}Z_\rho(\frac ab,\beta)=-\frac{b^{m-1}}{m}\sum_{1\leq c\leq b}\rho(\zeta_{a/b}^c)\,B_m(\frac{c}{b})\in \cun.
\end{equation}
Thus we get
\begin{equation}\label{startequ}
\sum_{1\leq c\leq b}\rho(\zeta_{a/b}^c)\,B_m(\frac{c}{b})=\sum_{1\leq c\leq b}\rho'(\zeta_{a/b}^c)\,B_m(\frac{c}{b})\qquad\forall a/b\in \mup.
\end{equation}
Since both $\rho$ and $\rho'$ are isomorphisms of the group of roots of unity in $\qcp$ with the group   of roots of unity in $\C_p$ of order prime to $p$, there exists an automorphism $\theta\in \Aut(\mup)$ such that $\rho'(\zeta_{a/b})=\rho(\zeta_{\theta(a/b)})$
for all $a/b\in \mup$. One has
\[
\sum_{1\leq c\leq b}\rho'(\zeta_{a/b}^c)\,B_m(\frac{c}{b})=\sum_{1\leq c\leq b}\rho(\zeta_{\theta(c/b)}^{a})\,B_m(\frac{c}{b})=\sum_{1\leq c\leq b}\rho(\zeta_{a/b}^{c})\,B_m(\theta^{-1}(\frac{c}{b})).
\]
 By uniqueness of the Fourier transform for the finite group $\Z/b\Z$,  \eqref{startequ} yields the equality
\begin{equation}\label{sequ}
    B_m(\theta^{-1}(\frac{c}{b}))=B_m(\frac{c}{b})\qquad\forall c/b\in \mup.
\end{equation}

\begin{lem}\label{symbreak}
Let $p>2$ and let $\theta\in \Aut(\mup)$. If $\theta\in\{ \pm 1\}$ one has
\begin{equation}\label{diff}
   Z_\rho(\frac ab,\beta)=Z_{\theta\circ\rho}(\frac ab,\beta)\qqq a/b\in \mup, \ \beta\in D_p.
\end{equation}
If $\theta \notin\{ \pm 1\}$ and $\beta=1-m=1-k\varphi(q)$, $k>0$, then the functionals $Z_\rho(\cdot,\beta)$ and
$Z_{\theta\circ\rho}(\cdot,\beta)$ are distinct.
\end{lem}
\proof To prove \eqref{diff} we can assume that $\theta=-1$ \ie that $\theta(\zeta_{a/b})=\zeta_{a/b}^{-1}$ for all $a/b\in \mup$. Then we have, with $\rho'=\theta\circ\rho$: ~ $\rho'(\zeta_{a/b}^c)=\rho(\zeta_{a/b}^{b-c})$.
Let first $\beta=1-m=1-k\varphi(q)$. One has
\begin{equation*}
\sum_{1\leq c\leq b}\rho'(\zeta_{a/b}^c)\,B_m(\frac{c}{b})=\sum_{0\leq c\leq b-1}\rho(\zeta_{a/b}^c)\,B_m(\frac{b-c}{b})\,.
\end{equation*}
Since $m=k\varphi(q)$ is even, the Bernoulli polynomial $B_m$ fulfills the equality
\begin{equation}\label{onemin}
    B_m(1-x)=B_m(x)\qqq m\in 2\N\,.
\end{equation}
Thus  \eqref{diff} follows for all values $\beta=1-m=1-k\varphi(q)$. Since these values admit $0$ as an accumulation point, one derives the equality of the analytic functions on their domain $D_p$.

Now, we assume that $\theta \notin\{ \pm  p^\Z\}$. Then it follows from Lemma  \ref{dense} that the graph of $\theta$ is dense in $\R/\Z\times \R/\Z$. Thus \eqref{sequ} implies that $B_m(x)$ is constant which is a contradiction. It remains to show that for non-zero powers $p^a$ of $p$ one cannot have an equality of the form
\[
B_m(x)=B_m(p^ax-[p^ax])\qqq x\in [0,1]
\]
where $[p^ax]$ is the integral part of $p^ax$. In fact, this would imply that $B_m(x)-B_m(p^ax)$ has infinitely many zeros, thus $B_m(x)=B_m(p^ax)$ which is a contradiction.
\endproof

\subsection{Extension of the KMS$_\beta$ theory to the covering of $\C_p$}\label{sectexten}

In this section we show that the construction of the KMS$_\beta$ states $\varphi_{\beta,\rho}$, for $\beta\in D_p$, extends naturally to the covering of $\C_p$ defined by the following group homomorphism
\begin{equation}\label{cover1}
     M=D(1,1^-)\ni\lambda\mapsto \beta=\ell(\lambda)=\frac{\log_p\lambda}{\log_p(1+q)}\in \C_p
\end{equation}
where $M=D(1,1^-)$ is the open unit disk in $\C_p$ with radius $1$, viewed as a multiplicative group. Up to the normalization factor $\log_p(1+q)$, this group homomorphism coincides with the definition of the Iwasawa logarithm, it is surjective with kernel  the subgroup of roots of unity of order a $p$-power (\cf \cite{Rob}, Theorem p. 257)  and it defines by restriction  a bijection
\begin{equation}\label{bij}
    \ell: \{\lambda\in M\mid |\lambda-1|_p<p^{-1/(p-1)}\}\stackrel{\sim}{\to} D_p
\end{equation}
whose inverse is given by the map
\begin{equation}\label{invbij}
     D_p\ni\beta\mapsto \psi(\beta)=(1+q)^\beta={\rm exp}(\beta\log_p(1+q)).
\end{equation}
By construction, this local section is a group homomorphism which allows one to view the additive group $D_p$ as a subgroup of $M$.

 We start by extending the definition of the functions $r^{(\beta)}$ as in \eqref{defnrbeta} which were implemented in the construction  of the automorphisms $\sigma^{(\beta)}\in \Aut(\cH_{\C_p}^{(p)})$ (\cf Proposition \ref{kmsauto}). For $r\in \Z_{(p)}^\times$ the equality
\begin{equation}\label{ipr}
   i_p(r)=\frac{\log_p(r)}{\log_p(1+q)}\in \Z_p
\end{equation}
defines a group homomorphism from $\Z_{(p)}^\times$ to the additive group $\Z_p$.
\begin{lem}\label{agreelem}
For $\beta\in D_p$, $r\in \Z_{(p)}^\times$ and $\lambda=(1+q)^\beta$ one has
\begin{equation}\label{lambdapr}
 \langle r\rangle^\beta=\lambda^{i_p(r)}\,, \ \   r^{(\beta)}=\omega(r)\,\lambda^{i_p(r)}.
\end{equation}
\end{lem}
\proof One has $\log_p(r)=i_p(r)\log_p(1+q)\in q\Z_p$. Thus $|\beta \log_p(r)|_p<p^{-1/(p-1)}$ and
\begin{equation*}
    \langle r\rangle^\beta={\rm exp}(\beta \log_p(r))=
{\rm exp}(\beta i_p(r)\log_p(1+q))=(1+q)^{\beta i_p(r)}=\lambda^{i_p(r)}.
\end{equation*}
The second equality follows from the definition \eqref{defnrbeta}.\endproof
Proposition \ref{kmsauto} and its proof thus extend from $D_p$ to $M$. This means that for $\lambda \in M$  there exists a unique automorphism $\sigma[\lambda] \in \Aut(\cH_{\C_p}^{(p)})$ such that
\begin{equation}\label{autocondbis}
   \sigma[\lambda] (\tilde\mu_a e(\gamma)\mu_b^*) =\omega(b/a)\lambda^{i_p(b/a)} \tilde\mu_a e(\gamma)\mu_b^*\qqq a,b\in I(p),\ \gamma\in \mup.
\end{equation}
 Next, we extend the construction of the linear forms $\varphi_{\beta,\rho}$ given in \S \ref{sectlinform}. It is sufficient to extend the definition of the functions of Lemma \ref{indlem} (which we proved to be independent of the choice of $f\neq 0$ multiple of $bq$)
\begin{equation}\label{padidef1bis}
    Z_\rho(\frac ab,\beta):=\frac{1}{f}\sum_{1\leq c<f\atop c\notin p\N}\rho(\zeta_{a/b}^c)\,
    \frac{\langle c\rangle^{1-\beta}}{\beta-1}\sum_{j=0}^\infty \binom{1-\beta}{j}\left(\frac{f}{c}\right)^{j}B_j,\qquad \beta\in D_p.
\end{equation}
To define the sought for extension it is convenient to express the above function in terms of the $p$-adic $L$-functions
$L_p(\beta,\chi)$ associated to even Dirichlet characters   of conductor $f_\chi$ prime to $p$.
By definition, a Dirichlet character $\chi$ is a character of the multiplicative group $\hat\Z^*$ and its conductor $f_\chi$ is the integer such that the kernel of $\chi$ is the kernel of the projection $\hat\Z^*\to (\Z/f_\chi\Z)^*$.
The definition of $L_p(\beta,\chi)$ is similar to \eqref{padidef1bis} precisely as follows
\begin{equation}\label{llsim}
    L_p(\beta,\chi):=\frac{1}{f}\sum_{1\leq c<f\atop c\notin p\N}\chi(c)\,
    \frac{\langle c\rangle^{1-\beta}}{\beta-1}\sum_{j=0}^\infty \binom{1-\beta}{j}\left(\frac{f}{c}\right)^{j}B_j
\end{equation}
where $f$ is any multiple of $pf_\chi$ and where $\chi$ has been extended to a periodic function of period $f_\chi$ vanishing outside $(\Z/f_\chi\Z)^*$. We recall that the $L$-function $L_p(\beta,\chi)$ is identically zero when the character $\chi$ is odd, \ie when $\chi(-1)=-1$ (\cf \cite{Wa} Remarks p. 57). Moreover when $\chi$ is even, non-trivial, and its conductor is prime to $p$, there exists an analytic function $H_\chi$ on $M$ such that (\cf \cite{Wa} Theorem 7.10)
\begin{equation}\label{LLL}
    L_p(\beta,\chi)=H_\chi((1+q)^\beta)\qqq \beta\in D_p\,.
\end{equation}
The extension of the functions $ Z_\rho(\frac ab,\beta)$ to $M$ is a consequence of the following
\begin{lem}\label{comblem}
For   any $a/b\in \mup$ there exists coefficients $c(d,\chi)\in \C_p$ such that
\begin{equation}\label{cdchi}
   Z_\rho(\frac ab,\beta)=\sum_{d|b,\,\chi}c(d,\chi)L_p(\beta,\chi)d^{-1}\langle d\rangle^{1-\beta}\prod (1-\chi(\ell)\ell^{-1}\langle  \ell\rangle^{1-\beta})
\end{equation}
where $d$ varies among the divisors of $b$, and, for fixed $d$, $\chi$ varies among the set of Dirichlet characters whose conductor $f_\chi$ divides $m=b/d$. The integers $\ell$ are the primes which divide $m/f_\chi$ but not $f_\chi$.
\end{lem}
\proof Let $b$ be an integer prime to $p$, and $g\in C(\Z/b\Z,\C_p)$. The expression
\begin{equation}\label{padidef15}
    Y(g,\beta):=\frac{1}{f}\sum_{1\leq c<f\atop c\notin p\N}g(c)\,
    \frac{\langle c\rangle^{1-\beta}}{\beta-1}\sum_{j=0}^\infty \binom{1-\beta}{j}\left(\frac{f}{c}\right)^{j}B_j,\qquad \beta\in D_p
\end{equation}
is independent of the choice of the multiple $f\neq 0$ of $bq$.
Let $\chi$ be a Dirichlet character (with values in $\C_p$) with conductor $f_\chi$ and let $m$ be a multiple of $f_\chi$. Then the following defines a multiplicative map from $\Z/m\Z$ to $\C_p$
\begin{equation}\label{zchi}
    z(\chi,m)(c)=\left\{
           \begin{array}{ll}
             \chi(c) & \hbox{if}\, c\in (\Z/m\Z)^*\\
             0 & \hbox{otherwise.}
           \end{array}
         \right.
\end{equation}
If $m$ divides $b$ and one replaces $\chi$ with $z(\chi,m)$ in \eqref{llsim} one obtains instead of $L_p(\beta,\chi)$ the function
\begin{equation}\label{zchi3}
 Y(z(\chi,m),\beta)=   L_p(\beta,\chi)\prod (1-\chi(\ell)\ell^{-1}\langle  \ell\rangle^{1-\beta})
\end{equation}
where the integers $\ell$ are the primes which divide $m/f_\chi$ without dividing $f_\chi$. Next, define for any divisor $d$ of $b$ and any function $h\in C(\Z/m\Z,\C_p)$, $m=b/d$,
\begin{equation*}
 e_d(h)(a)=\left\{
           \begin{array}{ll}
             h(a/d) & \hbox{if}\, d|a\\
             0 & \hbox{otherwise.}
           \end{array}
         \right.
\end{equation*}
One then gets
\begin{equation}\label{eded}
    Y(e_d(h),\beta)=d^{-1}\langle d\rangle^{1-\beta} Y(h,\beta).
\end{equation}
Thus using \eqref{zchi3} and \eqref{eded}  it is enough to prove that for any function $g\in C(\Z/b\Z,\C_p)$ there exists coefficients $c(d,\chi)\in \C_p$ such that
\begin{equation*}
    g(c)=\sum_{d|b,\,\chi}c(d,\chi)e_d(z(\chi,b/d)).
\end{equation*}
It is in fact enough to check this for $g=\delta_a$ where $a\in \Z/b\Z$. Let then $d$ be the gcd of $a$ and $b$. One has $\delta_a=e_d(\delta_c)$ where $c=a/d$ is prime to $m=b/d$. Moreover for any element
  $c\in (\Z/m\Z)^*$ one has
\begin{equation*}
    \delta_c(x)=\frac{1}{\varphi(m)}\sum_{\chi,\, f_\chi|m}\chi(c)^{-1}z(\chi,m)(x)\qqq x\in \Z/m\Z,
\end{equation*}
which gives the required equality.\endproof
We thus obtain the following extension of Theorem \ref{thmkms}
\begin{thm}\label{thethmkms} There exists an analytic family of functionals $\psi_{\lambda,\rho}$,
$\lambda\in M$, on $\cH_\Z^{(p)}$ such that\vspace{.05in}

$\bullet$~$\psi_{\lambda,\rho}(1)=1.$\vspace{.05in}

$\bullet$~$\psi_{\lambda,\rho}$ fulfills the KMS condition
\begin{equation}\label{kmscnt2}
  \psi_{\lambda,\rho}(x\sigma[\lambda](y))=\psi_{\lambda,\rho}(y\,x)\qqq x,y\in \cH_{\C_p}^{(p)}.
\end{equation}

$\bullet$~For $\beta\in D_p$ and $\lambda=(1+q)^\beta$ one has
\begin{equation*}
    \psi_{\lambda,\rho}=Z(\beta)^{-1}\varphi_{\beta,\rho}\,.
\end{equation*}
\end{thm}
\proof It follows from \eqref{iwas} that there exists an analytic function $z(\lambda)$ of $\lambda\in M$ such that
\begin{equation*}
  Z(\beta)^{-1}=(1+q-\lambda)z(\lambda)\,, \ \ \lambda=(1+q)^\beta.
\end{equation*}
By applying \eqref{LLL}, Lemma \ref{agreelem} and Lemma \ref{comblem}, we see that there exists, for $b\in I(p)$ and $a/b\notin \Z$, an analytic function $H_{a,b}(\lambda)$ of $\lambda\in M$ such that
\begin{equation*}
  Z_\rho(\frac ab,\beta)=H_{a,b}(\lambda)\,, \ \ \lambda=(1+q)^\beta.
\end{equation*}
This proves the existence of the analytic family of functionals $\psi_{\lambda,\rho}$ fulfilling the required conditions. \endproof

\section{Extension of the $p$-adic valuation to $\qcy$}\label{sectval}

For a global field $\K$ of positive characteristic (\ie a function field  associated to a projective, non-singular curve $C$ over a finite field $\F_q$) it is a well known fact that the space of valuations of the maximal abelian extension
$\K^{\rm ab}$ of $\K$  has a geometric meaning. In fact, for each finite extension
$E$ of $\bar\F_q\otimes_{\F_q}\K\subset \K^{\rm ab}$ the space
${\rm Val}(E)$ of (discrete) valuations of $E$ is turned into an algebraic, one-dimensional scheme
  whose non-empty open sets are the complements of finite subsets $F\subset {\rm Val}(E)$. The structure sheaf  is locally defined by the intersection $\bigcap_F R$ of the valuation rings inside $E$. Then the space ${\rm Val}(\K^{\rm ab})$ is the projective limit of the schemes ${\rm Val}(E)$, $E\subset \K^{\rm ab}$. \vspace{.05in}

For the global field $\K=\Q$ of rational numbers, one can consider its maximal abelian extension  $\qcy$ as an abstract  field (\cf Definition \ref{acf}) and try to follow a similar idea. In Section \ref{curve}, we will see however that the space ${\rm Val}(\qcy)$ provides only a rough analogue,  in characteristic zero, of  ${\rm Val}(\K^{\rm ab})$. This section develops the preliminary step of presenting 5 different but equivalent descriptions of the space ${\rm Val}_p(\qcy)$ of extensions of the $p$-adic valuation of $\Q$ to the abstract cyclotomic field $\qcy$. The field  $\qcy$ is the composite of the field generated by roots of unity of order a $p$-power and the field $\qcp$ generated by the  roots of unity of order prime to $p$. We describe canonical isomorphisms of ${\rm Val}_p(\qcy)$ with\footnote{while we believe these results may be known we
give the complete proofs for completeness}:\vspace{.05in}

$(1)$~The space of sequences  of irreducible polynomials $P_n(T)\in\F_p[T]$, $n\in \N$,  fulfilling the basic conditions of the Conway polynomials (\cf Theorem \ref{conway0}).\vspace{.05in}

   $(2)$~The space $\Sigma_p$ of bijections of the monoid $\tilmup=\emup\cup \{0\}$ of roots of unity of order prime to $p$ which commute with their conjugates, as in Definition \ref{defnsp} (\cf Proposition \ref{hope}).\vspace{.05in}

    $(3)$~The space $\Hom(\qfr, \Q_p)$ of field homomorphisms,    where $\qfr\subset \qcp$ is the fixed field under the Frobenius automorphism (\cf Proposition \ref{homfield}).\vspace{.05in}

 $(4)$~The quotient of the space $X_p$ of Definition \ref{defnxp} by the action of $\Gal(\bar\F_p)$ (\cf Proposition \ref{try}).\vspace{.05in}

 $(5)$~The algebraic spectrum of the quotient algebra $\F_p[\mup]/J_p$, where $J_p$ is the reduction modulo $p$ of the cyclotomic ideal (\cf Definition \ref{acf} and Proposition \ref{cpcp1}).\vspace{.05in}

Incidentally, we notice that $(1)$ describes the link between  ${\rm Val}_p(\qcy)$  and the explicit construction of an algebraic closure $\bar\F_p$ of $\F_p$,  by means of a sequence  of irreducible polynomials  over $\F_p$, fulfilling the basic conditions of the Conway polynomials\footnote{Conway polynomials provide a particular  example of such sequence, they are selected  using a lexicographic ordering}. Theorem \ref{conway0} states that the map which associates to a valuation $v\in{\rm Val}_p(\qcy)$ the sequence $\{P_n\}$ of  characteristic polynomials for the action (by multiplication) of the primitive root $\xi_{\frac{1}{p^n-1}}\in \qcp$ on the residue field of the restriction of $v$  to $\qcp$, determines a bijection between ${\rm Val}_p(\qcy)$ and sequences of polynomials in $\F_p[T]$ fulfilling the basic conditions of the Conway polynomials.

\begin{defn}\label{acf} The abstract cyclotomic field $\qcy$ is the quotient of the group ring $\Q[\Q/\Z]$ by the ideal $J$ generated by the idempotents
\begin{equation}\label{enidem}
    \pi_n=\frac 1n \sum_{j=0}^{n-1}e(\frac jn),\qquad n\ge 2.
\end{equation}
\end{defn}

 In general, if we let
\begin{equation}\label{sigmakx}
\sigma_k(x)= \sum_{j=0}^{k-1}\,x^j ,
\end{equation}
then one knows that the $n$-th cyclotomic
polynomial $\Phi_n(x)$ is the gcd~of the polynomials
$\sigma_m(x^d)$, for $m>1$, $m \vert n$ and $d=n/m$. For
$x=e(1/n)$, and $n=md$ one has
\[
\sigma_m(x^d)=\sum_{j=0}^{m-1}\,e(j/m)=m\,\pi_m\in J
\]
thus $\Phi_n(e(1/n))\in J$. It follows that the homomorphism
\begin{equation}\label{isocyc}
    \rho_0:\Q[\Q/\Z]/J\to \C\,, \ \ \rho_0(e(\gamma))=e^{2\pi i\gamma}
\end{equation}
induces an isomorphism of $\qcy$ with the subfield of $\C$ generated by roots of unity.\vspace{.05in}

Using the identification $\Q/\Z=\A_{\Q}^f/\hat\Z$ the group $\hat\Z^*$ acts by automorphisms of $\Q/\Z$ and hence by automorphisms of the group ring $\Q[\Q/\Z]$. This action preserves globally the $n$-torsion in $\Q/\Z$ and hence fixes each of the projection $\pi_n$. It follows that it leaves the ideal $J$ globally invariant and hence it induces an action on the quotient field $\qcy$. This action gives the Galois group $G={\rm Gal}( \qcy:\Q)\simeq \hat\Z^*$ which acts on roots of unity as it acts on $\Q/\Z$. For each prime $p$, one has ($\ell=$ rational prime)
\begin{equation}\label{ggp}
G=\prod_{\ell}\Z_\ell^*=\Z_p^*\times\prod_{\ell\neq p}\Z_\ell^*=\Z_p^*\times G_p.
\end{equation}
One lifts $\Z_p^*$ to the subgroup $\Z_p^*\times 1\subset G$, with all components equal to $1$ except at $p$. This subgroup acts trivially on $\mup$. Its fixed subfield $\qcp\subset \qcy$ is the subfield of $\qcy$ generated over $\Q$ by the group $\emup\subset \qcy$ of roots of unity of order prime to $p$. It coincides with the {\em inertia subfield}
\begin{equation}\label{inertia}
    \qcy\cap \Q_p^{\rm ur}\subset \qcy
\end{equation}
for any extension
 $v\in{\rm Val}_p(\qcy)$ of the $p$-adic valuation to $\qcy$. More precisely let $(\qcy)_v$ be the completion of $\qcy$ for the valuation $v$. Then one knows that the composite subfield $\Q_p\cdot \qcy\subset(\qcy)_v$  is the maximal abelian extension $\Q_p^{\rm ab}$ of $\Q_p$.  This extension is the composite (\cf \cite{Serre})
\begin{equation}\label{composite}
    \Q_p^{\rm ab}=\Q_p^{\rm ur}\cdot \Q_{p^\infty}
\end{equation}
where $\Q_p^{\rm ur}$ denotes the maximal unramified extension of $\Q_p$ and $\Q_{p^\infty}$ is obtained by adjoining to $\Q_p$ all roots of unity of order a $p$-power. The translation Theorem of Galois theory gives a canonical isomorphism (by restriction) of  Galois groups
\begin{equation}\label{isogal}
    \Gal(\Q_p^{\rm ab}:\Q_p)\stackrel{\sim}{\to} \Gal(\qcy:\qcy\cap\Q_p), \qquad \alpha\mapsto \alpha|_\qcy.
\end{equation}
The {\em decomposition subfield}: $\qcy\cap\Q_p$ is independent of the choice of the valuation $v\in{\rm Val}_p(\qcy)$ since $G$ is abelian and acts transitively on ${\rm Val}_p(\qcy)$, more precisely one has the following classical result
\begin{prop}\label{Ggalois}
$(1)$~The  group $G_p=\prod_{\ell\neq p}\Z_\ell^*$ is the group of automorphisms of the
group $\mup$.\vspace{.05in}

$(2)$~The inertia subfield $\qcp$ is the fixed subfield of $\Z_p^*\subset G$ and its Galois group is canonically isomorphic to $G_p$ acting on $\emup\subset\qcp$ as it acts on $\mup$.\vspace{.05in}

$(3)$~Let $f_p\in G_p$ be the element of $G_p=\prod_{\ell\neq p}\Z_\ell^*$ with all components equal to $p$. Then the associated automorphism $\fr\in \Aut(\qcp)$ is the unique automorphism which acts by $x\mapsto x^p$ on the multiplicative group $\emup\subset \qcp$.\vspace{.03in}

$(4)$~The fixed subfield $\qfr\subset \qcp$ of $\fr$ is the decomposition subfield $\qcy\cap\Q_p$.\vspace{.05in}

$(5)$~The  group $G={\rm Gal}( \qcy:\Q)$ acts transitively on ${\rm Val}_p(\qcy)$
 with isotropy $\Z_p^*\times f_p^{\hat
\Z}$, where $f_p^{\hat
\Z}\subset G_p$ is the closure of $f_p^\Z$.
\end{prop}
\proof $(1)$~Let $\Gamma=\mup$ viewed as a
discrete group. The Pontrjagin dual $\hat \Gamma$ is the product $\prod_{\ell\neq
p}\Z_\ell$. We claim that the group of automorphisms of $\Gamma$ is
\begin{equation}\label{autc}
\Aut(\Gamma)=\prod_{\ell\neq p}\Z_\ell^*\,.
\end{equation}
Indeed, one has $\Gamma=\prod^{res}_{\ell\neq p}\Q_\ell/\Z_\ell$, so that the dual of
$\Gamma$ is  $\prod_{\ell\neq p}\Z_\ell$. This  is
a compact ring which contains $\Z$ as a dense subring. Thus an automorphism $\theta$ of the
additive group is characterized by the assignment $a=\theta(1)$ and is given by multiplication
by $a$. Invertibility shows that $a\in\prod_{\ell\neq p}\Z_\ell^*$. This proves \eqref{autc}.

$(2)$~Under the isomorphism \eqref{isogal} the Galois group ${\rm Gal}( \Q_{p^\infty}:\Q_p)\simeq \Z_p^*$ becomes the subgroup $\Z_p^*\times 1\subset G$. The fixed subfield of this subgroup is $\qcp\subset \qcy$ and is the inertia subfield of $\qcy$. The quotient $G/\Z_p^*$ is canonically isomorphic to $G_p$.

$(3)$~Under the isomorphism $\Gal(\qcp:\Q)=G_p$ the action of $\fr$ on $\emup$ corresponds to the
multiplication by $p$ in $\mup$.

$(4)$~The Galois group ${\rm Gal}( \Q_p^{\rm ur}:\Q_p)\simeq \hat\Z$ is topologically generated by the Frobenius automorphism $\fr_p$  whose action on the roots of unity of order prime to $p$ is given by $\fr_p(\xi)= \xi^p$. Under the isomorphism \eqref{isogal} this automorphism restricts to the automorphism $\fr\in \Aut(\qcp)$.  Notice that the fields $\qcp$ and $\Q_p$ are linearly disjoint over their intersection
\begin{equation}\label{interiso}
   K= \qcp\cap \Q_p =\qcy \cap \Q_p\,.
\end{equation}
Then, the translation theorem in Galois theory shows that, by restriction to $\qcp$, one has an isomorphism
\begin{equation*}
    \Gal(\Q_p^{\rm ur}:\Q_p)\stackrel{\sim}{\to} \Gal(\qcp:K), \qquad \fr_p\mapsto\fr.
\end{equation*}
This shows that $K$ is the fixed subfield $\qfr\subset \qcp$ of $\fr$.

$(5)$~It is well known that the Galois group acts transitively on extensions of a valuation. Moreover the isotropy subgroup is the subgroup of the Galois group corresponding to the decomposition subfield and is hence given by $\Z_p^*\times f_p^{\hat
\Z}$. \endproof

\begin{cor}\label{restrict}
The natural map ${\rm Val}_p(\qcy)\to {\rm Val}_p(\qcp)$ given by restriction of valuations is equivariant and bijective.
\end{cor}
\proof The restriction map is equivariant for the action of $G$ on both spaces, these actions are transitive and have the same isotropy group so the restriction map is bijective.\endproof

In fact it is worth giving explicitly the unique extension of a valuation $v\in {\rm Val}_p(\qcp)$ to $\qcy$. The latter field is obtained by adjoining to $\qcp$ primitive roots of unity of order a power of $p$, \ie a solution $z$ of an equation of the form
\[
z^{(p-1)p^{m-1}}+z^{(p-2)p^{m-1}}+\ldots +1=0.
\]
One writes $z=1+\pi$ and finds that the equation fulfilled by $\pi$ is of Eisenstein type, the constant term being equal to $p$,  and reduces to $\pi^{\varphi(n)}=0$,  modulo $p$. This shows that
\[
v(\pi)=\frac{v(p)}{\varphi(n)}\,, \ \ \varphi(n)=(p-1)p^{m-1}.
\]
Then the valuation $v$, normalized so that $v(p)=1$, extends uniquely to elements of the extension $\qcp[z]$ by setting
\begin{equation}\label{valext1}
v(a_0+a_1\pi+\ldots+a_{\varphi(n)-1}\pi^{\varphi(n)-1})=\inf\{v(a_j)+\frac{j}{\varphi(n)}\}.
\end{equation}

\begin{rem}\label{infinite} {\rm
The decomposition subfield $\Q_p\cap \qcy$ is an infinite extension of $\Q$ which contains for instance $\sqrt n$ for $n$ a quadratic residue modulo $p$. Its Galois group $\Gal(\Q_p\cap \qcy:\Q)$ is
the quotient of $G_p$ by the closure of the group of powers of $f_p$ and is a compact group which contains for each prime $\ell\neq p$ the cyclic group of order $\ell-1$ coming from the torsion part of $\Z_\ell^*$.
}\end{rem}

\begin{defn} \label{defnsp} Let $\tilmup=\{0 \} \cup \emup$ be the monoid obtained by adjoining a zero element to the multiplicative group $\emup$. We denote by $\Sigma_p$  the set of bijections $s:\tilmup\to \tilmup$  which
commute with all their conjugates $R\circ s\circ R^{-1}$ under rotations $R$ by elements of $\emup$,
and fulfill the relations: $s(0)=1$, $s^p=s\circ s\circ \ldots\circ s=id$.
\end{defn}
The maps $s$ encode the addition of $1$ on $\tilmup$, when one enriches the multiplicative structure of the monoid $\tilmup$ with an additive structure turning it to a field of characteristic $p$ (\ie an algebraic closure of $\F_p$). Notice that using distributivity the addition of $1$ encodes the full additive structure (\cf \cite{jamifine}).

\begin{lem}\label{Gaut}
The  group $G_p=\prod_{\ell\neq p}\Z_\ell^*$  acts transitively on $\Sigma_p$ with isotropy $f_p^{\hat
\Z}\subset G_p$.
\end{lem}

\proof
We check that $G_p$ acts transitively on $\Sigma_p$. Let $s_j\in \Sigma_p$, for $j=1,2$ and let $\K(s_j)$ be the two corresponding field structures on  $\tilmup$. Then
the two fields $\K(s_j)$ are algebraic closures of $\F_p$ and hence they are
isomorphic. We let $\theta:\K(s_1)\to \K(s_2)$ be such an isomorphism. By
construction $\theta$ is an automorphism of the multiplicative group $\emup$ and it
transports the operation $s_1$ of addition of $1$ in $\K(s_1)$ into the operation
$s_2$ of addition of $1$ in $\K(s_2)$. Since the Galois group of $\bar\F_p$ is topologically generated by the Frobenius $x\mapsto x^p$   one gets, using Galois theory, that the isotropy of
any $s\in \Sigma_p$ is the closure of the group of powers of $f_p$, \ie the subgroup $f_p^{\hat
\Z}\subset G_p$.\endproof

We are now ready to state the main result of this section

\begin{thm}\label{conway0} An element $v\in {\rm Val}_p(\qcy)$ is entirely characterized by
a sequence  of polynomials $P_n(T)\in\F_p[T]$ of degree $n\ge 1$,  such that\vspace{.05in}

$\bullet$~each $P_n(T)$ is monic and irreducible.\vspace{.05in}

$\bullet$~$T\in \F_p[T]/(P_n(T))$ is a generator of the multiplicative group of the quotient field.\vspace{.05in}

$\bullet$~For any integer $m|n$ and for $d=(p^n-1)/(p^m-1)$,  $P_m(T^d)$ is a multiple of
  $P_n(T)$.
\end{thm}

\proof The first step in the proof is to construct a natural map ${\rm Val}_p(\qcy)\ni v\mapsto s_v\in \Sigma_p$.
 We know that $\Q_p\subset (\qcy)_v$ and that $\emup\cup\{0\} \subset\qcy$, thus we consider the  valuation ring $\Z_p^{\rm ur}\subset(\qcy)_v$ of $\Q_p^{\rm ur}$. It contains  $\Z_p$ and $\emup$. Note that the ring generated by $\Z$ and $\emup$ is the ring of integers of the subfield $\qcp\subset \qcy$  generated over $\Q$ by $\emup$. One has the diagram of inclusions

\begin{gather}
\label{fieldext3}
 \,\hspace{45pt}\raisetag{57pt} \xymatrix@C=25pt@R=25pt{
&(\qcy)_v &
 \qcy \ar[l] & \\
\bar\F_p &\Z_p^{\rm ur} \ar[u]\ar[l]_{ \epsilon }  & \Z_p^{\rm ur}\cap \qcy\ar[u]^{ }\ar[l]^{  }&\emup\cup\{0\} \ar[l]^{  }\\
\F_p \ar[u]&\Z_p  \ar[l]_{ \epsilon }\ar[u]^-{ }  & \Z_p\cap \qcy\ar[u]^{ }\ar[l]_{  }&\tau(\F_p) \ar[u]^{ }\ar[l]^{  } \\
}\hspace{140pt}
\end{gather}

\bigskip

where $\tau:\F_p\to \Z_p$ is the \te lift. Note that $\tau(\F_p)\subset \Z_p\cap \qcy$ since this lift is formed of roots of unity (of order $p-1$). In the middle line of the above diagram, the composite map $\epsilon$  from
$\emup\cup\{0\}$ to $\bar\F_p$ is an isomorphism of multiplicative monoids. Indeed, the \te lift $\bar\F_p\ni x\mapsto \tau(x)\in \Z_p^{\rm ur}$ gives the inverse map. Since $\bar\F_p$ is a field one can transport its additive structure using
$\epsilon$ and one obtains a unique element $s_v\in \Sigma_p$.

\begin{prop} \label{hope} The map ${\rm Val}_p(\qcy)\ni v\mapsto s_v\in \Sigma_p$ is a bijection and is equivariant for the action of $G_p=\Gal(\qcp:\Q)$.
\end{prop}
\proof The action of $G_p$ on the  subset $\emup$ is the one described in Lemma \ref{Gaut}. This shows that the map $v\mapsto s_v$ is equivariant. Since both spaces
${\rm Val}_p(\qcy)$ and $\Sigma_p$ are homogeneous spaces over $G_p$ with the same isotropy groups $p^{\hat
\Z}\subset G_p$ as follows from  Lemmas \ref{Ggalois} and \ref{Gaut},  the map $v\mapsto s_v$ is bijective. \endproof

We can produce a concrete construction of the valuation $v$ associated to the map $s_v$. One first determines $v$ on the subfield $\qcp\subset \qcy$. It is enough to determine the valuation $v$ on elements of the form
\begin{equation*}
    x=\sum n_j \xi_j\,, \ \ n_j\in \Z\,, \ \ \xi_j\in \emup  \subset\qcy.
\end{equation*}
Let $K=\bar\F_p$ be the algebraic closure of $\F_p$ obtained by endowing the multiplicative monoid $\emup\cup\{0\}$ with the addition associated to $s_v$. One then has
\begin{equation}\label{valcal}
    v(x)=w_p\left(\sum n_j \tau(\xi_j)\right)
\end{equation}
where  $w_p$ is the $p$-adic valuation in the Witt ring $\winf(K)$ and $\tau$ the \te lift. Finally since the field $\qcy$ is the composite of the subfields $\qcp$ and the fixed field of the action of $G_p\subset \hat\Z^*=\Gal(\qcy:\Q)$ which is generated by roots of unity of order a $p$-power, one can use \eqref{valext1} to extend the valuation $v$ uniquely to $\qcy$.\vspace{.05in}

  \begin{figure}
\begin{center}
\includegraphics[scale=0.9]{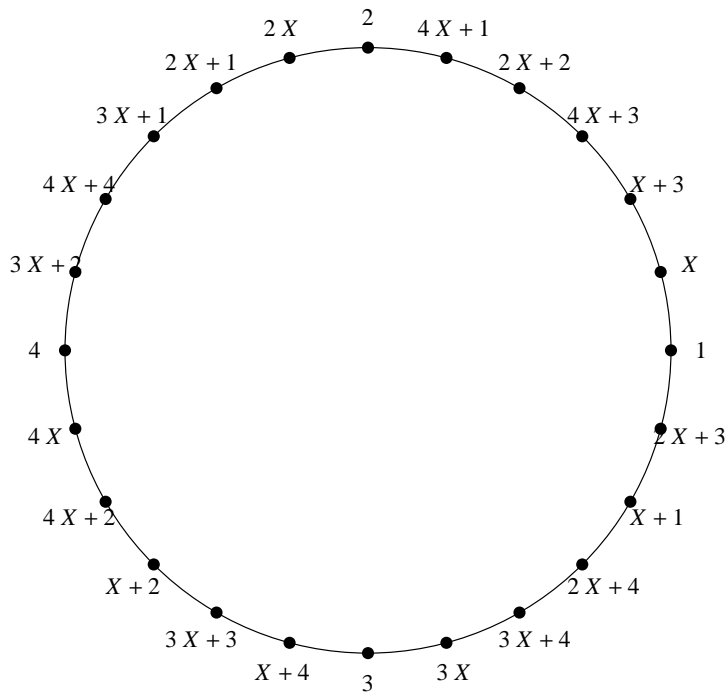}
\end{center}
\caption{The elements of $\F_{25}$ and roots of unity \label{conway} }
\end{figure}

We are now ready to complete the proof of Theorem \ref{conway0}, \ie we prove that:\vspace{.05in}

\begin{lem} \label{hope1}
An element $s\in \Sigma_p$ is entirely characterized by
a sequence $P_n(T)$ of polynomials of $\F_p[T]$  fulfilling the Conway conditions as in
  Theorem \ref{conway0}.\vspace{.05in}
\end{lem}
\proof
Let $s\in \Sigma_p$. For each $n\in \N$, let $\K_n(s)$ be the
corresponding field structure on the union $\{0\}\cup \emup(n)$,  where $\emup(n)$ is the group of  roots of unity of order $p^n-1$ in $\qcy$ generated by $\xi=e(\frac{1}{p^n-1})$. The
$\F_p$ vector space $\K_n(s)$ is of dimension $n$ since its cardinality is
$p^n$. The canonical generator $\xi$ of $\emup(n)$  acts on
the $\F_p$ vector space $\K_n(s)$ by the multiplication $M_\xi$. We let
$P_n(T)$ be its characteristic polynomial \ie the determinant
$P_n(T)=\det(T-M_\xi) $. It is a monic polynomial of degree $n$ with
coefficients in $\F_p$. In the field $\K_n(s)$ one has $P_n(\xi)=0$, since
$M_\xi$ fulfills its characteristic equation. Thus we derive a homomorphism
of algebras $\rho: \F_p[T]/(P_n(T))\to\K_n(s)$ which sends $T\mapsto
\xi$. It is surjective since any non-zero element of $\K_n(s)$ is a power of
$\xi$. Since $P_n(T)$ has degree $n$, the two algebras have the same dimension
over $\F_p$ and thus $\rho$ is an isomorphism. It follows that $P_n(T)$ is
irreducible over $\F_p$. The second property of $P_n(T)$ also follows, since $\xi$ is a
generator of the multiplicative group. Now let $m|n$ be a divisor of $n$. Then
$r=p^m-1$ divides $k=p^n-1$ and the group $\emup(m)$ is a subgroup of $\emup(n)$.
Thus one has a field inclusion $\K_m(s)\subset\K_n(s)$, where the canonical
generator $\xi_m=e(\frac{1}{p^m-1})$ of $\K_m(s)$ is sent to $\xi_n^d$, where $\xi_n$
is the canonical generator $\xi_n=e(\frac{1}{p^n-1})$ of $\K_n(s)$ and
$d=(p^n-1)/(p^m-1)$. One has $P_m(\xi_m)=0$ and hence $P_m(\xi_n^d)=0$ so that, using the above isomorphism $\rho$, it follows that the polynomial $P_m(T^d)$ is a multiple of
  $P_n(T)$.

  Conversely, given a sequence $P_n(T)$ of polynomials fulfilling the   conditions of the theorem, one constructs an algebraic closure $\bar\F_p$ and an isomorphism
  \[\bar \F_p^*\stackrel{j}{\longrightarrow}\emup\] as follows. One lets
  for each $n$, $\K_n=\F_p[T]/(P_n(T))$ and one gets an inductive system using
  for $m|n$ the field homomorphism which sends the generator $T_m$ of $\K_m$ to
  $T_n^d$, $d=(p^n-1)/(p^m-1)$. The inductive limit $\K=\varinjlim \K_n$ is an
  algebraic closure $\bar\F_p$ of $\F_p$ and the map $T_n\mapsto e^{2\pi i/k}$, $k=p^n-1$,
  defines an isomorphism $j$ of $\bar \F_p^*$ with $\emup$. Note that this construction makes sense also
  for $n=1$ and that the first  polynomial is of degree one and thus picks a specific generator of the multiplicative group of $\F_p$. One checks that
  the sequence of polynomials associated to the pair $(\bar\F_p,j)$ is the
  sequence $P_n(T)$. Thus there is a complete equivalence between elements $s\in
  \Sigma_p$ and sequences of polynomials fulfilling the Conway conditions of the Theorem.
\endproof

To make the above map from $\Sigma_p$ to sequences of polynomials more explicit we introduce the ``trace invariant" of an element  $s\in\Sigma_p$. We continue to denote by $\K_n(s)$ the field structure on the union $\{0\}\cup \emup(n)$,  where $\emup(n)$ is the group of  roots of unity generated by $\xi=e(\frac{1}{p^n-1})$. In particular,  $\K_1(s)$ is a field uniquely isomorphic to  $\F_p$. Let $\eta\in \emup$, then the orbit $\cO=\{\fr^k(\eta)\mid k\in\N\}$ of the map $x\mapsto \fr(x)=x^p$ is a finite set, let $|\cO|$ be its cardinality. Then the following sum
\begin{equation}\label{tramap}
    \tr_s(\cO)=\sum_{\cO}\eta
\end{equation}
computed in any $\K_n(s)$, for $|\cO||n$ is the same and determines an element of $\K_1(s)=\F_p$.
\begin{defn}\label{tracemap} Let $\orb(p)$ be the space of orbits of the map $x\mapsto \fr(x)=x^p$ acting on $\emup$.
Let $s\in \Sigma_p$.  We call the map
\begin{equation}\label{tramap1}
   \tr_s:\orb(p)\to \F_p\,,\ \ \cO \mapsto \tr_s(\cO)
\end{equation}
 the trace invariant of $s$.
\end{defn}
The trace invariant characterizes $s$ as shown by the next Proposition.

\begin{prop}\label{tracemapprop}
Let $s\in \Sigma_p$. Then for each $n\in\N$ the polynomial $P_n(T)\in \F_p[T]$ associated to $s$ by Lemma \ref{hope1} is given by
\begin{equation}\label{symsym}
   P_n(T)=T^n+\sum_{k=1}^{n-1}(-1)^k\sigma_k T^{n-k}
\end{equation}
for
\begin{equation}\label{sigmm}
    \sigma_k=\sum_{\cO\subset D_k}\tr_s(\cO),
\end{equation}
where $D_k\subset \mup$ is the set of fractions $\frac{a}{p^n-1}$ where $1\leq a\leq p^n-1$ and the  digits of $a$ in  base $p$ are all zeros except for $k$ of them which are equal to $1$.
\end{prop}
\proof  In the field $\K_n(s)$ the $n$ roots of the polynomial $P_n(T)$ are  the elements
$e(\frac{p^j}{p^n-1})$, for $j=0,\ldots,n-1$. For each $k=1,\ldots,n$, the set of products of $k$ distinct roots is the set of elements of the form
\begin{equation*}
   e\left( \sum_{j\in Y}\frac{p^j}{p^n-1}\right)\,, \ \ Y\subset\{0,1,\ldots,n-1\}\,, \ |Y|=k.
\end{equation*}
One thus gets that the $k$-th symmetric function $\sigma_k$ of the roots of $P_n(T)$ is given by the sum \eqref{sigmm}, over orbits $\cO$  satisfying the prescribed condition $\cO\subset D_k$.\endproof

We now give a third equivalent description of the space ${\rm Val}_p(\qcy)$. We recall that the decomposition subfield $\Q_p\cap\qcy$ is independent of the choice of $v\in {\rm Val}_p(\qcy)$ and is equal to $\qfr\subset \qcy$.
\begin{prop}\label{homfield} The map
\[
\beta: {\rm Val}_p(\qcy) \to \Hom(\qfr, \Q_p),\qquad \beta(v) = \beta_v:\qfr\subset \Q_p
\]
where the fields inclusion $\beta_v$ derives from \eqref{interiso}, determines a canonical and equivariant isomorphism of sets.
\end{prop}
\proof  Notice that the inclusion $\beta_v:\qfr\subset \Q_p$ depends upon the choice of the valuation $v$. One has
$\beta_v\in \Hom(\qfr, \Q_p)$ and the map $v\mapsto \beta_v$ is equivariant for the action of $G_p/f_p^{\hat\Z}$ on ${\rm Val}_p(\qcy)$ and on the space $\Hom(\qfr, \Q_p)$ by
\begin{equation}\label{actact}
    \Hom(\qfr, \Q_p)\ni\beta\mapsto \beta\circ \gamma \qqq \gamma \in G_p=\Gal(\qcp:\Q).
\end{equation}
Since for both spaces the action of $G_p/f_p^{\hat\Z}$ is free and transitive, it follows that the map $\beta$ is bijective. \endproof

We let $\Z[\mup]$ be the group ring of $\mup$ and let $\fr\in \Aut(\Z[\mup])$ be the Frobenius automorphism given  by
the natural linearization of the group automorphism $\mup\to\mup$, of multiplication by $p$ (\cf Corollary \ref{noncan}). The natural ring homomorphism
\begin{equation}\label{ringhom}
   \delta: \Z[\mup]\to \qcp
\end{equation}
is equivariant for the action of $\fr$, its image is the subring of integers of $\qcp$ while the kernel is described by the intersection $\Z[\mup]\cap J$, where $J$ is the ideal of Definition \ref{acf}.
The $\F_p$-algebra
\begin{equation}\label{fpalg}
    \Z[\mup]\otimes_\Z\F_p=\F_p[\mup]
\end{equation}
is perfect since the group $\mup$ is uniquely $p$-divisible. By restriction to the fixed points of $\fr$ and composition with the residue map $\epsilon:\Z_p\to \F_p$, one obtains the map
\begin{equation}\label{resrez}
    \Hom(\qfr, \Q_p)\to \Hom(\F_p[\mup]^\fr,\F_p),\quad \alpha\mapsto {\rm res}(\alpha)=\epsilon\circ\alpha\circ \delta.
\end{equation}
Note that elements of $\Hom(\F_p[\mup]^\fr,\F_p)$ are finitely supported maps from $\orb(p)$ to $\F_p$,  thus they can be lifted to elements of $\Z[\mup]^\fr$. One derives
\begin{equation*}
    \F_p[\mup]^\fr=\Z[\mup]^\fr\otimes_\Z\F_p.
\end{equation*}
Next, we show that the map res as in \eqref{resrez} is injective.
\begin{prop}\label{homfield1} Let $v\in{\rm Val}_p(\qcy)$. We denote by $s_v\in \Sigma_p$ and $\beta_v:\qfr\to \Q_p$ the corresponding elements as in  Lemma \ref{hope} and Proposition \ref{homfield}. Then the  trace invariant map of $s_v$ has the following description
\begin{equation}\label{svsv}
    \tr_{s_v}={\rm res}(\beta_v).
\end{equation}
The map ${\rm res}$ as in \eqref{resrez} is injective.
\end{prop}
\proof The additive structure $s_v$ on $\tilmup=\{0\}\cup \emup$ is the same as that of the residue field of the completion $\qcp$ for the restriction of $v$. It follows that on each orbit $\cO$ of the action of $\fr$ on $\mup$, the sum $\tr_{s_v}(\cO)$ coincides with the residue
\begin{equation*}
 \epsilon(\beta_v(u))\,, \ \  u= \sum_\cO \xi\in \qfr.
\end{equation*}
Since $u=\delta(w)$ where $w=\sum_\cO \xi\in \Z[\mup]^\fr$ one gets  \eqref{svsv}. Then, it follows from Proposition \ref{tracemapprop}  that the map ${\rm res}$ is injective.\endproof

We now briefly explain how one can reconstruct $\alpha\in\Hom(\qfr, \Q_p)$ from its residue $ {\rm res}(\alpha)$, using the Witt functor $\winf$.
Given  $\varsigma\in \Hom(\F_p[\mup]^\fr,\F_p)$, the Witt functor $\winf$ yields a homomorphism
\begin{equation}\label{winfhelp}
    \winf(\varsigma)\in \Hom(\winf(\F_p[\mup]^\fr),\Z_p).
\end{equation}
If $\varsigma={\rm res}(\alpha)$, one can reconstruct $\alpha$ directly using  $\winf(\varsigma)$. This gives a direct proof of the injectivity of the map ${\rm res}$. Indeed, for an orbit $\cO$ of the action of $\fr$ on $\mup$, the element ($\tau=$ \te lift)
\begin{equation*}
   \nu(\cO)= \sum_\cO\tau(\upsilon)\in \winf(\F_p[\mup])
\end{equation*}
is fixed by the Frobenius, \ie $ \nu(\cO)\in\winf(\F_p[\mup]^\fr)$. One then sees that
\begin{equation}\label{finalfine}
    \alpha\left(\sum_\cO \upsilon\right)=\winf(\varsigma)(\nu(\cO))\,.
\end{equation}

\vspace{.05in}

We end this section by giving the relation between $\Sigma_p={\rm Val}_p(\qcy)$ and the space $X_p$ of all injective group
homomorphisms $\sigma:\bar\F_p^\times \to (\qcy)^\times$ (\cf Definition \ref{defnxp}).

We recall that  the Galois group $\Aut(\bar\F_p)$ is the closure $f_p^{\hat\Z}$ of the group generated by the Frobenius $f_p$.
\begin{prop}\label{try}
Let $\bar\F_p$ be a fixed algebraic closure of $\F_p$. Then
\vspace{.05in}

$(1)$~ $G_p$ acts freely and transitively on $X_p$.\vspace{.05in}

$(2)$~ The quotient of $X_p$ by $f_p^{\hat\Z}$ is isomorphic to $\Sigma_p={\rm Val}_p(\qcy)$.
  \end{prop}

 \proof Let $\sigma\in X_p$. The range of $\sigma$ is the group $\emup$ of all roots of unity in $\qcy$ of order prime to $p$. Thus for a pair $\sigma_j\in X_p$, $j=1,2$, one has $\sigma_1\circ \sigma_2^{-1}\in \Aut(\mup)=G_p$. This proves the first statement.
  For any  isomorphism
$
    \sigma:\bar \F_p^*\stackrel{}{\longrightarrow}(\qcy)^\times
$
of the multiplicative group of the algebraic closure $\bar \F_p$ with the group
$\mup$, the following defines an element $s\in \Sigma_p$,
\begin{equation}\label{rootinv}
s(x)=\sigma(\sigma^{-1}(x)+1)\qqq x\neq -1\,, \ s(-1)=0.
\end{equation}
   All elements of $\Sigma_p$ arise this way. Two pairs $(\bar \F_p,\sigma_j)$, $j=1,2$ whose associated
  $ s_j\in \Sigma_p$ are the same are easily seen to be related by an automorphism $\theta\in \Aut(\bar\F_p)$ \ie
  $
  \sigma_2=\sigma_1\circ \theta$.  The second statement thus follows.\endproof

  Proposition \ref{homfield} suggests a more appropriate equivalent description of $X_p$ using a chosen algebraic closure $\bar\Q_p$ of the $p$-adic field and its completion $\C_p$.
  \begin{cor}\label{cpcp} The map
  \begin{equation}\label{isohom}
    i: X_p\to \Hom(\qcp,\C_p),\qquad \sigma\mapsto\tau\circ\sigma^{-1}\
  \end{equation}
where $\sigma^{-1}$ is composed with the \te lift  to determine a field homomorphism from $\qcp$ to $\C_p$, is a bijection  of sets.

  The canonical surjection  $X_p\to \Sigma_p$ of Proposition \ref{try} $(2)$ is the restriction map
  \begin{equation}\label{isohom1}
    \Hom(\qcp,\C_p)\to \Hom(\qfr,\Q_p).
  \end{equation}
  \end{cor}
  \proof Let $\sigma\in X_p$, then $\sigma^{-1}:\mup\to \bar\F_p^\times$ composed with the \te lift $\tau: \bar\F_p^\times\to \zun\subset \C_p$ extends to a unique homomorphism $i(\sigma)\in \Hom(\qcp,\C_p)$.
The  map $i$ is equivariant for the action of $G_p$ on $X_p$ as in Proposition \ref{try} and on
$\Hom(\qcp,\C_p)$ by composition with elements of $G_p=\Gal(\qcp:\Q)$.
Since both actions are free and transitive, $i$ is bijective.

For any $\gamma\in \Hom(\qcp,\C_p)$, the range of $\gamma$ is the subfield of the maximal unramified extension $\Q_p^{\rm ur}\subset \C_p$ generated over $\Q$ by roots of unity of order prime to $p$. One has by construction $\gamma\circ \fr=\fr_p\circ \gamma$. Thus the image $\gamma(\qfr)$ is contained in the fixed subfield $\Q_p$ for the action of $\fr_p$ on $\Q_p^{\rm ur}$. This shows that the restriction map \eqref{isohom1} is well defined. For $j=1,2$, let $\gamma_j\in \Hom(\qcp,\C_p)$, then $\gamma_2^{-1}\circ \gamma_1\in \Gal(\qcp:\Q)$ and this automorphism fixes $\qfr$ pointwise if and only the restrictions $\gamma_j|_\qfr$ are equal. Since $\Gal(\qcp:\qfr)$ is topologically generated by $\fr$, this happens if and only if the $\gamma_j$ are the same in the quotient of $X_p$ by $f_p^{\hat\Z}$.\endproof

We implement the homomorphism $\delta: \Z[\mup]\to \qcp$ of \eqref{ringhom} to associate to an element $\rho\in\Hom(\qcp,\C_p)$ its residue
\begin{equation}\label{resrho}
    {\rm res}(\rho)=\epsilon\circ \rho\circ \delta\in \Hom(\F_p[\mup],\bar\F_p).
\end{equation}
The image of $\delta$ is the ring of integers of $\qcp$, thus the image of $\rho\circ \delta$ in $\C_p$ is contained in $\zun$ and the composite $\epsilon\circ \rho\circ \delta$ is well defined.
Moreover, since  $\text{Ker}(\delta)=J\cap \Z[\mup]$, it follows that $\text{Ker}({\rm res}(\rho))$ contains the ideal $J_p$ reduction of $\text{Ker}(\delta)$ modulo $p$.

\begin{prop}\label{cpcp1} Let $\fpmp$ be the quotient algebra $\F_p[\mup]/J_p$. Then\vspace{.05in}

$(1)$~The map
\begin{equation}\label{resresmap}
    {\rm res}:\Hom(\qcp,\C_p)\to \Hom(\fpmp,\bar\F_p),\qquad {\rm res}(\rho)=\epsilon\circ \rho\circ \delta
\end{equation}
is a bijection of sets.\vspace{.05in}

$(2)$~The algebraic spectrum $\Spec(\fpmp)$ is in canonical  bijection with the set $\Sigma_p$.\vspace{.05in}

$(3)$~ The canonical surjection  $X_p\to \Sigma_p$ of Proposition \ref{try} $(2)$ corresponds to
the natural map
\begin{equation}\label{homspec}
    \Hom(\fpmp,\bar\F_p)\to \Spec(\fpmp).
\end{equation}
\end{prop}
\proof $(1)$~For any integer $m$ prime to $p$, the ideal $J_p$ contains  the projection (\cf Definition \ref{acf})
 $\pi_m=\frac 1m \sum_{j=0}^{m-1}e(\frac jm)$.
 Thus an element $\rho\in\Hom(\fpmp,\bar\F_p)$ is given by a group homomorphism $\rho:\mup\to\bar\F_p^\times$ such that (for $m>1$ prime to $p$)
 $\sum_{j=0}^{m-1}\rho(e(\frac jm))=0$. Notice that this equality holds if and only if $\rho$ is injective and hence, by restriction to the finite level subgroups in the projective limit $\mup$, if and only if it is bijective. Thus $(1)$ follows from the first statement of Corollary \ref{cpcp}.

$(2)$~Consider the finite field $\F_{p^n}$. Two generators of the multiplicative group $\F_{p^n}^\times$ have the same characteristic polynomial if and only if they are conjugate under the action of the Galois group $\text{Gal}(\F_{p^n}:\F_p)$. This shows that the cardinality of the set $I_n$ of irreducible primitive polynomials of degree $n$ over $\F_p$ is
$\varphi(p^n-1)/n$, where $\varphi$ is the Euler totient function. Each of these polynomials $P(X)$ divides the reduction modulo $p$ of the cyclotomic polynomial $\Phi_{p^n-1}(X)$,  thus one derives, modulo $p$, the following equality
\begin{equation}\label{redcyclopol}
    \Phi_{p^n-1}(X)=\prod_{I_n} P(X)
\end{equation}
since the degrees of the polynomials are the same and the right hand side divides the left one. Moreover one also has
\begin{equation*}
\Phi_{p^n-1}\left(e(\frac{1}{p^n-1})\right)\in J_p\,.
\end{equation*}
This determines a canonical isomorphism
\begin{equation}\label{caniso1}
\fpmp_n=\F_p[\emup(n)]/(J_p\cap \F_p[\emup(n)])\to\prod_{I_n}\F_{p^n}
\end{equation}
and thus a canonical bijection of sets $\Spec(\fpmp_n)\to I_n$. Since $\fpmp$ is the inductive limit of the $\fpmp_n$,  $\Spec(\fpmp)$ is the projective limit of the $I_n$, \ie the space of sequences of Conway polynomials as in Theorem \ref{conway0}. This space is in canonical bijection with $\Sigma_p$.

$(3)$~follows from the  proof of $(2)$. \endproof

The restriction to the fixed points of the Frobenius automorphism $\fr\in\Aut(\fpmp)$ does not change the algebraic spectrum as a set, thus we derive a canonical bijection of sets
\begin{equation}\label{bijspec}
    \Spec(\fpmp)\stackrel{\sim}{\to} \Spec(\fpmp^\fr).
\end{equation}

Finally, we characterize the image of the map ${\rm res}$ as in \eqref{resrez}.
\begin{cor}\label{imageres} Let  $\varsigma\in \Hom(\F_p[\mup]^\fr,\F_p)$. Then $\varsigma$ belongs to the image of the map ${\rm res}$ as in \eqref{resrez} if and only if $\text{Ker}(\varsigma)$ contains $\F_p[\mup]^\fr\cap J_p$.
\end{cor}
\proof By \eqref{bijspec}, and Proposition \ref{cpcp1}, $(2)$, one has natural bijections of sets
\begin{equation*}
    \Sigma_p\simeq \Spec(\fpmp)\simeq \Spec(\fpmp^\fr)\simeq \Hom(\fpmp^\fr,\F_p).
\end{equation*}
Then, the statement follows by noticing that the elements of $\Hom(\fpmp^\fr,\F_p)$ are the elements of $\Hom(\F_p[\mup]^\fr,\F_p)$ whose kernel contains $\F_p[\mup]^\fr\cap J_p$.\endproof

 \section{The base point problem and the ``curve" for the global field $\Q$}\label{curve}

In this section we compare the space ${\rm Val}_p(\qcy)$ of extensions of the $p$-adic valuation to $\qcy$ (studied at length in section \ref{sectval}), with the fiber over a prime $p$ of a space $Y$ which represents, in this set-up, the analogue  of the curve that, for function fields, plays a fundamental role in A. Weil's proof of the Riemann Hypothesis. Our results show that for each place $v\in \Sigma(\Q)$, there is a natural model $Y_v$ for the fiber over $v$  and an embedding of this model in a noncommutative space $X(\C_v)$ which is a $v$-adic avatar of the ad\`ele class space $\H_\Q=\A_\Q/\Q^*$.

We shall denote by $\K$ a global field. To motivate our constructions we first recall a few relevant facts holding for function fields.

\subsection{Adelic interpretation of the loop groupoid $\Pi_1^{\rm ab}(X)'$ for function fields}\label{subsectprel}

In this subsection we assume that $\K$ is a function field. We let $\F_q\subset \K$ be the field of constants. Let $\bar \K$ be a fixed separable closure  of $\K$ and let
 $\K^{\rm ab}\subset \bar \K$ be the maximal abelian extension of $\K$. We denote by $\bar \F_q$ the algebraic closure of the finite field $\F_q$ inside $\K^{\rm ab}$.

A main result holding for function fields is that for each finite field extension $E$ of $\bar\F_q\otimes_{\F_q}\K$ the space of (discrete) valuations ${\rm Val}(E)$ inherits the structure of  an algebraic, one-dimensional scheme $X_{E}$
  whose non-empty open sets are the complements of the finite subsets and whose structure sheaf is defined by considering the intersection of the valuation rings inside  $E$. More precisely,  ${\rm Val}(E)$ coincides with the set of (closed) points of the unique projective, nonsingular algebraic curve $X_{E}$ with function field $E$.

We recall  (\cf \cite{Hart} Corollary 6.12) that the category of nonsingular, projective algebraic curves and dominant morphisms is equivalent to the category of function fields of dimension one over $\bar \F_q$. Thus, one associates (uniquely) to $\K^{\rm ab}=\varinjlim_E E$  the projective limit $X^{\rm ab}=\varprojlim_E X_{E}$ which is the abelian cover $X^{\rm ab}\to X$ of the non singular projective curve $X$ over $\F_q$ with function field $\K$. By restricting valuations, one also derives a natural projection map
\[
\pi:  X^{\rm ab}={\rm Val}(\K^{\rm ab})\to \Sigma(\K)
\]
onto the space $\Sigma(\K)$ of valuations of $\K$. Thus, in the function field case one derives a geometric interpretation for the natural fibration associated to the space of valuations of  the field extension $\K^{\rm ab}\supset \K$.

In \cite{wagner} we have given an adelic description of the loop groupoid $\Pi_1^{\rm ab}(X)'$ of the abelian cover $X^{\rm ab}\to X$.
We recall that the \ad class space $\A_\K/\K^*$ of any global field $\K$ has a natural structure of hyperring $\ads$ (\cf\cite{wagner}) and that the prime elements $P(\ads)$ of this hyperring determine a groupoid. The units of this groupoid form the set $\Sigma(\K)$ of places of $\K$ and the source and range maps coincide with the map \[
s:P(\ads)\to \Sigma(\K)\] which associates to a prime element of $\ads$ the principal prime ideal of $\ads$ it generates (and thus the associated place).
When $\K$ is a function field, the groupoid $P(\ads)$ is canonically isomorphic to the loop groupoid $\Pi_1^{\rm ab}(X)'$ of the abelian cover $X^{\rm ab}\to X$, and the isomorphism is equivariant for the respective actions of the
 abelianized  Weil group $\cW^{\rm ab}$ (\ie the subgroup of elements of ${\rm Gal}( \K^{\rm ab}:\K)$ whose restriction to $\bar \F_q$ is an integral power of the Frobenius), and of the id\`ele class group $C_\K = \A_\K^*/\K^*$.

It follows that, as a group action on a set, the action of $\cW^{\rm ab}$ on ${\rm Val}(\K^{\rm ab})$ is isomorphic to the action of the id\`ele class group $C_\K$ on $P(\ads)$. In other words,
 by choosing a set theoretic section $\xi$ of the projection
\begin{equation}\label{mapp}
    \pi\;:\; {\rm Val}(\K^{\rm ab})\to \Sigma(\K)\,, \ \ \pi(v)=v|_\K,
\end{equation}
one  obtains an equivariant set theoretic bijection $P(\ads)\simeq_\xi {\rm Val}(\K^{\rm ab})$ which depends though, in a crucial manner, on the choice of the base point $\xi(w)$, for each place $w\in \Sigma(\K)$. This dependence  prevents one from transporting the  algebraic geometric structure of $X^{\rm ab}$ onto $P(\ads)$, and it also shows that the adelic space $P(\ads)$ carries only the information on the curve $X^{\rm ab}$ given in terms of a set with a group action. \vspace{.05in}

\subsection{Fiber over a finite place of $\Q$}\label{subsectfinite}

Now, we turn to the global field $\K=\Q$. A natural starting point for the construction of a replacement of the covering $X^{\rm ab}$ in this number field case is  to consider the maximal abelian extension of $\Q$, \ie the cyclotomic field $\qcy$  as  analogue of  $\K^{\rm ab}$. Then, the space ${\rm Val}_p(\qcy)$ of extensions of the $p$-adic valuation to $\qcy$  appears as the first candidate for the analogue of the fiber, over a finite place,  of the abelian cover $X^{\rm ab}\to X$. Thus, the first step is evidently that to compare ${\rm Val}_p(\qcy)$ with the fiber $P_p(\H_\Q)$ of the fibration $s: P(\H_\Q)\to \Sigma(\Q)$ over a rational prime $p\in\Sigma(\Q)$. At the level of sets with group actions,  this process shows that  ${\rm Val}_p(\qcy)$ is not yet the correct fiber. The following discussion indicates that one should consider instead  the total space of a principal bundle, with base ${\rm Val}_p(\qcy)$ and structure group a connected compact solenoid $S$ whose definition is given in Proposition \ref{solenoid}. Then, a natural construction of the fiber is provided by the mapping torus $Y_p$ of the action of the Frobenius on the space $X_p$ of Definition \ref{defnxp}.
\begin{prop}\label{compare}
Let $P_p(\H_\Q)$ be the fiber of the groupoid $P(\H_\Q)$  over a non-archimedean, rational prime $p\in \Sigma_\Q$. Then, the following results hold.\vspace{.05in}

$(1)$~The id\`ele class group $C_\Q=\A^*_\Q/\Q^*$ acts transitively on $P_p(\H_\Q)$. The isotropy group of any element of $P_p(\H_\Q)$ is the cocompact subgroup $W_p=\Q_p^*\subset C_\Q$ of classes of  id\`eles $(j_v)$ such that $j_v=1$, $\forall v\neq p$.

$(2)$~Under the class field theory isomorphism
\begin{equation}\label{clfiso}
    \Gal(\qcy:\Q)\simeq C_\Q/D_\Q\,, \ \ D_\Q=\text{connected component of}\; 1,
\end{equation}
$C_\Q$ acts transitively on ${\rm Val}_p(\qcy)$ and  the isotropy group of any element of ${\rm Val}_p(\qcy)$ is \begin{equation}\label{isot}
I_p=\Z_p^*\times H\times \R_+^*\subset \hat\Z^*\times \R_+^*=C_\Q.
\end{equation}
$H\subset G_p=\prod_{\ell\neq p}\Z_\ell^*$ is the closed subgroup $p^{\hat\Z}\subset G_p$ generated by $p$ in $G_p=\prod_{\ell\neq p}\Z_\ell^*$.
\end{prop}
\proof $(1)$ follows from Theorem 7.10 of \cite{wagner}. $(2)$ follows from Lemma \ref{Ggalois}. \endproof

Notice that if $\K$ is a function field and $v$ is a valuation of $\K^{\rm ab}$ extending the valuation $w$ of $\K$, any $g\in\cW^{\rm ab}\subset {\rm Gal}( \K^{\rm ab}:\K)$ such that $g(v)=v$,
belongs to the local Weil group  $\cW^{\rm ab}_w\subset \cW^{\rm ab}$. This is due to the fact that  the restriction of $g$ to an automorphism of $\bar\F_q$ is an integral power of the Frobenius.

When $\K=\Q$,  the isotropy group of the valuation $v$ is instead  larger than the local Weil group $W_p$. The difference is determined by the presence of the quotient $I_p/W_p$ of the isotropy group $I_p$ by the local Weil group  $W_p=\Q_p^*=\Z_p^*\times (\tilde p)^\Z$. Here, $\tilde p$ is represented by the id\`ele all of whose components are $1$ except at the place $p$ where it is equal to $p^{-1}$. By multiplying with the principal id\`ele $p$, one gets the same class as the element of $\hat\Z^*\times \R_+^*$ which is equal to $p$ everywhere except at the place $p$ where it is equal to $1$. Thus, its image in $G_p=\prod_{\ell\neq p}\Z_\ell^*$  is $p$. The quotient group
\begin{equation}\label{connected}
I_p/W_p=(H\times \R_+^*)/(\tilde p)^\Z\simeq  (\hat \Z\times \R)/\Z=S
\end{equation}
is a compact connected solenoid which is described in the following Proposition \ref{solenoid}. The presence of the connected piece $S$ is due to the fact that the connected component of the identity in the id\`ele class group acts trivially,  at the Galois level, on $\qcy$.
 \begin{prop}\label{solenoid} The group $S$ is compact and connected and is canonically isomorphic to the projective limit of the compact groups $\R/n\Z$, under divisibility of the labels $n$.
\end{prop}
\proof We consider first the factor
\begin{equation*}
S_n=( (\Z/n\Z)\times \R)/\Z
\end{equation*}
of the projective limit $S$, where $\Z$ acts diagonally, \ie by the element $(1,1)$, on $ (\Z/n\Z)\times \R$. One has a natural map $p_n: S_n\to \R/n\Z$ given by
\begin{equation*}
p_n(j,s)=s-j\qqq s\in \R, \ j\in \Z/n\Z,
\end{equation*}
 where one views $\Z/n\Z$ as a subgroup of $\R/n\Z$. The map $p_n$ is an isomorphism of groups. When $n$ divides $m$, the subgroup $m\Z\subset \Z$ is contained in $n\Z\subset\Z$ and this inclusion corresponds to the projection $\Z/m\Z\to \Z/n\Z$. Under the isomorphisms $p_n$, this corresponds to the projection $\R/m\Z\to \R/n\Z$. Thus the projective system defining $S$ is isomorphic to the projective system of the projections $\R/m\Z\to \R/n\Z$ and the projective limits are isomorphic.\endproof

Next, we describe a general construction of mapping torus which yields,  when applied to the groups
\begin{equation}\label{special}
    X=G_p\,, \ \ Z=G_p/p^{\hat
\Z},
\end{equation}
the fiber $P_p(\H_\Q)$ of the groupoid $P(\H_\Q)$  over a finite, rational prime $p\in \Sigma(\Q)$.

\begin{prop}\label{fibration} Let $G_p=\prod_{\ell\neq p}\Z_\ell^*$ be the group of automorphisms of the multiplicative
group $\emup$ of roots of unity in $\qcy$ of order prime to $p$ and let $f_p\in G_p$ be the element $\xi\mapsto \xi^p$. Let $G_p$ act freely and transitively on a compact space $X$. Let $Y$ be the quotient space
\begin{equation}\label{quot}
    Y=\left(X\times (0,1)\right)/\sigma^\Z,
\end{equation}
where $\sigma^\Z$ acts on the product $X\times (0,1)$ as follows
\begin{equation}\label{quot1}
    \sigma(x,\rho)=(f_p x,\rho^p)\qqq x\in X,\ \rho\in (0,1).
\end{equation}
Then, the following results hold.\vspace{.05in}

$(1)$~The space $Y$ is compact  and is  an $S$-principal bundle over the quotient $Z$ of $X$ by $f_p^{\hat
\Z}\subset G_p$, where $S$ is the solenoid group of Proposition \ref{solenoid}.\vspace{.05in}

$(2)$~Let $X$ and $Z$ be as in \eqref{special}, then $Y$ is canonically isomorphic to the fiber $P_p(\H_\Q)$.
\end{prop}
\proof $(1)$~We first look at the action of $\Z$ on the open interval $(0,1)$ given by $\rho\mapsto \rho^p$. We consider the map $\psi:(0,1)\to \R$ given by
\begin{equation}\label{psimap}
    \psi(\rho)=\log(-\log(\rho))\qqq \rho\in (0,1).
\end{equation}
One has
\begin{equation*}
\psi(\rho^p)=\log(-\log(\rho^p))=\log(-p\log(\rho))=\log(-\log(\rho))+\log(p)=\psi(\rho)+\log(p)
\end{equation*}
which shows that the action of $\Z$ on $(0,1)$ given by $\rho\mapsto \rho^p$ is isomorphic to the action of $\Z$ on $\R$ given by translation by $\log(p)$.

By construction $G_p=\prod_{\ell\neq p}\Z_\ell^*$ is a compact, totally disconnected group. Next, we show that the map which associates to $n\in \Z$ the element $f_p^n\in G_p$ extends to a bijection of $\hat\Z$ with the closed subgroup of $G_p$ generated by $f_p$. In fact, the isomorphism follows from the isomorphism between $G_p$ and $\text{Gal}(\bar\F_p:\F_p)$, with $f_p$ being the Frobenius. The result follows by applying \eg  \cite{Bourbaki} (Chapitre V, Appendice II, Exercice 5). This gives a natural inclusion $\hat\Z\subset G_p$, $a\mapsto f_p^a$, as a closed subgroup. We now consider the action of the product group $\hat\Z\times \R_+^*$ on $X\times (0,1)$ given by
\begin{equation}\label{prodact}
    (a,\lambda). (x,\rho)=(f_p^a\,x,\rho^\lambda).
\end{equation}
By construction the element $(1,p)\in \hat\Z\times \R_+^*$ acts as $\sigma$ (\cf \eqref{quot1}).
The quotient group
\begin{equation}\label{quogr}
    (\hat\Z\times \R_+^*)/s^\Z\,, \ \ s=(1,p)
\end{equation}
is isomorphic to the solenoid $S$, by using the isomorphism of the group $\R_+^*$ with $\R$ given by the logarithm in base $p$. To see that $Y$ is a principal bundle over $S$ one uses the map $\psi$ of \eqref{psimap} to check that $S$ acts freely on $Y$. The quotient of $Y$ by the action of $S$ is the quotient of $X$ by the action of $\hat\Z$.

$(2)$~The fiber $P_p(\H_\Q)$  has a canonical base point given by the idempotent $u\in P_p(\H_\Q), u^2=u$. Hence by applying Proposition \ref{compare}, this fiber is canonically isomorphic to the quotient $C_\Q/W_p$. By identifying
$C_\Q$ with $\hat\Z^*\times  \R_+^*$, this quotient coincides with the quotient of $G_p\times \R_+^*$ by the powers of the element $(p,p)\in G_p\times \R_+^*$. Under the bijection $\rho\mapsto -\log(\rho)$ from $(0,1)$ to $\R_+^*$, one obtains the same action as in \eqref{quot1} and hence the desired isomorphism.
\endproof

In order  to obtain the analogue, for the global field $\K=\Q$, of the fiber  of the algebraic curve $X^{\rm ab}$, we should  apply the construction of Proposition \ref{fibration}  to
a compact space $X_p$  so that the following requirements are satisfied\vspace{.05in}

$(1)$~$G_p$ acts freely and transitively on $X_p$\vspace{.05in}

$(2)$~The quotient of $X_p$ by $f_p^{\hat\Z}$ is canonically isomorphic to ${\rm Val}_p(\qcy)$.\vspace{.05in}

 Proposition \ref{try} provides a natural candidate for $X_p$. Moreover, equation \eqref{isohom} shows that one can equivalently describe  $X_p$ as the space $\Hom(\qcp,\C_p)$ and that the canonical identification of $X_p/f_p^{\hat\Z}$ with ${\rm Val}_p(\qcy)$ is given by the restriction map to the fixed points of $\fr$ as in \eqref{isohom1}. We derive the definition of the following model
 for the fiber $Y_p$ over a finite prime $p$
 \begin{equation}\label{quotmod}
    Y_p=\left(\Hom(\qcp,\C_p)\times (0,1)\right)/\sigma^\Z.
\end{equation}

 \subsection{Fiber over the archimedean place of $\Q$}\label{subsectarch}

We move now to the discussion of the analogues of the spaces ${\rm Val}_p(\qcy)$, $X_p$ and $Y_p$, when $p$ is the archimedean prime $p=\infty$ (\ie the archimedean valuation). The space ${\rm Val}_\infty(\qcy)$ is the space of multiplicative norms on $\qcy$ whose restriction to $\Q$ is the usual absolute value. For $v\in {\rm Val}_\infty(\qcy)$, the field completion $(\qcy)_v$ is isomorphic to $\C$, thus  one derives
\begin{equation}\label{infty}
    {\rm Val}_\infty(\qcy)=\Hom(\qcy,\C)/\{\pm 1\}
\end{equation}
where $\{\pm 1\}\subset \hat\Z^*=\Gal(\qcy:\Q)$ corresponds to complex conjugation. It follows that for $p=\infty$ the space $X_p$ is simply
\begin{equation}\label{xpinclinfty}
    X_\infty=\Hom(\qcy,\C).
\end{equation}
 On the other hand, the fiber $P_\infty(\H_\Q)$ is the quotient $C_\Q/W_\infty$, where $W_\infty=\R^*$ is the cocompact subgroup of $C_\Q$ given by classes of id\`eles, whose components are all $1$ except at the archimedean place. Then, we derive that
\begin{equation}\label{infty1}
    P_\infty(\H_\Q)=\hat\Z^*/\{\pm 1\}.
\end{equation}
This discussion shows that at $p=\infty$ there is no need for a mapping torus, and that the expected fiber is simply given by \begin{equation}\label{infty2}
   Y_\infty= {\rm Val}_\infty(\qcy)=\Hom(\qcy,\C)/\{\pm 1\}= X_\infty/\{\pm 1\}.
\end{equation}

 \subsection{Ambient noncommutative space}\label{subsectwhync}
The model \eqref{quotmod}
for the fiber over a rational prime $p$ is only a preliminary step toward the global construction of the ``curve'' which we expect to replace, when $\K=\Q$, the geometric cover $X^{\rm ab}$.  In fact, one still needs to suitably combine these
models into a  noncommutative space to account for the presence of transversality factors in the explicit formulas. We explain why in some details below.

In \cite{announc3} we showed how to determine the counting function $N(q)$ (a distribution on $[1,\infty)$) which replaces, for $\K=\Q$, the classical Weil counting function for a field $\K$ of functions of an algebraic curve $Y$ over $\F_p$ (\cf \cite{Manin,Soule}). The Weil counting function determines the number of rational points on the curve $Y$ defined over  field extensions $\F_q$ of $\F_p$
 \[
\#Y(\F_q)=N(q)=q-\sum_{\alpha} \alpha^r+1,\qquad  q=p^r.
\]
The numbers $\alpha$'s are the complex roots  of the characteristic polynomial of the Frobenius endomorphism acting on the \'etale cohomology $H^1(Y\otimes\bar\F_p,\Q_\ell)$, for $\ell\neq p$.
In \cite{announc3} we have  shown that the distribution $N(q)$ associated to the (complete) Riemann zeta function is described by the similar formula
\begin{equation}\label{countingfunct}
 N(u)=u-\frac{d}{du}\left(\sum_{\rho\in Z}{\rm order}(\rho)\frac{u^{\rho+1}}{\rho+1}\right)+1.
\end{equation}
where $Z$ is the set of non trivial zeros of the Riemann zeta function.
This distribution is positive on $(1,\infty)$ and  fulfills all the expected properties of a counting function. In particular, it takes the correct value
$N(1)=-\infty$ in agreement with the (expected) value of the Euler characteristic. In \cite{japanese} we pushed these ideas further and we explained
 how to implement the trace formula understanding of the explicit formulas in number-theory,
to express the distribution $N(q)$ as an {\em intersection number} involving
the scaling action of the id\`ele class group on the \ad class space. This development involves  a Lefschetz formula whose geometric side corresponds to the following expression of the counting distribution $N(u)$
  \begin{equation}\label{Nu}
    N(u)=\frac{d}{du}\varphi(u)+ \kappa(u), \ \    \varphi(u)=\sum_{n<u}n\,\Lambda(n).
\end{equation}
Here, $\Lambda(n)$ is the von-Mangoldt function taking the value $\log p$ at prime powers $p^\ell$  and zero otherwise and $\kappa(u)$ is the distribution
defined, for any test function $f$, as
\begin{equation}\label{kappadu}
\int_1^\infty\kappa(u)f(u)d^*u=\int_1^\infty\frac{u^2f(u)-f(1)}{u^2-1}d^*u+cf(1)\,, \qquad c=\frac12(\log\pi+\gamma)
\end{equation}
where $\gamma=-\Gamma'(1)$ is the Euler constant. The distribution $\kappa(u)$ is positive on $(1,\infty)$ and in this domain it is equal to the function $\kappa(u)=\frac{u^2}{u^2-1}$.
The contribution in the counting distribution $N(u)$ coming from the term $\frac{d}{du}\varphi(u)$ in \eqref{Nu} can be understood geometrically as arising from a counting process performed on the fibers $Y_p$ (each of them accounting for the delta functions located on the powers of $p$). The value $\log(p)$ coming from the von-Mangoldt function $\Lambda(n)$ corresponds to the length of the orbit in the mapping torus (\cf \cite{japanese}, \S2.2).
On the other hand, as explained in \cite{japanese}, the contribution of the archimedean place cannot be understood in a naive manner as a simple counting  process of points and its expression involves a transversality factor  measuring the transversality of the action of the id\`ele class group with respect to periodic orbits. This shows that the periodic orbits cannot be considered in isolation and must be thought as (suitably) embedded in the ambient \ad class space. This development supplies a precious hint toward  the final construction of the ``curve'' and  shows that the role of ergodic theory and noncommutative geometry is indispensable. \vspace{.05in}

\subsection{The BC-system over $\Z$ and $\F_{1^\infty}\otimes_{\F_1}\Z$}\label{subsectendo}
 Next, we shall  explain how  the BC-system over $\Z$  gives, for each $p$, a natural embedding of the fiber $Y_p$ (\cf\eqref{infty2}) into a noncommutative space constructed using the set $\cE(\C_p)$ of the $\C_p$-rational points  of the affine group scheme $\cE$ which describes the abelian part of the system (\cf\cite{ccm}). Since the fields $\C_p$ are abstractly pairwise isomorphic the obtained spaces are also abstractly isomorphic, but in a non canonical manner.
In \cite{ccm}, following a proposal of C. Soul\'e for the meaning of the ring $\F_{1^n}\otimes_{\F_1}\Z$, we noted that the inductive limit
\begin{equation}\label{missext0}
\F_{1^\infty}\otimes_{\F_1}\Z:=\varinjlim_n\F_{1^n}\otimes_{\F_1}\Z=\Z[\Q/\Z]
\end{equation}
coincides with the abelian part of the algebra defining the integral BC-system. The description given in that paper of the BC-system as an affine pro-group scheme $\cE$ over $\Z$ together with the
dynamic of the action of a semigroup of endomorphisms, allows one to consider its rational points over any ring $A$
\begin{equation}\label{points}
    \cE(A)=\Hom(\Z[\Q/\Z],A)\,.
\end{equation}
Then, one can implement, for each rational prime $p$, the canonical inclusion
\begin{equation}\label{xpincl}
    X_p=\Hom(\qcp,\C_p)\subset \Hom(\Z[\Q/\Z],\C_p)=\cE(\C_p)\,.
\end{equation}
The next result shows that  the space
\begin{equation}\label{algadclass2}
 X(\C_p):=  \left( \cE(\C_p)\times (0,\infty)\right)/(\N\times \{\pm 1\})
\end{equation}
 matches,  for any $p$ including $p=\infty$,  the definition of the \ad class space $\H_\Q$.
The action of $m=\pm n$ (in the semigroup $\N\times \{\pm 1\}$) is the product of the linearization  of the action $e(\gamma)\mapsto e(m\gamma)$ on the ($\C_p$-rational points of the) scheme $\cE$, with the action on  $(0,\infty)$ given by the map $x\mapsto x^m$.

\begin{prop}\label{caseC} $(1)$~The space $X(\C)$ is canonically isomorphic to the \ad class space $\H_\Q$.

$(2)$~The subspace of the \ad class space made by classes whose archimedean component vanishes corresponds to the quotient
\begin{equation}\label{infty12}
   \cE(\C)/(\N\times \{\pm 1\})=\hat\Z/(\N\times \{\pm 1\}).
\end{equation}
\end{prop}
\proof $(1)$~The space $\cE(\C)$ is the space of complex characters of the abelian group $\Q/\Z$ and is canonically isomorphic to $\hat\Z$. We use the map $\rho\mapsto -\log(\rho)$ to map the interval $(0,\infty)$ to $\R$. Under this map the transformation $x\mapsto x^m$ becomes the multiplication by $m$. The action $e(\gamma)\mapsto e(m\gamma)$ on the scheme $\cE$ corresponds to the multiplication by $m$ in $\hat\Z$. Since any \ad class is equivalent to an element of $\hat\Z\times \R$, \eqref{algadclass2} gives, for $p=\infty$
\begin{equation}\label{adclassinfty}
   X(\C)=\left(\hat\Z\times \R\right)/(\N\times \{\pm 1\})=\A_\Q/\Q^*=\H_\Q.
\end{equation}
$(2)$~follows from the identification \eqref{adclassinfty}.\endproof

 Note that by using the inclusion $(0,1)\subset (0,\infty)$, one derives a natural inclusion
\begin{equation*}
    Y_p=\left(\Hom(\qcp,\C_p)\times (0,1)\right)/\sigma^\Z\to
    \left( \cE(\C_p)\times (0,\infty)\right)/(\N\times \{\pm 1\})=X(\C_p).
\end{equation*}
For $p=\infty$ one has the natural inclusion
\begin{equation}\label{natmapyinfty}
Y_\infty=\Hom(\qcy,\C)/\{\pm 1\}\to
    \left( \cE(\C)\times (0,\infty)\right)/(\N\times \{\pm 1\})=X(\C)
\end{equation}
which is obtained by using the canonical inclusion \eqref{xpincl} for $p=\infty$ and  the fixed point $1\in (0,\infty)$.
\vspace{.05in}

The group ring $\Z[\Q/\Z]$ is a Hopf algebra for the coproduct
\begin{equation}\label{copr}
    \Delta(e(\gamma))=e(\gamma)\otimes e(\gamma)\qqq \gamma\in \Q/\Z
\end{equation}
and the antipode $e(\gamma)\mapsto e(-\gamma)$, thus
 $\cE$ is a group scheme.
\begin{prop} \label{groupscheme} Let $A$ be a commutative ring.\vspace{.05in}

$(1)$~The abelian group $\cE(A)$ is torsion free.\vspace{.05in}

$(2)$~The space
\begin{equation}\label{quotX}
X(A)=\left(\cE(A)\times (0,\infty)\right)/(\N\times \{\pm 1\})
\end{equation}
is a module over the hyperring $\H_\Q$.\vspace{.05in}

$(3)$~For any rational prime $p$, $X(\C_p)$ is a free module of rank one over $\H_\Q$.
\end{prop}
\proof
$(1)$~One has
\begin{equation*}
    \cE(A)=\Hom(\Z[\Q/\Z],A)=\Hom(\Q/\Z,A^\times)
\end{equation*}
where the second $\Hom$ is taken in the category of abelian groups. Since the group $\Q/\Z$ is divisible
the group $\Hom(\Q/\Z,H)$ has no torsion, for any abelian group $H$.

$(2)$~We first show that $X(A)$ is a hypergroup and in fact a vector space over the Krasner hyperfield $\kras=\{0,1\}$ (\cf \cite{wagner}).  The two abelian groups $\cE(A)$ and $(0,\infty)$ are both torsion free, thus one gets
\begin{equation}\label{hyperdefn}
   \left(\cE(A)\times (0,\infty)\right)/(\N\times \{\pm 1\})=
   \left(\left(\cE(A)\times (0,\infty)\right)\otimes_\Z\Q\right)/\Q^\times
\end{equation}
which is a projective space, hence a vector space over $\kras$ (\cf \cite{wagner}). Next we show that $X(A)$ is a module over $\H_\Q$. We use the canonical ring isomorphism $\hat\Z=\End_\Z(\Q/\Z)$ to define the following ring homomorphism from $\hat\Z$ to the ring $\End_\Z(\cE(A))$
\begin{equation}\label{end}
   c_A: \hat\Z\to \End_\Z(\cE(A))\,, \ \ c_A(\alpha)\xi=\xi\circ \alpha \qqq \xi\in \Hom(\Q/\Z,A^\times).
\end{equation}
 The map
 \begin{equation}\label{pospos}
 p: \R\to  \End_\Z(\R_+^*),\quad p(\lambda)x= x^\lambda
 \end{equation}
 is a ring homomorphism, thus $c_A\times p$ defines a ring homomorphism from $\hat\Z\times \R$ to the endomorphisms of the abelian group $\cE(A)\times (0,\infty)$. For any $m\in \Z\subset \hat\Z\times \R$, one has
 \begin{equation}\label{submon}
    (c_A\times p)(m)((e(\gamma),x))=(e(m\gamma), x^m),
 \end{equation}
 thus the restriction of $c_A\times p$ to the monoid of non-zero elements of $\Z$ gives the equivalence relation which defines $X(A)$ as in \eqref{quotX}. It follows an action  of the  hyperring
 \begin{equation*}
    \left((\hat\Z\times \R) \otimes_\Z\Q\right)/\Q^\times=\A_\Q/\Q^\times=\H_\Q
 \end{equation*}
on the hypergroup \eqref{hyperdefn}.

$(3)$~It is easy to see that, once one fixes an embedding $\rho:\qcy\to\C_p$ and an $x\in (0,\infty)$ and a real number $x\neq 1$, the element $(\rho,x)\in X(\C_p)$ is a generator of $X(\C_p)$ as a free module over $\H_\Q$. \endproof

The next result displays some interesting arithmetic-geometric properties of the scheme $\cE$.

\begin{prop} \label{ramified}
$(1)$~Let $\Q_p^{\rm ab}\subset \C_p$ be the maximal abelian extension of $\Q_p$. Then the natural map $\cE(\Q_p^{\rm ab})\to \cE(\C_p)$ is a bijection of sets.\vspace{.05in}

$(2)$~Let $\Q_p^{\rm ur}\subset \Q_p^{\rm ab}$ be the maximal unramified extension of $\Q_p$ and
$\Z_p^{\rm ur}\subset \Q_p^{\rm ur}$ the valuation ring of the $p$-adic valuation. Then the natural map $\cE(\Z_p^{\rm ur})\to \cE(\Q_p^{\rm ur})$ is a bijection of sets.\vspace{.05in}

$(3)$~Let $\epsilon : \Z_p^{\rm ur}\to \bar\F_p$ be the residue homomorphism.  Then the associated map $\cE(\Z_p^{\rm ur})\to \cE(\bar\F_p)$ is a bijection.\vspace{.05in}

$(4)$~The $\H_\Q$-module
\begin{equation*}
   X(\bar\F_p)\simeq X(\Z_p^{\rm ur})\simeq X(\Q_p^{\rm ur})\subset X(\Q_p^{\rm ab})\simeq X(\C_p)
\end{equation*}
is described as
\begin{equation}\label{substuff}
   X(\bar\F_p)=\ffp_p X(\C_p)
\end{equation}
where $\ffp_p\in\Sp(\H_\Q)$ is the prime ideal of ad\`ele classes whose $p$-component vanishes.
\end{prop}

\proof $(1)$~Let $\rho\in \cE(\C_p)=\Hom(\Z[\Q/\Z],\C_p)$. Then the image of $\rho$ is contained in the subfield of $\C_p$ generated over $\Q$ by roots of unity, which is contained in $\Q_p^{\rm ab}$.

$(2)$~Let $\rho\in \cE(\Q_p^{\rm ur})=\Hom(\Z[\Q/\Z],\Q_p^{\rm ur})$. Then the image of $\rho$ is contained in the subring of $\Q_p^{\rm ur}$ generated over $\Z$ by roots of unity  (of order prime to $p$) which is contained in $\Z_p^{\rm ur}$.

$(3)$~Let $\rho\in \cE(\Z_p^{\rm ur})=\Hom(\Z[\Q/\Z],\Z_p^{\rm ur})$. Then $\rho$ is entirely characterized by the group homomorphism
\begin{equation*}
    \rho:\Q/\Z\to G
\end{equation*}
where $G$ is the group of roots of unity in $\Z_p^{\rm ur}$, which is non canonically isomorphic to the group $\emup$ of abstract  roots of unity of order prime to $p$. Similarly an element of
$\cE(\bar\F_p)=\Hom(\Z[\Q/\Z],\bar\F_p)$ is entirely characterized by the associated  group homomorphism from $\Q/\Z$ to $\bar\F_p^*$. Since the residue morphism $\epsilon$ gives an isomorphism $G\stackrel{\simeq}{\to} \bar\F_p^*$ one obtains the conclusion.

$(4)$~One has $\hat\Z=\Hom(\Q/\Z,\Q/\Z)$. Let, as above, $\mup\subset \Q/\Z$ be the  subgroup of elements of denominator prime to $p$. Then the subset
$\Hom(\Q/\Z,\mup)\subset \Hom(\Q/\Z,\Q/\Z)$ is given by
\begin{equation*}
    \{(a_\ell)\in \prod \Z_\ell=\hat\Z \mid \ \ a_p=0\}
\end{equation*}
which corresponds to the prime, principal ideal $\ffp_p$ of the hyperring structure $\H_\Q$ inherent to the \ad class space (\cf \cite{wagner}).\endproof

\section{The standard model of $\bar\F_p$ and the BC-system}\label{sectlens}

As shown in section \ref{sectval}, the space ${\rm Val}_p(\qcy)$ is intimately related to the space of sequences  of irreducible polynomials $P_n(T)\in\F_p[T]$, $n\in \N$,  fulfilling the basic conditions of the Conway polynomials (\cf Theorem \ref{conway0}) and hence to the explicit construction of an algebraic closure of $\F_p$.  The normalization condition using the lexicographic ordering just specifies a particular element $v_c$ of ${\rm Val}_p(\qcy)$. Since the explicit computation of the sequence $P_n(T)\in\F_p[T]$, $n\in \N$, associated to $v_c$ has been proven to be completely untractable, B. de Smit and H. Lenstra have recently devised a more efficient construction of $\bar\F_p$ (\cf \cite{smitlenstra}). Our goal in this section is to make explicit the relation between their construction, the BC-system and the sought for ``curve". \vspace{.05in}

When $\K$ is a global field of positive characteristic \ie the function field  of an algebraic curve over a finite field $\F_q$, the intermediate extension $\K\subset\bar\F_q\otimes_{\F_q}\K\subset \K^{\rm ab}$ plays an important geometric role since it corresponds to working over an algebraically closed field. For $\K=\Q$, it is therefore natural to ask for an intermediate extension $\Q\subset L\subset \qcy$ playing a similar role. One  feature of the former extension is that the residue fields are algebraically closed.

In their construction, de Smit and Lenstra use the intermediate extension $\Q\subset \qsl\subset \qcy$ which comes very close to fulfill the expected properties. For each prime $\ell$, let us denote by $\Delta_\ell\subset \Z_\ell^*$
the torsion subgroup. For $\ell=2$ one has $\Delta_2=\{\pm 1\}$, while for $\ell\neq 2$ one gets $\Delta_\ell=\tau(\F_\ell^*)$, where $\tau:\F_\ell\to \Z_\ell$ is the \te lift. The product
\begin{equation}\label{deltasub}
    \Delta:=\prod_\ell \Delta_\ell\subset \prod_\ell \Z_\ell^*
\end{equation}
is a compact group, and  a subgroup of the Galois group $\hat\Z^*=\Gal(\qcy:\Q)$. By Galois theory, one can thus associate to $\Delta$ a (fixed) field extension
\begin{equation}\label{Lext}
    L=\qsl\subset \qcy.
\end{equation}
Notice that one derives a subsystem of the BC-system given by the fixed points of the action of $\Delta$. At the rational level and by implementing the cyclotomic ideal $J$ of Definition \ref{acf}, one obtains the exact sequence of algebras
\begin{equation}\label{bcsl}
   0\to J\cap \qqsl\to \qqsl\stackrel{q}{\to} \qsl\to 0
\end{equation}
The image of the restriction to $\bcsl$ of the homomorphism $q$  is contained in the integers of $\qsl$ and one has
\begin{equation}\label{galqsl}
    \Gal(\qsl:\Q)\simeq\hat\Z^*/\Delta\simeq\prod_\ell \Z_\ell^*/\Delta_\ell.
\end{equation}
The space ${\rm Val}_p(\qcy)$  is the total space of a principal bundle whose base is the space ${\rm Val}_p(\qsl)$ of valuations on $\qsl$ extending the $p$-adic valuation.  The group of the principal bundle is the quotient of $\Delta$ by its intersection $\Delta_p$ with the isotropy group of elements of ${\rm Val}_p(\qcy)$. The projection ${\rm Val}_p(\qcy)\to {\rm Val}_p(\qsl)$ is given by restriction of valuations from $\qcy$ to $\qsl$.
For $w\in {\rm Val}_p(\qsl)$, the isotropy group $\Pi_p$ of $w$ for the action of $\Gal(\qsl :\Q)$ is the image of the isotropy group of $v$ in $\Gal(\qcy :\Q)$ for any extension $v$ of $w$ to $\qcy$.
It follows from  Lemma \ref{Ggalois} that the  isotropy subgroup of $v$ is $\Z_p^*\times f_p^{\hat\Z}\subset \Z_p^*\times G_p$, thus  one gets
\begin{equation}\label{imiso}
   \Pi_p\simeq \Z_p^*/\Delta_p\times \overline{f_p^\Z}\,, \ \ \overline{f_p^\Z}\subset
\prod_{\ell\neq p}\Z_\ell^*/\Delta_\ell.
\end{equation}
\begin{lem}\label{subopen} For each prime $\ell$ the group $\Z_\ell^*/\Delta_\ell$ is canonically isomorphic to the additive group $\Z_\ell$.
Moreover, for each prime $p\neq \ell$ the closed subgroup of $\Z_\ell^*/\Delta_\ell$ generated by $p$ is open and of finite index $\ell^{u(p,\ell)}$ where
\begin{equation}\label{upl}
    u(p,\ell)=\begin{cases}
                 v_\ell(p^{\ell-1}-1)-1, &\text{for $\ell>2$}\\
                  v_2(p^2-1)-3, &\text{for $\ell=2$.}
                \end{cases}
\end{equation}
\end{lem}
\proof For each prime $\ell$ there is a canonical isomorphism of groups
\begin{equation}\label{caniso}
    \Z_\ell^*\stackrel{\sim}{\to} \Delta_\ell\times \Z_\ell, \quad  x\mapsto (\omega(x),i_\ell(x))
\end{equation}
where the group $\Z_\ell$ is viewed as an additive group. For $\ell$ odd,  $\omega(x)$  is the unique $\ell-1$ root of unity which is congruent to $x$ modulo $\ell$ and  $i_\ell(x)$, as  in \eqref{ipr}, is the ratio $\displaystyle{\log_\ell x/\log_\ell(1+\ell)}$. For $\ell=2$, $\omega(x)=\pm 1$ is congruent to $x$ modulo $4$ and $i_2(x)=\log_2 x/\log_2(1+4)$.
The first statement thus follows. The second statement follows since one has
\begin{equation*}
    v_\ell(i_\ell(p))=u(p,\ell)
\end{equation*}
and
the closed subgroup of $\Z_\ell$ generated by $i_\ell(p)$ is $\ell^{u(p,\ell)}\Z_\ell$.\endproof

Under the isomorphisms
\begin{equation}\label{muladdiso}
   \Gal(\qsl:\Q)\simeq \prod_\ell \Z_\ell^*/\Delta_\ell\simeq\prod_\ell\Z_\ell\simeq\hat\Z
\end{equation}
one gets, by the Chinese remainder theorem, that
\begin{equation}\label{imiso2}
   \Pi_p\simeq \Z_p\times \prod_{\ell\neq p}\ell^{u(p,\ell)}\Z_\ell\subset \hat\Z\,.
\end{equation}
Notice the independence of the places $\ell$ in the above formula which makes the group $\Pi_p$  a cartesian product and allows one to express ${\rm Val}_p(\qsl)$ as an infinite product of finite sets.

To label concretely these finite sets consider,
for each prime $\ell$ the $\Z_\ell$-extension $\B_\infty(\ell)$ of $\Q$. One has $\B_\infty(\ell)=\cup_k \B_k(\ell))$ where, for $k\in \N$, the finite extension $\B_k(\ell)$ of $\Q$ is associated to $\ell^{-k}\in\Q/\Z$ viewed as a character of $\hat\Z\simeq\Gal(\qsl:\Q)$. For $\ell$ odd, $\B_k(\ell)$ is the fixed subfield for the action of $\Delta_\ell$ on the extension of $\Q$ generated by a primitive root of unity of order $\ell^{k+1}$. For $\ell=2$ one uses a primitive root of unity of order $2^{k+2}$.
 We  denote by $\B(\ell,p)=\B_{u(p,\ell)}(\ell)$: this is a cyclic extension of $\Q$ of degree $\ell^{u(p,\ell)}$. The Artin reciprocity law shows that, for  $p$ a prime $p\neq \ell$, the reduction modulo $p$ of the integers of $\B(\ell,p)$,   decomposes as a product of $\ell^{u(p,\ell)}$ copies of $\F_p$, parameterized by the set  ${\rm Val}_p(\B(\ell,p))$ of extensions of the $p$-adic valuation to $\B(\ell,p)$, which  is a finite set of cardinality $\ell^{u(p,\ell)}$.

The following result is  a  consequence of the construction of the ``standard model" of de Smit and Lenstra for the algebraic closure of a finite field

\begin{thm}\label{BCSL} Let $p$ be a rational prime.

$(1)$~For  $\ell\neq p$ a prime,  the restriction map ${\rm Val}_p(\B_\infty(\ell))\to{\rm Val}_p(\B(\ell,p))$ is bijective.

 $(2)$~The restriction maps from $\qsl$ to $\B_\infty(\ell)$ give  a bijection
\begin{equation}\label{prodzpl}
  {\rm Val}_p(\qsl)= \prod_{\ell\neq p} {\rm Val}_p(\B_\infty(\ell)).
\end{equation}

$(3)$~The restriction of $v\in {\rm Val}_p(\qsl)$ to $\qcp_\Delta$ is unramified and the residue field is isomorphic to
\begin{equation}\label{primetop}
    \bigcup_{n \in I(p)} \F_{p^n}\subset \bar\F_p
\end{equation}
where $I(p)\subset\N$ denotes the subset of positive integers which are prime to $p$.
\end{thm}
\proof $(1)$~It is enough to show that the image of the isotropy group $\Z_p^*\times p^{\hat\Z}\subset \Z_p^*\times G_p$ of Lemma \ref{Ggalois}, maps surjectively onto the Galois group $\Gal(\B_\infty(\ell):\B(\ell,p))$. This follows from Lemma \ref{subopen}.

$(2)$~The restriction maps determine an equivariant map
\begin{equation}\label{biject}
{\rm Val}_p(\qsl)\to \prod_{\ell\neq p} {\rm Val}_p(\B_\infty(\ell))
\end{equation}
for the action of the Galois group $\Gal(\qsl:\Q)$. By \eqref{muladdiso} and \eqref{imiso2}, the isotropy groups are the same so that the map \eqref{biject} is bijective.

$(3)$~By extending $v$ to an element of ${\rm Val}_p(\qcy)$ one gets that the restriction to $\qcp$ and hence to $\qcp_\Delta$ is unramified. Moreover the residue field is determined by the topology on the closure set of the action of the  Frobenius, \ie on $\prod_{\ell\neq p}\Z_\ell$. The result follows.\endproof

Next, we shall explain the link with the notations used by de Smit and Lenstra and their construction. First, we recall that the additive group $\Q/\Z$ is the direct sum of its $\ell$-torsion components
\begin{equation}\label{hp}
    H_\ell=\{\alpha\in \Q/\Z\mid \exists n \ \ell^n\alpha=0\}\simeq\Q_\ell/\Z_\ell.
\end{equation}
Thus the group ring $\Z[\Q/\Z]$ can be written as a tensor product
\begin{equation}\label{inftens}
   \Z[\Q/\Z]=\bigotimes_{\ell\,\rm prime} \Z[H_\ell].
\end{equation}
The natural action of $\hat\Z^*$ on $\Z[\Q/\Z]$ by automorphisms of the group $\Q/\Z$ factorizes in the individual actions of $\Z_\ell^*=\Aut(H_\ell)$.

One lets $A_\ell$ be the ring $\Z[X_0,X_1,\ldots]$ modulo the ideal generated by
\begin{equation}\label{rels}
\sum_{j=0}^{\ell-1} X_0^j\,,\ \  X_{k+1}^\ell -X_k\qqq k\geq 0.
\end{equation}
Thus one has $X_0^\ell=1$ in $A_\ell$ and $X_{k+1}^\ell =X_k$ for all $k\geq 0$.  The algebra $B_\ell$ of de Smit and Lenstra is defined as $B_\ell=A_\ell^{\Delta_\ell}$.
The next lemma shows that the algebra $B_\ell$ is intimately related to the fixed point algebra $\Z[H_\ell]^{\Delta_\ell}$.
\begin{lem} \label{smitlenstra} One has
\begin{equation}\label{sl}
    B_\ell\simeq(\Z[H_\ell]/J)^{\Delta_\ell}
\end{equation}
where $J$ is the ideal generated by the relations $\sum_{\ell\gamma=0}e(\gamma)\in \Z[H_\ell]$.
\end{lem}
\proof  It follows from the relations \eqref{rels} that
$X_k^{\ell^{k+1}}=1$ for all $k$. Moreover, the map $   \theta(e(\ell^{-k}))=X_{k-1}$ extends to a surjective homomorphism $\Z[H_\ell]\to A_\ell$ with kernel $J$,  one thus gets \eqref{sl}.\endproof

One has the trace map
\begin{equation}\label{condexp}
    \trr:A_\ell\to B_\ell,\qquad \trr(x)=\sum_{\sigma\in\Delta_\ell}\sigma(x)
\end{equation}
and  natural ring homomorphisms $B_\ell\to E(\ell)$.
De Smit and Lenstra (\cf \cite{smitlenstra}) lift the natural generator of $E_k(\ell)$ as an extension of $E_{k-1}(\ell)$, and the Galois conjugates under $\Gal(E_k(\ell):E_{k-1}(\ell))$ as the following elements of $B_\ell$
\begin{equation}\label{etakl}
    \eta_{\ell,k,i}=\trr(e(\frac{1}{\ell^{k+1}}+\frac i\ell)), \ \ i=0,\ldots \ell-1.
\end{equation}
When $\ell=2$, one has simply ${\Delta_2}=\{\pm 1\}\subset \Z_2^*$ and in this case the above list of elements reduces to
\begin{equation}\label{etaktwo}
    \eta_{2,k}= \trr(e(\frac{1}{2^{k+2}})).
\end{equation}
The two authors show that the prime ideals $\ffp$ of $B_\ell$ which contain $p$, are uniquely specified by a finite system of elements $a(\ffp,j)\in \F_p$, $0\leq j<\ell  u(p,\ell)$. More precisely, $\ffp$ is generated by $p$ and by the  $\eta_{\ell,k+1,i} -a(\ffp,i+k\ell)$ for $0\leq k< u(p,\ell)$ and
$0\leq i<\ell$.

To complete the dictionary with the notations of de Smit and Lenstra, we leave to the reader as an exercise to show that\vspace{.05in}

$\bullet$~The prime ideals $\ffp$ of $B_\ell$ which contain $p$ correspond to the valuations ${\rm Val}_p(\B_\infty(\ell))$.\vspace{.05in}

$\bullet$~The subfield $\Q_p\cap \B_\infty(\ell)\subset \B_\infty(\ell)_v$ is equal to $\B(\ell,p)$.\vspace{.05in}

$\bullet$~The system of elements $a(\ffp,j)\in \F_p$ corresponds, as in Proposition \ref{homfield1}, to the residue of the inclusion $\gamma_v: \B(\ell,p)\to \Q_p$,
defined as in Proposition \ref{homfield}.\vspace{.05in}

Theorem \ref{BCSL} does not yield the full algebraic closure of $\F_p$ but only the subfield
\begin{equation}\label{primetopbis}
    \bigcup_{n \in I(p)} \F_{p^n}\subset \bar\F_p.
\end{equation}
Thus it remains to understand how to produce naturally the missing part
\begin{equation}\label{primetopter}
    \bigcup_{n} \F_{p^{p^n}}\subset \bar\F_p
\end{equation}
in such a way  that the tensor product over $\F_p$ yields $\bar\F_p$.

De Smit and Lenstra construction of $\F_{p^{p^\infty}}=\varinjlim_{n} \F_{p^{p^n}}$ is performed using the following Artin--Schreier equations
\begin{equation}\label{slas}
   y_0^p-y_0=1\,, \ \  y_{n+1}^p-y_{n+1}+\frac{y_n}{y_n+1}=0\qqq n\geq 0
\end{equation}
which have the advantage of simplicity.  E. Witt gave in \cite{witt37} a conceptual construction of $\F_{p^{p^\infty}}$ based on the Witt functor $\winf$ and its finite  truncations
$\W_{p^n}$.  The addition of two Witt vectors $x=(x_j)$ and $y=(y_j)$ is a vector whose components $S_j(x,y)$ were proven by Witt to be polynomials with integer coefficients.
Note also that for $p\neq 2$ the Witt components of $-x$ (the additive inverse of $x$) are simply $-x_j$, but this result does not hold for $p=2$. Recall also that in terms of the Witt vectors, the Frobenius $F$ is given in characteristic $p$, by
$(F(x))_j=x_j^p\qqq j$.

 From \cite{witt37}, one derives the following result

\begin{thm} \label{wittextension} Let $n\in \N$. Let $R_n=\F_p[x_0,x_1,\ldots,x_{n-1}]$ be the ring of polynomials in $n$ variables and  $J_n\subset R_n$ the ideal generated by the components of the Witt vector $F(x)-x-1$, where $x\in \W_{p^{n-1}}(R)$ is the Witt vector with components $x_j$. Then  $J_n$ is a prime ideal and the quotient field of the integral ring $R_n/J_n$ defines the field extension $E_n\simeq\F_{p^{p^n}}$.

As an extension of $E_{n-1}$, $E_n$ is given by an Artin--Schreier equation of the form
\begin{equation}\label{asequ}
    X^p=X+\alpha, \ \ \alpha\in E_{n-1}.
\end{equation}
\end{thm}

One derives, for instance, that the first extensions for $p=2$ are given by the equations with  coefficients in $\F_2$
\begin{eqnarray}
 x_0^2 &=&  1+x_0\nonumber \\
 x_1^2  &=& x_0+x_1 \nonumber \\
 x_2^2  &=&  x_0+x_0^3+x_0 x_1+x_2 \nonumber \\
 x_3^2  &=& x_0+x_0^3+x_0^5+x_0^7+x_0^2 x_1+x_0^3 x_1+x_0^4 x_1+x_0 x_1^3+x_0 x_2+x_0^3 x_2+x_0 x_1 x_2+x_3\nonumber
\end{eqnarray}
 For $p=3$ one gets the following equations  with coefficients in $\F_3$
\begin{eqnarray}
 x_0^3 & = & 1+x_0 \nonumber\\
 x_1^3 & = & 2 x_0+2 x_0^2+x_1\nonumber \\
 x_2^3 & = & 2 x_0+2 x_0^2+2 x_0^4+2 x_0^5+2 x_0^7+2 x_0^8+2 x_0^2 x_1+x_0^3 x_1+2 x_0^4 x_1+x_0 x_1^2+x_0^2 x_1^2+x_2\nonumber
\end{eqnarray}

In this way one obtains a completely canonical construction of the field $\F_{p^{p^\infty}}$ by
simply  writing the equation $F(X)=X+1$ in the ring of Witt vectors $\winf$.

\end{document}